\documentclass[10pt,twoside]{article}

\usepackage{amsfonts,amsmath,amsthm,oldgerm,amscd,amssymb}

\usepackage[all]{xypic}

\usepackage[french]{babel}

\usepackage [latin1]{inputenc}

\usepackage{mathrsfs}

\usepackage[all,cmtip]{xy}

\def\YEAR{\year}\newcount\VOL\VOL=\YEAR\advance\VOL by-1995
\def\firstpage{1}\def\lastpage{1000}
\def\received{}\def\revised{}
\def\communicated{}

\makeatletter
\def\magnification{\afterassignment\m@g\count@}
\def\m@g{\mag=\count@\hsize6.5truein\vsize8.9truein\dimen\footins8truein}
\makeatother

\font\eightrm=cmr8
\font\caps=cmcsc10                    
\font\Caps=cmcsc10 scaled \magstep1   

\makeatletter
\@specialpagefalse\headheight=8.5pt
\def\DocMath{}
\renewcommand{\@evenhead}{%
    \ifnum\thepage>\lastpage\rlap{\thepage}\hfill%
    \else\rlap{\thepage}\slshape\leftmark\hfill{\caps\SAuthor}\hfill\fi}%
\renewcommand{\@oddhead}{%
    \ifnum\thepage=\firstpage{\DocMath\hfill\llap{\thepage}}%
    \else{\slshape\rightmark}\hfill{\caps\STitle}\hfill\llap{\thepage}\fi}%
\makeatother

\def\TSkip{\bigskip}
\newbox\TheTitle{\obeylines\gdef\GetTitle #1
\ShortTitle  #2
\SubTitle    #3
\Author      #4
\ShortAuthor #5
\EndTitle
{\setbox\TheTitle=\vbox{\baselineskip=20pt\let\par=\cr\obeylines%
\halign{\centerline{\Caps##}\cr\noalign{\medskip}\cr#1\cr}}%
        \copy\TheTitle\TSkip\TSkip%
\def\next{#2}\ifx\next\empty\gdef\STitle{#1}\else\gdef\STitle{#2}\fi%
\def\next{#3}\ifx\next\empty%
    \else\setbox\TheTitle=\vbox{\baselineskip=20pt\let\par=\cr\obeylines%
    \halign{\centerline{\caps##} #3\cr}}\copy\TheTitle\TSkip\TSkip\fi%
\centerline{\caps #4}\TSkip\TSkip%
\def\next{#5}\ifx\next\empty\gdef\SAuthor{#4}\else\gdef\SAuthor{#5}\fi%
\ifx\received\empty\relax
    \else\centerline{\eightrm Received: \received}\fi%
\ifx\revised\empty\TSkip%
    \else\centerline{\eightrm Revised: \revised}\TSkip\fi%
\ifx\communicated\empty\relax
    \else\centerline{\eightrm Communicated by \communicated}\fi\TSkip\TSkip%
\catcode'015=5}}\def\Title{\obeylines\GetTitle}
\def\Abstract{\begingroup\narrower
    \parskip=\medskipamount\parindent=0pt{\caps Abstract. }}
\def\EndAbstract{\par\endgroup\TSkip}

\long\def\MSC#1\EndMSC{\def\arg{#1}\ifx\arg\empty\relax\else
     {\par\narrower\noindent%
     2000 Mathematics Subject Classification: #1\par}\fi}

\long\def\KEY#1\EndKEY{\def\arg{#1}\ifx\arg\empty\relax\else
        {\par\narrower\noindent Keywords and Phrases: #1\par}\fi\TSkip}

\newbox\TheAdd\def\Addresses{\vfill\copy\TheAdd\vfill
    \ifodd\number\lastpage\vfill\eject\phantom{.}\vfill\eject\fi}
{\obeylines\gdef\GetAddress #1
\Address #2 
\Address #3
\Address #4
\EndAddress
{\def\xs{4.3truecm}\parindent=0pt
\setbox0=\vtop{{\obeylines\hsize=\xs#1\par}}\def\next{#2}
\ifx\next\empty 
     \setbox\TheAdd=\hbox to\hsize{\hfill\copy0\hfill}
\else\setbox1=\vtop{{\obeylines\hsize=\xs#2\par}}\def\next{#3}
\ifx\next\empty 
     \setbox\TheAdd=\hbox to\hsize{\hfill\copy0\hfill\copy1\hfill}
\else\setbox2=\vtop{{\obeylines\hsize=\xs#3\par}}\def\next{#4}
\ifx\next\empty\ 
     \setbox\TheAdd=\vtop{\hbox to\hsize{\hfill\copy0\hfill\copy1\hfill}
                \vskip20pt\hbox to\hsize{\hfill\copy2\hfill}}
\else\setbox3=\vtop{{\obeylines\hsize=\xs#4\par}}
     \setbox\TheAdd=\vtop{\hbox to\hsize{\hfill\copy0\hfill\copy1\hfill}
                \vskip20pt\hbox to\hsize{\hfill\copy2\hfill\copy3\hfill}}
\fi\fi\fi\catcode'015=5}}\gdef\Address{\obeylines\GetAddress}

\begin{document}

\Title

G-structures entières et modules de Wach
\ShortTitle 
\SubTitle   
\Author Lionel DORAT
\ShortAuthor 
\EndTitle
\Abstract 

In this paper, we study the tannakian properties of the Fontaine-Laffaille functor $\mathop{\bf V_{cris}}$ thanks to the theory of Wach's modules. We construct a point of the torsor linking cristalline representations and weakly admissible filtered modules, preserving the lattices in the sens of the Fontaine-Laffaille correspondance.

\EndAbstract
\MSC 
11F80, 11F85, 11S20, 11S23
\EndMSC
\KEY 
Représentations galoisiennes, représentations cristallines, représentations entières, modules filtrés, $(\phi,\Gamma)$-modules  
\EndKEY
\Address 
IRMA
Université Louis Pasteur
7 rue René Descartes
67084 Strasbourg Cedex
\Address
\Address
\Address
\EndAddress

\renewcommand{\refname}{Bibliographie}

\newcommand{\transposee}[1]{{\vphantom{#1}}^{\mathit t}{#1}}

\def\dsp{\displaystyle}

\def\ptg{\hat{\otimes}}

\def\z{\zeta}
\def\a{\alpha}

\def\r{\rho}
\def\rr{\rho_{{}_{\mathop{\bf MF}}}}
\def\rrr{\rho_{{}_{mr}}}
\def\t{\tau}

\def\b{\beta}

\def\d{\partial}

\def\g{\gamma}

\def\l{\lambda}

\def\s{\sigma}

\def\o{\omega}

\def\p{\varphi}

\def\v{\wedge}

\def\n{\mu}

\def\m{\mathsf{m}}

\def\x{\xi}

\def\w{\mathsf{w}}

\def\f{\mathsf{f}}

\def\I{\mathcal{I}}

\def\U{\mathcal{U}}

\def\P{\mathcal{P}}

\def\W{\mathcal{W}}

\def\N{\mathbb{N}}

\def\V{\mathcal{V}}

\def\C{\mathbb{C}}

\def\A{\mathfrak{A}}

\def\B{\mathcal{B}}

\def\O{\mathcal{O}}

\def\lap{\Delta}

\def\G{\mathbb{G}}

\def\D{{\cal D}}

\def\F{\mathbb{F}}

\def\K{\mathcal{K}}

\def\phi{\varphi}

\def\e{\varepsilon}

\def\R{\mathbb{R}}

\def\S{\mathfrak{S}}

\def\H{\mathcal{H}}

\def\L{\Lambda}

\def\M{\mathcal{M}}

\def\E{{\mathcal E}}

\def\pt{\forall\;}

\def\Z{\mathbb{Z}}

\def\Q{\mathbb{Q}}
\def\Fil{\mathop{\rm Fil}\nolimits}
\def\MF{\mathop{\bf MF_W^{-h}}\nolimits}
\def\V{\mathop{\bf V_{cris} }\nolimits}
\def\Rep{\mathop{\bf  Rep_{\Z_p}(\Gamma_{\K}) }\nolimits}
\def\Gammaphi{\mathop{\bf  \Gamma \Phi M_{\O_{\E}}^{\text{ét}}}\nolimits}
\def\Gam{\mathop{\bf  \Gamma_0 \Phi M_{S_0}^{h}}\nolimits}
\def\Mf{\mathop{\bf MF_W^{\pm h}}\nolimits}
\def\D{\mathop{\bf D_{cris} }\nolimits}
\def\Spec{\mathop{ \rm Spec }\nolimits}
\def\DD{\mathop{\bf \overline D_{cris} }\nolimits}
\def\Ker{\mathop{\rm Ker}\nolimits}
\def\v{\mathop{\bf V_{cris}' }\nolimits}
\def\rg{\mathop{\rm rg }\nolimits}
\def\Tnr{\mathop{\rm \mathcal{T}_{mr}}\nolimits}
\def\T{\mathop{\rm \mathcal{T}_{}}\nolimits}
\def\End{\mathop{ \rm End }\nolimits}
\def\Tan{\mathop{ \bf Tan }\nolimits}
\def\ad{\mathop{ \rm ad }\nolimits}
\def\Hom{\mathop{ \rm Hom }\nolimits}
\def\Ad{\mathop{ \rm Ad }\nolimits}
\def\Id{\mathop{ \rm Id }\nolimits}
\def\Im{\mathop{ \rm Im }\nolimits}

\def\X{\mathcal{X}}

\def\Y{\mathcal{Y}}

\newtheorem{defi}[equation]{D\'efinition}

\newtheorem{theo}[equation]{Th\'eor\`eme}

\newtheorem{prop}[equation]{Proposition}

\newtheorem{lem}[equation]{Lemme}

\newtheorem{cor}[equation]{Corollaire}

\newtheorem{rmq}[equation]{Remarque}
\newtheorem*{rmqq}{Remarque}

\newtheorem{exo}{Exercice}
\newtheorem{Conj}{Conjecture}

\newtheorem{theor}{Th\'eor\`eme}
\newtheorem*{theore}{Th\'eor\`eme 1'}
\newtheorem{coroll}{Corollaire}
\newtheorem*{theorm}{Th\'eor\`eme}

\section*{Introduction}
Dans tout ce travail, $p$ est un nombre premier impair, $\K$ un corps de caractéristique $0$, complet pour une valuation discrète, absolument non ramifié et de corps résiduel $k$ parfait de caractéristique $p$. Nous noterons $W$ l'anneau des vecteurs de Witt à coefficients dans $k$, c'est donc l'anneau des entiers de $\K$. Tous trois sont munis d'une action de Frobenius, notée $\s$. Fixons $\overline{\K}$ une cloture algébrique de $\K$, et posons $\Gamma_{\K}=\mathop{\rm Gal}(\overline{\K},\K)$. Nous noterons $\C$ le complété de $\overline{\K}$ et $\X : \Gamma_{\K}\to \Z_p^*$ désignera le caractère cyclotomique de $\Gamma_{\K}$ (c'est-à-dire que $g(z)=z^{\X(g)}$ pour tout $g\in \Gamma_{\K}$ et pour toute racine de l'unité $z\in \overline{\K}$ d'ordre une puissance de $p$). Nous allons étudier les représentations continues de $\Gamma_{\K}$ dans des $\Q_p$-espaces vectoriels de dimension finie.

Nous nous restreindrons aux représentations cristallines, condition vérifiée dans bien des cas issus de la géométrie (par exemple, pour le module de Tate ou la cohomologie étale à coefficients dans $\Q_p$ d'une variété abélienne ayant bonne réduction). L'avantage de ces représentations est que J.-M. Fontaine et P. Colmez ont montré dans \cite{Font} et \cite{Fofa} qu'elles forment une catégorie tannakienne, qui est $\otimes$-équivalente à la catégorie tannakienne des $\Phi$-modules filtrés sur $\K$ faiblement admissibles (c'est à dire ceux qui ont des réseaux fortement divisibles). 

Le foncteur qui induit cette équivalence de catégories se décrit de la manière suivante : si $V$ est une représentation $p$-adique cristalline, le $\Phi$-module filtré associé est $\mathop{\bf D_{cris,p}}(V)=(V\otimes_{\Q_p} B_{cris})^{\Gamma_{\K}}$ (le quasi-inverse est donné par : pour $D$ un $\Phi$-module filtré faiblement admissible, $\mathop{\bf V_{cris,p}}(D)=\Fil^0(D\otimes_{\K} B_{cris})^{\Phi}$). De plus, l'application 
$$ \mathop{\bf V_{cris,p}}(D) \otimes_{\Q_p} B_{cris}\to D\otimes_{\K} B_{cris}$$
issue de la multiplication de $B_{cris}$ est un isomorphisme (préservant l'action de $\Gamma_{\K}$, la filtration, et le morphisme $\Phi$). Cela peut se traduire de la façon suivante : en notant $w_V$ le foncteur oubli qui à la $\Q_p$-représentation cristalline $V$ associe le $\Q_p$-espace vectoriel sous-jacent à $V$, et $w_D$ celui qui associe le $\K$-espace vectoriel sous-jacent à $ \mathop{\bf D_{cris,p}}(V)$, alors les $\otimes$-isomorphismes du foncteur fibre $w_V\otimes_{\Q_p}\K$ sur le foncteur fibre $w_D$, $\mathop{\bf Isom}(w_V\otimes_{\Q_p}\K,w_D)$, forment un torseur sous $\mathop{\bf Aut^{\otimes}}(w_V)_{|\K}$ et sous $\mathop{\bf Aut^{\otimes}}(w_D)$, qui est non vide sur $B_{cris}$.

Du côté des $\Phi$-modules filtrés sur $\K$, nous disposons de la notion de réseaux fortement divisibles (dont l'existence est une condition nécessaire et suffisante pour que le module soit faiblement admissible), qui sont des $\Phi$-modules filtrés sur $W$ (cf. paragraphe \ref{modu}). J.-M. Fontaine et G. Laffaille ont montré dans \cite{FL} que, si la longueur de la filtration est strictement plus petite que $p-1$, il existe une équivalence de catégories abéliennes entre réseaux fortement divisibles d'un module filtré faiblement admissible, et les réseaux stables de la représentation cristalline associée. 

Plus précisément, à $M$ un $\Phi$-module filtré sur $W$ vérifiant $\Fil^1(M)=\{0\}$ et $\Fil^{2-p}(M)=M$, ils associent le réseau $\V(M)=\Fil^0(M\otimes_{W} A_{cris})^{\phi^0}$, et cette construction induit un foncteur exact, pleinement fidèle (dont nous noterons $\D$ un quasi-inverse). Deux problèmes apparaissent : la condition sur la filtration n'est pas stable par produit tensoriel, et l'application naturelle  
$$\V(M) \otimes_{\Z_p} A_{cris}\to M\otimes_{W} A_{cris}$$
n'est pas un isomorphisme (le déterminant est une puissance de $t$, non inversible dans $A_{cris}$). De plus, une question naturelle se pose : est-ce qu'il existe un point $f$ de $\mathop{\bf Isom}(w_V\otimes_{\Q_p}\K,w_D)$ qui envoie un réseau galoisien sur celui qui lui correpond d'après la correspondance de Fontaine-Laffaille ? Répondre à ces questions revient à étudier les propriétés tannakiennes de $\V$.

L'idée va être d'introduire la théorie des modules de Wach de L. Berger (voir \cite{Berger}), qui à un réseau d'une $\Q_p$-représentation cristalline associe un $(\phi,\Gamma)$-module dont un quotient redonne le $\Phi$-module filtré sur $W$ correspondant à la théorie de Fontaine-Laffaille. Le problème se ramène alors à : pouvons-nous à partir d'un $\Phi$-module filtré sur $W$ reconstruire le module de Wach correspondant ? Pouvons-nous le faire de manière à ce que cette construction soit fonctorielle ?

Le résultat technique principal de cet article est la construction à partir des idées de N. Wach d'un foncteur de la catégorie des modules de Fontaine-Laffaille vers la catégorie des modules de Wach. Plus précisément, notons $\mathop{\bf MF_W^{-h}}$ la catégorie des $\Phi$-modules filtrés $N$ libres sur $W$ tels que $\Fil^{-h}(N)=N$, $\Fil^{1}(N)=\{0\}$ (cf. paragraphe \ref{modu} pour plus de détails) et $\mathop{\bf MF_W<-h>}$ la catégorie engendrée par $\mathop{\bf MF_W^{-h}}$ pour les opérations de sous-objets, objets quotients, produit tensoriel et somme directe, $\V$ le foncteur de Fontaine-Laffaille, $\mathop{\bf \Gamma\Phi M_S^{-h}}$ la catégories des duaux des modules de Wach de hauteur $h$ (ce qui correspond à des modules de Wach d'après la définition de \cite{Berger}) et $\mathop{\bf \Gamma\Phi M_S^-}$ la réunion des $\mathop{\bf \Gamma\Phi M_S^{-h}}$, $\mathop{\bf N}$ le foncteur ``module de Wach'', $\mathop{\bf V_{\O_{\E}}}$ le foncteur de Fontaine pour les $(\phi,\Gamma)$-module sur $\O_{\E}$, $\mathop{\bf V_{cris,p}}$ le foncteur de Fontaine pour les $\Phi$-modules filtrés sur $\K$ admissibles, et $j:S\to \O_{\E}$ qui induit le foncteur extension des scalaires $j^*$ de la catégorie des modules de Wach vers la catégorie des $(\phi,\Gamma)$-modules sur $\O_{\E}$. 
\begin{theor}
Soit $h$ un entier compris entre $0$ et $p-2$, alors il existe un foncteur $\mathop{\rm F^-}$ exact, préservant le produit tensoriel, fidèle et pleinement fidèle de $\mathop{\bf MF_W<-h>}$ vers $\mathop{\bf \Gamma\Phi M_S^-}$. Restreint à $\mathop{\bf MF_W^{-h}}$, ce foncteur est essentiellement surjectif sur $\mathop{\bf \Gamma\Phi M_S^{-h}}$. De plus, pour tout objet $M$ de $\mathop{\bf MF_W^{-h}}$, $\mathop{\rm F^-}(M)$ est fonctoriellement isomorphe à $\mathop{\bf N}(\V(M))$. Dans le cas général, $\mathop{\rm F^-}(M)$ s'interprète encore comme le module de Wach du réseau galoisien correspondant au $(\phi,\gamma)$-module sur $\O_{\E}$ engendré par $\mathop{\rm F^-}(M)$. En outre, $\mathop{\bf V_{\O_{\E}}}\circ j^*\circ \mathop{\rm F^-}$ est isomorphe (comme foncteur) à $\mathop{\bf V_{cris,p}}$ une fois $p$ rendu inversible, et à $\V$ une fois restreint à la catégorie $\mathop{\bf MF_W^{-h}}$.
\end{theor}

\begin{rmq}
Ce théorème est optimal, dans le sens où nous ne pouvons espérer que $\mathop{\rm F^-}$ soit essentiellement surjectif sans la restriction sur $h$.
\end{rmq}
Pour illuster, le théorème nous dit essentiellement que le diagramme suivant est commutatif (où bien sûr il faut resteindre la catégorie des réseaux des représentations cristallines à ceux à poids de Hodge-Tate dans $[\![0,h]\!]$) :


\[\xymatrix{
\mathop{\bf Rep_{\Z_p}^{cris,h}}(\Gamma_{\K}) \ar@<1ex>[rrr]^-{\mathop{\bf D_{cris}}}\ar@<1ex>[d]^-{\mathop{\bf D_{\O_{\E}}}}\ar[drrr]_-{\mathop{\bf N}}&&&\mathop{\bf MF_W^{-h}} \ar@{.>}[d]_-{\mathop{\rm F^-}}\ar[lll]^-{\V}\\
\mathop{\bf \Gamma\Phi M_{\O_{\E}}^{et}}\ar[u]^-{\mathop{\bf V_{\O_{\E}}}}&&&\mathop{\bf \Gamma\Phi M_S^{-}}\ar@<1ex>[lll]^-{j^*}\ar@<-1ex>[u]_-{\mathop{\rm mod} \pi}
}\]

et, une fois $p$ rendu inversible, 

\[\xymatrix{
\mathop{\bf Rep_{\Q_p}^{cris}}(\Gamma_{\K}) \ar@<1ex>[d]^-{\mathop{\bf D_{\E}}}\ar[drrr]_-{\mathop{\bf N}}&&&\mathop{\bf MF_W<-h>}\otimes\K \ar@{.>}[d]_-{\mathop{\rm F^-}}\ar[lll]^-{\mathop{\bf V_{cris,p}}}\\
\mathop{\bf \Gamma\Phi M_{\E}^{et}}\ar[u]^-{\mathop{\bf V_{\E}}}&&&\mathop{\bf \Gamma\Phi M_{S[\frac{1}{p}]}^{-}}\ar@<1ex>[lll]^-{j^*}\ar@<-1ex>[u]_-{\mathop{\rm mod} \pi}
}\]
où $\mathop{\rm F^-}$ est fidèle, pleinement fidèle, préserve le produit tensoriel, et suivant les cas, peut être essentiellement surjectif (et $\mathop{\bf MF_W<-h>}\otimes\K $ représente juste la catégorie formée des objets de $\mathop{\bf MF_W<-h>}$ où nous avons rendu $p$ inversible, c'est à dire la catégorie engendrée pour les opérations de produit tensoriel, somme directe, sous-objet et objet quotient, par les modules filtrés sur $\K$ admissibles à pente compris entre $0$ et $-h$).

De ce théorème, nous en déduisons le corollaire voulu : 
\begin{theor}\label{theotorseur}
Il existe un point du torseur $\mathop{\bf Isom}(w_V\otimes_{\Q_p}\K,w_D)$ à coefficient dans le corps $\widehat{\E}_{nr}$ qui préserve les réseaux de Fontaine-Laffaille, c'est à dire qui identifie les réseaux stables par Galois des représentations cristallines à poids de Hodge-Tate dans $[\![0,\frac{p-2}{2}]\!]$ au $W$-module filtré correspondant par la théorie de Fontaine-Laffaille. 
\end{theor}

Pour obtenir un résultat sur $\K$ plutôt que sur $\widehat{\E}_{nr}$, il faut modifier le problème. Considérons $G$ un groupe algébrique lisse sur $\Z_p$ et une représentation $\r :\Gamma_{\K}\to G(\Z_p)$. Supposons donnée une immersion fermée $\a$ de $G$ dans $GL_U$, pour $U$ un $\Z_p$-module libre de rang fini, telle que la représentation $\a\circ \r$ de $\Gamma_{\K}$ (dans $GL(U\otimes_{\Z_p}\Q_p)$) soit cristalline à poids de Hodge-Tate dans $[\![0,h]\!]$ avec $h$ un entier compris entre $0$ et $\frac{p-2}{2}$. Notons $V=U\otimes_{\Z_p}\Q_p$. Par un théorème de Chevalley, il existe un $\Q_p$-espace vectoriel $V_{G}$ dans $\bigoplus_i \End(V)^{\otimes i}$ (en faisant agir $GL_V$ naturellement sur $V^*$ et trivialement sur $V$, dans $\End(V)=V\otimes V^*$) tel que $G\times_{\Z_p} \Q_p$ soit le groupe algébrique formé de l'ensemble des éléments de $GL_{V}$ qui laissent stable $V_{G}$. Alors, par le foncteur de Fontaine-Laffaille, nous pouvons définir naturellement un groupe $G_D$ sur $D=\mathop{\bf D_{cris,p}}(V)$ comme l'ensemble des éléments de $GL_{D}$ laissant stable $\mathop{\bf D_{cris,p}}(V_{G})$. Un corollaire de la proposition 6.3.3 de \cite{Fo} nous donne l'existence d'un élément de $\mathop{\bf Isom}(w_V\otimes_{\Q_p}\K,w_D)(\K)$, donc en particulier d'un isomorphisme de $\K$-modules $$f:V\otimes_{\Q_p}\K\to D$$ qui identifie $G\times_{\Z_p} \K$ à $G_D$. Le comportement de $f$ vis-à-vis des réseaux est à priori inconnu. Pour l'étudier, nous introduisons un $G$-torseur $\mathop{\bf Isom} $ défini sur $W$, qui est heuristiquement le $G$-torseur obtenu à partir de $\mathop{\bf Isom}(w_V\otimes_{\Q_p}\K,w_D)$ (c'est à dire une forme sur $W$ du $G\times_W \K$ torseur obtenu à partir de $\mathop{\bf Isom}(w_V\otimes_{\Q_p}\K,w_D)$). Le résultat suivant se montre alors en montrant que $\mathop{\bf Isom} $ est un $G$-torseur trivial sur $W$ :

\begin{theor}\label{thth}
Sous les hypothèses précédentes, si $M=\D(U)$, il existe un sous-groupe algébrique $G_M$ de $GL_M$ sur $W$, avec $G_M\times_W \K=G_D$, et il existe $f$ un isomorphisme de $W$-modules $$f:U\otimes_{\Z_p}W\to M$$ qui identifie $G$ à $G_M$.

De plus, si $U'$ est un réseau de $U\otimes_{\Z_p}\Q_p$ laissé stable par l'action de $G$, alors $f[\frac{1}{p}]$ envoie $U'\otimes_{\Z_p}W$ sur $\D(U')$.
\end{theor}

\begin{rmq}
Ce théorème nous donne en particulier que les réseaux $U$ et $M$ ont la m\^eme position vis à vis du groupe $G$.
\end{rmq}

Avec des hypothèses plus fortes sur $\a$, nous pouvons affaiblir l'hypothèse sur $h$. Une application directe de ce résultat concerne la semi-simplifiée d'une représentation cristalline à poids de Hodge-Tate petits : le groupe algébrique $H$ engendré par l'image de Galois sur $\Q_p$ est alors connexe et réductif, donc en appliquant les résultats cités dans \cite{T} (paragraphe 3.2 et 3.4.1), il existe un groupe algébrique lisse $G$ défini sur $\Z_p$, tel que $G(\Z_p)$ contienne l'image de Galois, et dont la fibre générique est $H$. Le Théorème \ref{thth} s'applique alors.

\section{Rappels}
\subsection{Rappels sur les $(\phi,\Gamma)$-modules}\label{rappels}
\subsubsection{Définition de $\O_{\E}$}\label{gamma}
Soit $R$\index{$R$} l'ensemble des suites $x=(x^{(n)})_{n\in \N}$ formées d'éléments de $\O_{\bar \K}/p\O_{\bar \K}$ vérifiant $(x^{(n+1)})^p=x^{(n)}$ pour tout $n$ (cf. \cite{F82}, p. 535). C'est un anneau parfait de caractéristique $p$, muni d'une valuation; son corps résiduel s'identifie à $k$. Son corps des fractions $\mathop{\rm Fr}R$ est un corps algébriquement clos de caractéristique $p$, et $R$ est intégralement clos dans $\mathop{\rm Fr}R$.

Si $A$ est une $k$-algèbre, $W(A)$ désigne l'anneau des vecteurs de Witt à coefficients dans $A$. Notons $\Z_{p^s}=W(\F_{p^s})$\index{$\Z_{p^s}$}, $\Z_p^{nr}=W(\overline{\F_p})$\index{$\Z_p^{nr}$}, $W\index{$W$}=W(k)$, $W_{\K}(A)=\K\otimes_W W(A)=W(A)[\frac{1}{p}]$ et si $a\in A$, $[a]=(a,0,\cdots,0,\cdots)$ le représentant de Teichmüller de $a$ dans $W(A)$. Le Frobenius $x\in A\mapsto x^p\in A$ s'étend à $W(A)$ en $\p$ (encore appelé l'endomorphisme de Frobenius) par fonctorialité, ainsi qu'à $W_{\K}(A)$; nous noterons $\s$ le Frobenius sur $W$ et sur $\K$ (si $\l\in W$, $\s(\l):=\p(\l)$). En particulier ceci s'applique à $W(R)$, $W(\mathop{\rm Fr}R)$ et $W_{\K}(\mathop{\rm Fr}R)$.

D'autre part, le groupe $\Gamma_{\K}$ opère par fonctorialité sur $R$, $\mathop{\rm Fr}R$ et $W(\mathop{\rm Fr}(R))$, et les anneaux $W(R)$, $W(\mathop{\rm Fr}R)$ et $W_{\K}(R)$ s'identifient à des sous-anneaux de $W_{\K}(\mathop{\rm Fr}R)$ stables par $\p$ et $\Gamma_{\K}$.

Notons $\Z_p(1)=\displaystyle \varprojlim_{n\in\N}\mu_{p^n}(\bar \K)$ le module de Tate du groupe multiplicatif et pour tout $i\in \N$, $\Z_p(i)=\Z_p(1)^{\otimes i}$ et $\Z_p(-i)$ son dual. Pour tout $\Z_p$-module $T$, et pour tout $i\in\Z$, posons $T(i)=T\otimes_{\Z_p}\Z_p(i)$.

Le module de Tate $\Z_p(1)=T_p(\G_m)$ s'identifie au sous $\Z_p$-module du groupe multiplicatif des unités de $R$ congrues à $1$ modulo l'idéal maximal, formé des $x$ tels que $x^{(0)}=1$. Choisisons un générateur de ce module, c'est-à-dire un élément $\e=(\e^{(n)})_{n\in\N}\in R$ tel que $\e^{(0)}=1$ et $\e^{(1)}\neq 1$, et considérons l'élément $\pi=[\e]-1$ dans $W(R)$. Alors l'adhérence $S\index{$S$}$ de la sous $W$-algèbre de $W(R)$ engendrée par $\pi$ s'identifie à l'algèbre $W[\![\pi]\!]$ des séries formelles en $\pi$ à coefficients dans $W$; de plus $S$ est stable par $\phi$ et $\Gamma_{\K}$, et nous avons les relations suivantes :

$$\phi(\pi)=(1+\pi)^p-1$$
$$ g(\pi)=(1+\pi)^{\X(g)}-1$$ pour $g\in \Gamma_{\K}$.

Soit $\K_n$ le sous corps de $\bar\K$ engendré sur $\K$ par les racines $p^n$-ièmes de l'unité, et $\K_{\infty}=\cup_{n\in \N}\K_n$. Notons $\Gamma\index{$\Gamma$}=\mathop{\rm Gal}(\K_{\infty}/\K)$ et $H_{\K}$ le noyau de la projection de $\Gamma_{\K}$ sur $\Gamma$. Le groupe $H_{\K}$ agit trivialement sur $S$. Si $\Gamma_f$ est le sous-groupe de torsion de $\Gamma$, posons $S_0\index{$S_0$}=S^{\Gamma_f}$ ainsi que $\Gamma_0\index{$\Gamma_0$}=\Gamma/\Gamma_f$; J.-M. Fontaine a montré (cf. \cite{F91}, p. 268-273) que $S_0=W[\![\pi_0]\!]$, où $\pi_0=-p+\sum\limits_{a\in\F_p}[\e]^{[a]}$. Notons $q=p+\pi_0$\index{$q$}. $S_0$ est munie d'une action naturelle de $\Gamma_0$.

Notons $\O_{\E}\index{$\O_{\E}$}$ le complété pour la topologie $p$-adique de $S[\frac{1}{\pi}]$. C'est l'anneau des entiers d'un corps complet pour une valuation discrète, absolument non ramifié, noté $\E$. Comme $\pi$ est inversible dans $W(\mathop{\rm Fr}R)$, l'inclusion de $S$ dans $W(R)$ se prolonge en un plongement de $S[\frac{1}{\pi}]$ dans $W(\mathop{\rm Fr}R)$, et $\O_{\E}$ s'identifie à l'adhérence de $S[\frac{1}{\pi}]$ dans $W(\mathop{\rm Fr}R)$ pour la topologie $p$-adique, tandis que $\E=\O_{\E}[\frac{1}{p}]$ s'identfie à un sous-corps de $W_{\K}(\mathop{\rm Fr}R)$. Alors si $E=\O_{\E}/p$, $\O_E=S/pS=k[\![\tilde{\pi}]\!]$, où $\tilde{\pi}$ est la réduction modulo $p$ de $\pi$.

De plus, si $\hat{\E}_{nr}$ désigne l'adhérence dans $W_{\K}(\mathop{\rm Fr}R)$ de l'extension maximale non ramifiée $\E_{nr}$ de $\E$ contenue dans $W_{\K}(\mathop{\rm Fr}R)$ et $\O_{\hat{\E}_{nr}}\index{$\O_{\hat{\E}_{nr}}$}$ son anneau des entiers, $\O_{\hat{\E}_{nr}}/p$ est une clôture séparable $E^{sep}$ de $E$, avec une identification des groupes de Galois $$H_{\K}=\mathop{\rm Gal}(E^{sep}/E)=\mathop{\rm Gal}(\E_{nr}/\E).$$

\subsubsection{$(\p,\Gamma)$-modules et représentations galoisiennes}\label{isomfontaine}
Nous ne considérerons des $(\p,\Gamma)$-modules que sur $S$ ou $\O_{\E}$ (nous considérerons aussi des $(\p, \Gamma_0)$-modules définis sur $S_0$). Soit $A$ l'un des anneaux précédent. Un $(\phi, \Gamma)$-module sur $A$ est un $A$-module muni d'un endomorphisme $\p$, semi-linéaire par rapport à $\sigma$ muni en plus d'une action continue de $\Gamma$, semi-linéaire par rapport à l'action de $\Gamma$ sur $A$, cette action commutant avec l'endomorphisme $\p$. Nous les supposerons toujours {\it étale}, c'est à dire de type fini sur $A$ et tels que l'application linéaire $\Phi : M^{\sigma}\to M$, déduite de $\p$ en posant $\Phi(\l\otimes x)=\l\p(x)$ pour $\l\in A$ et $x\in M$ est bijective. Les $(\p,\Gamma)$-modules étales (avec comme morphismes les morphismes $A$-linéaires commutants à $\p$ et à $\Gamma$) définissent une $\otimes$-catégorie abélienne notée $\mathop{\bf  \Gamma\Phi M_{A}^{\text{ét}}}$ (cf. \cite{F91} p.273).

Appelons représentation $\Z_p$-adique de $\Gamma_{\K}$ la donnée d'un $\Z_p$-module de type fini muni d'une action linéaire et continue de $\Gamma_{\K}$. Un morphisme sera une application $\Z_p$-linéaire commutant à l'action de $\Gamma_{\K}$. Notons $\mathop{\bf  Rep_{\Z_p}(\Gamma_{\K})}$ la catégorie des représentations $\Z_p$-adique de $\Gamma_{\K}$. La catégorie $\mathop{\bf  Rep_{\Q_p}(\Gamma_{\K})}$ est défini de même.

J.-M. Fontaine a montré dans \cite{F91} (p. 274) qu'il existait une équivalence de catégories entre $\mathop{\bf  \Gamma\Phi M_{\O_{\E}}^{\text{ét}}}$ et $\mathop{\bf  Rep_{\Z_p}(\Gamma_{\K})}$ induite par le foncteur $\mathop{\bf D_{\O_{\E}}}(T)=(\O_{\widehat{\E}_{nr}}\otimes_{\Z_p}T)^{H_{\K}}$ pour $T$ une $\Z_p$-représentation de $\Gamma_{\K}$, et son quasi invers $\mathop{\bf V_{\O_{\E}}}(\mathcal N)=(\O_{\widehat{\E}_{nr}}\otimes_{\O_{\E}}\mathcal N)^{\p=1}$. La multiplication dans $\O_{\widehat{\E}_{nr}}$ induit alors une application naturelle et fonctorielle :
\[\xymatrix{  {\mathop{\bf V_{\O_{\E}}}(\mathcal{N})\otimes_{\Z_p}\O_{\widehat{\E}_{nr}}} \ar[r]^(0.55){\psi_{\mathcal{N}}}& {\mathcal{N}\otimes_{\O_{\E}} \O_{\widehat{\E}_{nr}}}
}\] 
pour $\mathcal{N}$ un objet de la catégorie $\Gammaphi$.

\subsection{Représentations cristallines}

\subsubsection{Représentations cristallines}\label{cristalline}
Pour la définition de $A_{cris}$ et de $t:=\mathop{\rm log}([\e])$, nous renvoyons à \cite{Fo2} par exemple. Nous noterons $B_{cris}=A_{cris}[\frac{1}{t}]$. Soit $V$ un $\Q_p$-espace vectoriel de dimension finie, et $\r : \Gamma_{\K}\to GL(V)$ une représentation continue de $\Gamma_{\K}$. Définissons $\mathop{\bf D_{cris,p}}$\index{$\mathop{\bf D_{cris,p}}$} par  $$\mathop{\bf D_{cris,p}}(V)=(V\otimes_{\Q_p} B_{cris})^{\Gamma_{\K}}$$
Alors $\mathop{\bf D_{cris,p}}(V)$ est un $\K$-espace vectoriel, et $\dim_{\K} \mathop{\bf D_{cris,p}}(V)\leq \dim_{\Q_p}V$. 

\begin{defi}
La représentation $(\r,V)$ est {\it cristalline}\index{cristalline} si $\dim_{\K} \mathop{\bf D_{cris,p}}(V)= \dim_{\Q_p}V$.
\end{defi}
Notons $\mathop{\bf  Rep_{\Q_p,cris}(\Gamma_{\K})}$ la sous-catégorie pleine de $\mathop{\bf  Rep_{\Q_p}(\Gamma_{\K})}$ formée par les représentations cristallines. Définissons $\mathop{\bf MF_{\K}}$ la catégorie des $\Phi$-modules filtrés sur $\K$ : un objet $D$ de $\mathop{\bf MF_{\K}}$ est un $\K$-espace vectoriel de dimension finie muni d'une filtration $(\mathop{\rm Fil}^i(D))_{i\in \Z}$ formée de sous-espaces vectoriels, fitration qui est décroissante, exhaustive séparée, et muni d'une application $\sigma$-semi-linéaire bijective $\Phi : D\to D$. $\mathop{\bf D_{cris,p}}(V)$ est alors naturellement un $\Phi$-module filtré. Un élément de l''image essentiel du foncteur $\mathop{\bf D_{cris,p}}(V)$ restreint à $\mathop{\bf  Rep_{\Q_p,cris}(\Gamma_{\K})}$ est appelé admissible. Notons $\mathop{\bf MF_{\K}^{ad}}$ la sous-catégorie pleine de $\mathop{\bf MF_{\K}}$ formée des modules admissibles.

$\mathop{\bf  Rep_{\Q_p}(\Gamma_{\K})}$ et $\mathop{\bf MF_{\K}^{ad}}$ sont deux catégories tannakiennes, le foncteur $\mathop{\bf D_{cris,p}}$ induit une équivalence de $\otimes$-catégories entre ces deux catégories, et un quasi-inverse est donné par le foncteur $\mathop{\bf V_{cris,p}}(D)=\begin{array}{c}\mathop{\rm Fil}^0(D\otimes_{\K}B_{cris})^{\Phi=1}\end{array}$. L'application naturelle (provenant de la multiplication dans $B_{cris}$) 
\begin{equation}\label{isomfont}\mathop{\bf V_{cris,p}}(D)\otimes_{\Q_p}B_{cris} \to D\otimes_{\K}B_{cris} \end{equation} est alors une bijection.

\subsubsection{Poids de Hodge-Tate}\label{poids HT}
Rappelons que pour $\r:\,\Gamma_{\K}\to GL_{\Q_p}(V)$ une représentation continue sur un $\Q_p$-espace vectoriel de dimension finie, l'action de $\Gamma_{\K}$ peut s'étendre à $V_{\C}=V\otimes_{\Q_p}\C$ via $g(v\otimes x)=\r(g)(v)\otimes g(x)$. Notons alors pour $i\in\Z$, $V_{\C}\{i\}=\{v\in V_{\C} |\forall g\in \Gamma_{\K},\,g(v)=\X(g)^iv\}$. $V_{\C}\{i\}$ est un $\K$-sous espace vectoriel de $V_{\C}$ tel que l'injection $V_{\C}\{i\}\to V_{\C}$ s'étend en une injection $\C$-linéaire $$\bigoplus_{i\in\Z}V_{\C}\{i\}\otimes_{\K}\C\to V_{\C}$$
Alors $V$ est dit de Hodge-Tate si cette injection est une bijection. Les poids de Hodge-Tate sont alors les $i\in\Z$ tels que $\mathop{\rm dim}_{\K}V_{\C}\{i\}\neq 0$. Si $V$ est cristalline, alors elle est de Hodge-Tate, et ses poids de Hodge-Tate sont les opposées des sauts de la filtration de $\mathop{\bf D_{cris,p}}(V)$.

\subsection{Rappels sur les $\Phi$-modules}\label{modu}
La catégorie qui va nous interesser est la catégorie $\mathop{\bf MF_{W,tf}}$ dite des $\Phi$-modules filtrés sur $W$, dont les objets sont les $W$-modules $N$ de type fini, muni
\begin{itemize}
\item d'une filtration décroissante exhaustive et séparée formée de sous-modules $(\Fil^i (N))_{i\in \Z}$ ;
\item pour tout $i\in\Z$, d'une application $\s$-semi-linéaire $\phi^i \,:\,\Fil^i(N)\to N$ telle que $\phi^i|_{\Fil^{i+1}(N)}=p\phi^{i+1}$ ; 
\item il existe $i\in \Z$ avec $\Fil^i(N)=\{0\}$ ;
\item les $\Fil^i(N)$ sont des facteurs directs dans $N$ ;
\item $\sum\limits_{i\in \Z} \phi^i(\Fil^i(N))=N$.
\end{itemize}
Les morphismes de cette catégorie sont donnés par les applications $W$-linéaires compatibles aux filtrations et commutants aux $\p^i$. C'est une $\otimes$-catégorie qui est abélienne, $\Z_p$-linéaire, qui possède des $\mathop{\rm Hom}$ internes (cf. \cite{Wi}).

Soit $X$\index{$X$} (respectivement $X_s$ pour $s\in \N^*$) le groupe additif des applications périodiques (respectivement ayant $s$ pour période) de $\Z$ dans $\Z$. Le Frobenius $\s$ agit sur $X$ par $\forall \x\in X,\forall i\in \Z,\,\s(\x)(i)=\x(i+1)$, et laisse donc stable les $X_s$. 

Pour tout objet $N$ de $\mathop{\bf MF_{W,tf}}$, si $(N_i)_{i\in \Z}$ est un scindage de $(\Fil^i(N))_{i\in \Z}$, posons pour $x\in N$ tel que $x=\sum\limits_i x_i$ avec $x_i\in N_i$, $f_N(x)=\sum\limits_i\phi_N^i(x_i)$. Soit pour tout $\x\in X$, le $W$-module $N\{\x\}:=\{x\in N|f_N^j(x)\in N_{\x(j)}$ pour tout $j\in \Z\}$. Le module $N$ est dit {\it élémentaire} si $N=\oplus_{\x\in X} N\{\x\}$. 
\begin{lem}\label{lemmebaseadaptee}
Si $N$ est un module élémentaire, dont le module sous-jacent est libre sur $W$ ou sur $k$, alors il existe une base $(e_{\x}^i)_{\x\in X,\,1\leq i\leq \rg(N_{\x})}$ telle que $\phi^{\x(0)}(e_{\x}^i)=e_{\s(x)}^i$.
\end{lem}

J.-P. Wintenberger a montré dans \cite{Wi} :
\begin{theo}\label{thwi}
Pour tout objet $N$ de $\mathop{\bf MF_{W,tf}}$, il existe un et un seul scindage de la filtration de $N$ tel que 
\begin{itemize}
\item il existe un (unique) $u_N\in Aut_W(N)$ tel que $(N,(N_i),u_N^{-1}\circ f_N)$ soit élémentaire ;
\item $N/pN$ ait une suite de composition dont les quotients successifs sont des modules élémentaires.
\end{itemize}
Ce scindage vérifie les propriétés de fonctorialité attendues.
\end{theo}


Posons enfin $\mathop{\bf MF_{W,tf}^{[a,b]}}$ (resp. $\mathop{\bf MF_W^{[a,b]}} $) la sous-catégorie pleine de $\mathop{\bf MF_{W,tf}}$ formée des $W$-modules $M$ (resp. modules libres) tels que $\Fil^{a}(M)=M$ et $\Fil^{b+1}(M)=\{0\}$. Notons $\mathop{\bf MF_{W,tf}^{-h}}=\mathop{\bf MF_{W,tf}^{[-h,0]}}$, $\mathop{\bf MF_{W,tf}^{h}}=\mathop{\bf MF_{W,tf}^{[0,h]}}$ et $\mathop{\bf MF_{W,tf}^{\pm h}}=\mathop{\bf MF_{W,tf}^{[-h,h]}}$ (de même sans le symbole $\mathop{\rm tf}$). Pour terminer, nous désignerons par $\mathop{\bf MF_{W,tf}<h>}$ la catégorie engendrée par  $\mathop{\bf MF_W^{ h}}$ dans la catégorie $\mathop{\bf MF_{W,tf}}$ pour les opérations de sous-objet, objet quotient, somme directe et produit tensoriel.

Soit $D$ un $\Phi$-module filtré sur $\K$ admissible. Alors il possède des sous-réseaux fortement divisible, $M$, c'est-à-dire un réseau $M$ vérifiant $\sum\limits_{i\in \Z} p^{-i}\Phi(\Fil^i(D)\cap M)=M$. En posant $\Fil^i(M)=\Fil^i(D)\cap M$, $\phi^i=p^{-i}\Phi|_{\Fil^i(M)}$, $M$ devient un $\Phi$-module filtré sur $W$. Réciproquement, si $M$ est un objet de $\mathop{\bf MF_{W,tf}}$ libre sur $W$, en posant $D:=\K\otimes_W M$, $\Fil^i(D):=\K\otimes_W \Fil^i(M)$, et pour $x_i\in \Fil^i(M),\,\Phi(x_i):=p^i\phi^i(x_i)$, l'objet $D$ ainsi construit est un $\Phi$-module filtré sur $\K$ faiblement admissible (et donc en fait admissible) dont $M$ est un réseau fortement divisible. Par contre, différents $M$ peuvent donner le même $D$. Nous noterons $D_M$ ce $\Phi$-module filtré sur $\K$ faiblement admissible construit à partir de $M$.

\subsubsection{Le théorème de Fontaine-Laffaille}
\begin{defi}
Pour tout objet $M$ de $\mathop{\bf MF_{W,tf}}$ tel que $\Fil^1(M)=\{0\}$, soit $\V(M)$\index{$\V$} la représentation galoisienne définie par :
$$\V(M)=\Fil^0(M\otimes_W A_{cris})^{\phi^0}$$
Si $M$ est libre comme $W$-module, $\V(M)$ est un $\Z_p$-module libre (c'est un sous-réseau de $\mathop{\bf V_{cris,p}}(D_M)$).\end{defi}

\begin{theo}[Théorème de Fontaine-Laffaille]\label{theofontainelaffaille}
Si nous nous restreignons à la sous-catégorie pleine des $M$ vérifiant $\Fil^{2-p}(M)=M$ et $\Fil^1(M)=\{0\}$, alors le foncteur $\V$ ainsi défini est exact et pleinement fidèle. De plus si $M$ est libre sur $W$, $\V(M)$ est un réseau de la représentation galoisienne associée à $D_M$ (c'est-à-dire que $\rg_{\Z_p}(\V(M))=\rg_W(M)$).
\end{theo}
 
Nous noterons $\mathop{\bf D_{cris} }$ un quasi-inverse à $\V$.

\section{Construction du foncteur}\label{fonc}

\subsubsection{Rappels sur $\Gam$}
Notons $\Gam$ ($\mathop{\bf \Gamma\Phi M_S^h}$ se définit de la même façon) la sous-catégorie pleine de la catégorie des $(\phi,\Gamma_0)$-modules sur $S_0$ (cf. paragraphe \ref{rappels}) formée des objets $\mathcal{N}$ vérifiant :
\begin{itemize}
\item le $S_0$-module sous-jacent est de type fini et sans $p'$-torsion (i.e. pour tout élément irréductible $\l$ de $S_0$ non associé à $p$, $\mathcal{N}$ est sans $\l$-torsion),
\item le $S_0$-module $\mathcal{N}/\Phi(\mathcal{N}\otimes_{\s}S_0)$ est annulé par $q^h$ (où $q=\pi_0+p$),
\item le groupe $\Gamma_0$ agit trivialement sur $\mathcal{N}/\pi_0\mathcal{N}$.
\end{itemize}
Elle est abélienne si $0\leq h\leq p-2$, et l'inclusion $j:S_0\to \O_{\E}$ induit un foncteur $j^* : \Gam \to \Gammaphi$\index{$j^*$} pleinement fidèle qui est une équivalence de catégorie pour $0\leq h\leq p-2$ sur son image essentielle (cf. \cite{F91}, p.301). Si $\mathcal{N}$ est un objet de $\Gam$, alors $j^*(\mathcal{N})$ a pour espace sous-jacent $\mathcal{N}\otimes_{S_0}\O_{\E}$. Nous ferons souvent l'abus de notation de n'écrire que l'espace sous-jacent pour désigner $j^*(\mathcal{N})$.

Si $0\leq h\leq p-2$ et $\mathcal{N}$ un objet de $\Gam$, N. Wach a montré qu'il est possible de munir $N=\mathcal{N}/\pi_0N$ d'une structure de $\Phi$-module filtré sur $W$ en posant
$$\Fil^r N=\{x\in N \text{ tels qu'il existe un relèvement }\widehat{x}\in\mathcal{N} \text{ de }x \text{ avec }\p(\widehat{x})\in q^r\mathcal{N}\}$$
et pour tout $x\in  \Fil^r N$, $\p^r(x)$ égal à l'image de $\frac{\p(\widehat{x})}{q^r}$ dans $N$. Elle a alors démontré le théorème suivant (cf. \cite{Wa}, p.231) : 
\begin{theo}
Soit $0\leq h\leq p-2$. Pour tout objet $\mathcal{N}$ de $\Gam$, le $\Phi$-module filtré $i^*(\mathcal{N})=\mathcal{N}/\pi_0\mathcal{N}$ est un objet de $\mathop{\bf MF_{W,tf}^h}$; le foncteur $i^*$ ainsi défini est exact et fidèle.
\end{theo}

\subsubsection{Foncteur entre $\mathop{\bf MF_W^{h}}$ et $\Gam$}
N. Wach a donné les idées pour construire un quasi-inverse à $i^*$ : à partir d'un objet $N$ de $\mathop{\bf MF_{W}^h}$ avec $0\leq h\leq p-2$ et d'une base adaptée à la filtration, elle a construit un objet  $\mathcal{N}$ tel que $i^*(\mathcal{N})=N$. Nous allons montrer qu'en se fixant un scindage fonctoriel de la filtration, nous rendons cette construction fonctorielle.

\begin{prop}\label{propositionbesoin1}
Soit $\mathop{\bf MF_{W,tf}^+}$ la sous-catégorie pleine formée de la réunion des $\mathop{\bf MF_{W,tf}^h}$ (même définition pour $\mathop{\bf  \Gamma_0 \Phi M_{S_0}^{+}}$). A tout scindage fonctoriel de la filtration des objets de $\mathop{\bf MF_{W,tf}^+}$ nous pouvons associer un foncteur de $\mathop{\bf MF_{W,tf}^+}$ vers $\mathop{\bf  \Phi M_{S_0}^{\text{ét}}}$ (la catégorie des $\phi$-modules sur $S_0$ dont l'extension à $\O_{\E}$ donne un $\phi$-module étale), qui soit fidèle, additif, exact, et qui préserve le produit tensoriel.
\end{prop}

\begin{proof}
Si $N$ est un objet de $\mathop{\bf MF_{W,tf}^+}$, et $N=\oplus N_i$ le scindage de la filtration, il suffit de construire sur $N\otimes_W S_0$ une structure de $\phi$-module par : l'application $\phi^i$ étant défini sur $\Fil^i(N)$, elle se restreint à $N_i$, permettant de poser $\p_N$ égal à $q^i\p^i$ sur $N_i$, c'est-à-dire 
$$\forall x\in N_i,\,\p_N(x)=q^i\phi^i(x)$$
Nous prolongeons cette définition à $N\otimes_WS_0$ tout entier en utilisant la semi-linéarité de $\p_N$. Les propriétés de fonctorialité découlent alors de celles du scindage de la filtration. Au niveau des flèches, ce foncteur est construit de la manière suivante : si $f:N\to N'$ est un morphisme de $\Phi$-modules filtrés, le foncteur lui associe $f\otimes\Id$.
\end{proof}

\begin{rmq}
Nous pouvons étendre ce foncteur de la même façon en un foncteur de la catégorie des $\Phi$-modules filtrés libres sur $W$ vers $\mathop{\bf  \Phi M_{\O_{\E}}^{\text{ét}}}$, qui préserve sous-objet, objet quotient, somme directe, produit tensoriel et dual.
\end{rmq}

N. Wach a montré la proposition suivante (cf. le lemme 3.1.6 p.233 de \cite{Wa}) :
\begin{prop}
Supposons $0\leq h\leq p-2$. Alors pour tout objet $N$ de $\mathop{\bf MF_{W}^h}$, il existe une unique action de $\Gamma_0$ sur $N\otimes_W S_0$ triviale modulo $\pi_0$ et commutant au $\phi_N$ construit comme dans la proposition \ref{propositionbesoin1}. Le module $N\otimes_W S_0$ est alors muni d'une structure de $(\phi,\Gamma_0)$-module sur $S_0$ et devient un objet de $\Gam$..
\end{prop}

C'est le point de départ pour montrer le théorème suivant :
  
\begin{theo}\label{propositionbesoin2}
Supposons $0\leq h\leq p-2$. Il existe un $\otimes$-foncteur $\mathop{\rm F}$\index{$\mathop{\rm  F}$} additif, exact, fidèle et pleinement fidèle de $\mathop{\bf MF_{W,tf}<h>}$ dans $\mathop{\bf  \Gamma_0 \Phi M_{S_0}^{+}}$, qui composé avec le foncteur oubli donne juste le foncteur extension des scalaires de $W$ à $S_0$. De plus, il induit une équivalence de catégories entre $\mathop{\bf MF_{W,tf}^h}$ et $\Gam$, dont un quasi-inverse est $i^*$. 
\end{theo}

\begin{proof}
La {\it première étape} consiste à construire $\mathop{\rm F}$ sur $\mathop{\bf MF_W^h}$. Soit $N$ un objet de $\mathop{\bf MF_W^h}$ (donc libre comme $W$-module). Considérons $N\otimes_W S_0$ : comme $0\leq h\leq p-2$, il existe une unique action de $\Gamma_0$ sur $N\otimes_W S_0$ qui commute à $\phi$ et est triviale modulo $\pi_0$ (c'est le lemme 3.1.6 p.233 de \cite{Wa}). Le $(\phi,\Gamma_0)$ module ainsi défini, noté $\mathop{\rm F}(N)$, est bien un objet de $\Gam$. Il faut voir que nous définissons bien ainsi un foncteur. Comme la structure de $\phi$-module provient d'un scindage de la filtration qui préserve le produit tensoriel, l'unicité de l'action de $\Gamma_0$ nous donnera bien que $\mathop{\rm F}$ préserve le produit tensoriel (tant que celui-ci reste dans $\mathop{\bf MF_{W}^h}$). L'exactitude provient de la même raison.

N. Wach a montré (lemme 3.1.1.2 de \cite{Wa}) qu'il existe un unique générateur topologique $g_0$ de $\Gamma_0$ tel que $\frac{g_0(q)-q}{q\pi_0}=1$ modulo $qS_0$. Il suffit donc d'étudier l'action de $g_0$. Choisissons une base adaptée à la graduation $(e_i)_{1\leq\,i\leq \,d}$ (c'est-à-dire : si $r_i$ est le plus grand entier tel que $e_i\in\Fil^{r_i}(N)$, alors pour tout $r$, $(e_i)_{r_i= r}$ est une base de $N_r$), et si $(a_{i,j})$ est la matrice des applications $\phi^r$ dans cette base, l'action de $\p$ est donné par :

$$\phi(e_j)=q^{r_j}\sum\limits_{1\leq\,i\leq \,d} a_{i,j}e_i$$

Avant de montrer que $\mathop{\rm F}$ préserve les sous objets, nous allons étudier plus en détail l'action de $g_0$.

N. Wach construit l'action de $g_0$ sur $N\otimes_W S_0$ par récurrence modulo $\pi_0^n$. Nous avons besoin de voir cette action d'une autre façon : soit $G=(g_{i,j})$ la matrice dans $GL_{\rg(N)}(S_0)$ définie par $g_0(e_j)=\sum\limits_ig_{i,j}e_i$, et $A=(a_{i,j})\in GL_{\rg(N)}(W)$ donnant l'action de $\phi^j$ sur $e_j$. Alors, en écrivant $\phi\circ g_0(e_j)=\sum\limits_{i,k}\phi(g_{i,j})a_{k,i}q^{r_i}e_k$ et $g_0\circ \phi(e_j)=\sum\limits_{i,k}g(a_{i,j})g_{k,i}g(q)^{r_j}e_k$, la commutativité $\phi\circ g_0=g_0\circ \phi$ nous donne pour $G$ l'équation $AQ\phi(G)=Gg_0(A)g_0(Q)$ avec $Q$ la matrice correspondant à $Q(e_j)=q^{r_j}e_j$ (et $g_0(A)=A$ puisque $A$ est à coefficients dans $W$). Donc $G$ est un point fixe de l'application $f:H\mapsto AQ\phi(H)g_0(Q^{-1})g_0(A^{-1})$ (et le lemme 3.1.6 p.233 de \cite{Wa} affirme juste l'unicité d'un tel point fixe à coefficients dans $S_0$, qui soit congru à $\Id$ modulo $\pi_0$). Notons $I$ la matrice identité dans $GL_{\rg(N)}$ et $G_n=f^{(n)}(I)$ (c'est-à-dire la composée $n$ fois de $f$ appliquée à $I$). Alors, en utilisant que $G-I\in\pi_0M_{\rg(N)}(S_0)$, nous allons montrer :
\begin{lem}
La matrice $G$ est la limite de la suite $G_n$.
\end{lem}

\begin{proof}
Notons $\phi^{(n)}$ la composée $n$ fois de $\phi$ et introduisons alors $B_n=AQ\phi(A)\phi(Q)\cdots \phi^{(n-1)}(A)\phi^{(n-1)}(Q)$ qui est une matrice à coefficients dans $S_0$. Nous avons $G_n=B_n\phi^{(n)}(I)g_0(B_n^{-1})$, et comme $G$ est un point fixe de $f$, $G=B_n\phi^{(n)}(G)g_0(B_n^{-1})$, d'où l'égalité $G_n-G=B_n\phi^{(n)}(I-G)g_0(B_n^{-1})$. Notons $G=I-\pi_0H$ avec $H\in M_{\rg(N)}(S_0)$, alors nous avons $G_n-G=\phi^{(n)}(\pi_0)B_n\phi^{(n)}(H)g_0(B_n^{-1})$. Or, comme $A$ est inversible (dans $GL_{\rg(N)}(W)$), les seuls dénominateurs possibles sont les puissances de $g_0(q)^{r_i}$, et comme $0\leq r_i\leq p-2$, nous pouvons écrire $G_n-G=\frac{\phi^{(n)}(\pi_0)}{g_0\big(q\phi(q)\cdots \phi^{n-1}(q)\big)^{p-2}}G_n'$ avec $G_n'=B_n\phi^{(n)}(H)g_0\big(\phi^{(n-1)}(q^{p-2}Q^{-1}) \phi^{(n-1)}(A^{-1})\cdots q^{p-2}Q^{-1}A^{-1}\big)$ qui est une matrice à coefficients dans $S_0$.

Donc tout revient à montrer que $\frac{\phi^{(n)}(\pi_0)}{g_0\big(q\phi(q)\cdots \phi^{n-1}(q)\big)^{p-2}}$ tend vers $0$. Nous avons $g_0(q)=v_gq$ avec $v_g$ inversible dans $S_0$, par conséquent le fait que $\phi$ et $g_0$ commutent nous donne l'égalité
$$\frac{\phi^{(n)}(\pi_0)}{g_0\big(q\phi(q)\cdots \phi^{n-1}(q)\big)^{p-2}}=\frac{(v_g\phi(v_g)\cdots\phi^{(n-1)}(v_g))^{2-p}}{(q\phi(q)\cdots \phi^{n-1}(q))^{p-2}}\phi^{(n)}(\pi_0)$$ En utilisant que $\phi(\pi_0)=u\pi_0q^{p-1}$ pour $u$ un certain inversible dans $S_0$, nous obtenons que $\phi^{(n)}(\pi_0)=(q\phi(q)\cdots \phi^{n-1}(q))^{p-1}\pi_0u\phi(u)\cdots \phi^{(n-1)}(u)$. Donc, 
$$\frac{\phi^{(n)}(\pi_0)}{g_0\big(q\phi(q)\cdots \phi^{n-1}(q)\big)^{p-2}}=\pi_0\frac{u\phi(u)\cdots \phi^{(n-1)}(u)}{(v_g\phi(v_g)\cdots\phi^{(n-1)}(v_g))^{p-2}}q\phi(q)\cdots\phi^{(n-1)}(q)$$
et, puisque $q\phi(q)\cdots\phi^{(n-1)}(q)$ tend vers $0$ dans $S_0$ ($q$ est dans l'idéal maximal de $S_0$, idéal qui est stable par $\phi$), nous pouvons conclure que $\frac{\phi^{(n)}(\pi_0)}{g_0\big(q\phi(q)\cdots \phi^{n-1}(q)\big)^{p-2}}$ tend vers $0$ dans $S_0$, c'est à dire que $G_n$ tend vers $G$.
\end{proof}
Montrons alors la proposition suivante (qui est le point technique clé de cet article) :
\begin{prop}\label{propositiondur}
Soit $N_{i,j}$ des objets de $\mathop{\bf MF_W^h}$ avec $0\leq h\leq p-2$, $L$ un sous-objet (dans $\mathop{\bf MF_W^+}$) de $M:=\displaystyle \bigoplus_i\otimes_j N_{i,j}$, alors l'action de $\Gamma_0$ sur $\displaystyle \bigoplus_i\otimes_j \mathop{\rm F}(N_{i,j})=M\otimes_W S_0$ laisse stable $L\otimes_W S_0$.
\end{prop}

\begin{proof}
Il suffit de le montrer pour l'action du générateur $g_0$ de $\Gamma_0$. Fixons pour chaque $N_{i,j}$ une base $(e_k^{(i,j)})$ adaptée à la graduation. Notons $G^{(i,j)}$ la matrice de l'action de $g_0$ sur cette base et $C^{(i,j)}$ la matrice donnant l'action de $\phi$ sur $N_{i,j}\otimes_W S_0$ (avec les notations précédentes, $C=AQ$). Alors, par le lemme précédent nous avons $\lim\limits_{n\to+\infty}G_n^{(i,j)}=G ^{(i,j)}$ avec $C^{(i,j)}\phi(G_n^{(i,j)})g_0(C^{(i,j)})^{-1}=G_{n+1}^{(i,j)}$ et $G_0 ^{(i,j)}=I ^{(i,j)}$. 

Prenons $(u[l])_l$ une base de $L$, et notons $(u[l]_k^{(i,j)})$ les coordonnées de $u[l]$ dans la base $(e_k^{(i,j)})$. Nous voulons montrer (par récurrence) que $\bigoplus_i\otimes_j G_{n+1}^{(i,j)}g_0(u_k^{(i,j)})$ est une combinaison linéaire (à coefficients dans $S_0$) des $(u[l]_k^{(i,j)})$, pour $u$ élément quelconque de $L\otimes_W S$ (et $(u_k^{(i,j)})$ ses coordonnées). Remarquons que par linéarité, il suffit de le montrer pour $u$ égal aux $u[l]$.

Comme $L$ est un sous-objet de $M$, nous avons $L\otimes_W S_0$ qui est stable par $\phi$. Or $\phi$ induit une bijection de $L\otimes_WS_0[\frac{1}{q}]$. Cela se traduit alors en disant $\bigoplus_i\otimes_jC^{(i,j)}\phi(u[l']_k^{(i,j)})$ est une combinaison linéaire (à coefficients dans $S_0$) des $(u[l]_k^{(i,j)})$, et qu'il existe $N\in\N$ tel que $q^N\big(\otimes_jC^{(i,j)}\big)^{-1}(u[l']_k^{(i,j)})$ est une combinaison linéaire (à coefficients dans $S_0$) des $\phi(u[l]_k^{(i,j)})$.

Par conséquent, $g_0(q)^Ng_0\big(\otimes_jC^{(i,j)}\big)^{-1}\big(g_0(u[l']_k^{(i,j)})\big)$ s'écrit comme une combinaison linéaire (à coefficients dans $S_0$) des $\big(g_0(\phi(u[l]_k^{(i,j)}))\big)$, ceci pour tout $l'$.

Puis, $\bigoplus_i\otimes_j \phi(G_n^{(i,j)})g_0\big(\phi(u[l']_k^{(i,j)})\big)=\phi\big(\bigoplus_i\otimes_j G_n^{(i,j)}g_0(u[l']_k^{(i,j)})\big)$ est pour tout $l'$ une combinaison linéaire (à coefficients dans $S_0$) des $\big(\phi(u[l]_k^{(i,j)})\big)$, cela provient de notre hypothèse de récurrence.

En reprenant que $\bigoplus_i\otimes_jC^{(i,j)}\phi(u[l']_k^{(i,j)})$ est une combinaison linéaire (à coefficients dans $S_0$) des $(u[l]_k^{(i,j)})$ pour tout $l'$, et en mettant bout à bout ces affirmations, nous obtenons que \begin{multline*}g_0(q)^N\bigoplus_i\otimes_j G_{n+1}^{(i,j)}g_0(u[l']_k^{(i,j)})=\\g_0(q)^N\bigoplus_i \otimes C^{(i,j)}\otimes\phi(G_n^{(i,j)})g_0\big(\otimes C^{(i,j)}\big)^{-1}(g_0(u[l']_k^{(i,j)}))\end{multline*} est pour tout $l'$ une combinaison linéaire (à coefficients dans $S_0$) des $(u[l]_k^{(i,j)})$.

Par conséquent, si $g^{[n]}$ désigne l'application $g_0$-linéaire construite à partir de la matrice $\bigoplus_i\otimes_j G_{n}^{(i,j)}$ (l'hypothèse de récurrence se traduisant par : $L\otimes_WS_0$ est stable par $g^{[n]}$), alors $g^{[n+1]}(L\otimes_WS_0)\subset \frac{1}{g_0(q)^N}L\otimes_WS_0=\frac{1}{q^N}L\otimes_WS_0$. Considérons alors $(f_r)_{1\leq r\leq \rg_W(M)}$ une base de $M$ telle qu'il existe $n_r\in\N\cup\{+\infty\}$ avec $(p^{n_r}f_r)$ base de $L$. Alors $g^{[n+1]}(p^{n_r}f_r)=\sum_s\frac{\a_s}{q^N}p^{n_s}f_s$ avec $\a_s\in S_0$ (qui dépend de $r$). Mais, par construction, $g^{[n+1]}(M\otimes_WS_0)\subset M\otimes_WS_0$, alors $g^{[n+1]}(p^{n_r}f_r)=\sum_sp^{n_r}\beta_s f_s$ avec $\beta_s\in S_0$ (qui dépend aussi de $r$). D'où $p^{n_r}\beta_s=\frac{\a_s}{q^N}p^{n_s}$, ce qui implique que $q^N$ divise $\a_s$ dans $S_0$, donc que $g^{[n+1]}(L\otimes_WS_0)\subset L\otimes_WS_0$, ce qui montre bien la récurrence. 

Pour initialiser la récurrence ($n=0$) nous avons $G_0 ^{(i,j)}=I ^{(i,j)}$ (où $I$ est la matrice identité), donc $\bigoplus_i\otimes_j G_{0}^{(i,j)}g_0(u_k^{(i,j)})=u_k^{(i,j)}$ pour tout $u$ dans $L$. D'où par récurrence la propriété est vraie pour tout $n$. En passant à la limite, la propriété est vrai pour $\bigoplus_i\otimes_j G^{(i,j)}$. Donc l'action de $g_0$ sur $M\otimes_W S_0$ laisse stable $L\otimes_W S_0$. 
\end{proof}

Cette proposition est le coeur du théorème. Elle nous donne en particulier que si $N'$ est un sous-objet de $N$ dans $\mathop{\bf MF_{W}^h}$, alors l'action de $g_0$ sur $N\otimes_WS_0$ laisse stable $N'\otimes_WS_0$. Elle est triviale modulo $\pi_0$ : si $(e_i)$ est une base de $N$, telle qu'il existe $(\a_i)\in\N\cup\{+\infty\}$ avec $(p^{\a_i}e_i)$ base de $N'$, alors il existe des coefficients $x_{i,j}$ et $y_{i,j}$ dans $S_0$ tels que $g_0(e_i)=e_i+\pi\sum_{j\neq i}x_{i,j}e_j$ et $g_0(p^{\a_i}e_i)=\sum_jy_{i,j}p^{\a_j}e_j$. En identifiant les coordonnées, nous avons $y_{i,i}=1$ et $y_{i,j}p^{\a_j}=p^{\a_i}x_{i,j}\pi$ si $j\neq i$, donc $\pi$ divise bien $y_{i,j}$ dans $S_0$ pour $j\neq i$. Donc l'action de $g_0$ sur $\mathop{\rm F}(N)$ se restreint en une action triviale modulo $\pi_0$ sur $N'\otimes_WS_0$ qui commute à $\phi$, donc par unicité cette action est celle de $\mathop{\rm F}(N')$.

La {\it deuxième étape} consiste alors à définir $\mathop{\rm F}$ sur tout $\mathop{\bf MF_{W,tf}<h>}$. Le point important est que pour tout objet $M$ de $\mathop{\bf MF_{W,tf}<h>}$, il existe des objets $N_{i,j}$ dans $\mathop{\bf MF_W^h}$ et $L$ un sous-objet de $\bigoplus_i \otimes_j N_{i,j}$ tels que $M$ est isomorphe à un quotient $M'$ de $L$. Considérons alors $N$ un sous-objet de $M$, et supposons que sur $M\otimes_WS_0$ nous ayons une structure de $(\phi,\Gamma_0)$-module qui le rende isomorphe à $M'\otimes_WS_0$ muni de la structure de $(\phi,\Gamma_0)$-module obtenu à partir de celle de $L\otimes_WS_0$ donnée par la proposition \ref{propositiondur}. Il faut voir que $N'\otimes_WS_0$ (où $N'$ est l'image de $N$ dans $M'$) est stable par $\Gamma_0$. En notant $\pi:L\to M'$ la projection naturelle, $\pi^{-1}(N')$ est un sous-objet de $L$ (car c'est le noyau du morphisme $L\to M'/N'$, donc par la remarque 1.4.2 et la proposition 1.4.1 de \cite{Wi}, c'est bien un sous-objet de $L$), donc la proposition \ref{propositiondur} nous donne bien que $\pi^{-1}(N')\otimes_WS_0$ est stable par l'action de $\Gamma_0$. Par conséquent, $N\otimes_WS_0$ sera bien laissé stable par l'action de $\Gamma_0$ de $M\otimes_WS_0$, donc sera un sous-$(\phi,\Gamma_0)$-module de $M\otimes_WS_0$.

Puis, $\mathop{\rm F}$ se construit par itération : notons $\mathop{\bf MF_n}$ la sous-catégorie pleine de $ \mathop{\bf MF_{W,tf}}$, construite en disant qu'un objet de $\mathop{\bf MF_{n+1}}$ est soit le sous-objet ou le quotient d'un objet de $\mathop{\bf MF_n}$, soit la somme directe de deux objets de $\mathop{\bf MF_n}$, soit le produit tensoriel de deux objets de $\mathop{\bf MF_n}$, soit un objet de $\mathop{\bf MF_n}$, et posons (pour initialiser la récurrence) $\mathop{\bf MF_0}=\mathop{\bf MF_W^h}$. Alors, $\mathop{\bf MF_{W,tf}<h>}=\bigcup \mathop{\bf MF_n}$, et si $\mathop{\rm F}$ est construit sur $\mathop{\bf MF_n}$, alors il s'étend naturellement à la somme directe ou le produit tensoriel de deux objets de $\mathop{\bf MF_n}$, et l'étude précédente montre qu'il s'étend au cas d'un sous-objet, et donc d'un objet quotient, d'un objet de $\mathop{\bf MF_n}$. Nous pouvons donc donner naturellement une structure de $(\phi,\Gamma_0)$-module à tout $M\otimes_WS_0$, pour $M$ un objet de $\mathop{\bf MF_{W,tf}<h>}$.

La {\it troisième étape} s'occupe des morphismes. Pour $M$ et $M'$ deux objets de $\mathop{\bf MF_{W,tf}<h>}$, et $f:M\to M'$ un morphisme, nous posons $\mathop{\rm F}(f)=f\otimes\Id$ (En particulier, le foncteur sera exact (car $S_0$ est plat sur $W$) et fidèle). C'est un morphisme de $(\phi,\Gamma_0)$ module : par construction, c'est un morphisme de $\phi$-modules. Puis, $f\oplus\Id:M\oplus M'\to M'\oplus M'$ est un morphisme de $\phi$-modules filtrés, donc par la proposition 1.4.1 de \cite{Wi}, $\Ker(f\oplus\Id)=\{x-y|x\in M,\,y\in M',\,y=f(x)\}$ est naturellement un $\Phi$-module filtré, sous-objet de $M\oplus M'$, donc par la proposition \ref{propositiondur}, $\Ker(f\oplus\Id)\otimes_WS_0$ est laissé stable par l'action naturelle de $\Gamma_0$ sur $(M\oplus M')\otimes_WS_0$ (obtenue à partir de celle sur $M\otimes_WS_0$ et $M'\otimes_WS_0$), donc $f$ commute à l'action de $\Gamma_0$, car si $x-y\in \Ker(f\oplus\Id)\otimes_WS_0$, dire que $g(x)-g(y)\in \Ker(f\oplus\Id)\otimes_WS_0$, c'est dire que $f(g(y))=g(x)=g(f(y))$. 

Montrons la pleine fidélité. Remarquons que si nous munissons $M\otimes_WS_0$ de la structure de $\Phi$-module filtré donnée par $\Fil^i(M\otimes_WS_0)=\{x\in M\otimes_WS_0|\phi(x)\in q^iM\otimes_WS_0\}$, et $\phi^i=\frac{1}{q^i}\phi$, alors la restriction modulo $\pi_0$ est un morphisme de $\Phi$-module filtré. Par conséquent, si $f:M\otimes_WS_0\to M'\otimes_WS_0$ est un morphisme de $(\phi,\Gamma_0)$-module, alors la réduction modulo $\pi_0$ induit $\bar f:M\to M'$ un morphisme de $\phi$-modules filtrés. Donc $\bar f\otimes\Id : M\otimes_WS_0\to M'\otimes_WS_0$ est un morphisme de $(\phi,\Gamma_0)$-module par le résultat précédent, donc $g=f-\bar f\otimes\Id $ aussi, et il se réduit modulo $\pi_0$ sur l'application nulle. Notons $\mathcal{M}''$ le noyau de $g$, c'est un sous $(\phi,\Gamma_0)$-module de $M\otimes_WS_0$. Nous avons $\pi_0\mathcal{M}''=\mathcal{M}''\cap \pi_0M\otimes_WS_0$ car $M'\otimes_WS_0$ est sans $\pi_0$-torsion. Donc $\mathcal{M}''/\pi_0\mathcal{M}''\subset M$, et par le lemme du serpent, nous avons égalité, car :
\[\xymatrix{
 0\ar[r]&\mathcal{M}''\ar[r]\ar[d]^-{u_1}& M\otimes_WS_0 \ar[r]^-{g}\ar[d]^-{u_2}&M\otimes_WS_0\ar[r]\ar[d]^-{u_3}&0\\
0\ar[r]&M\ar[r]& M \ar[r]^-{\bar g}&M\ar[r]&0
}\]
les lignes horizontales sont exactes, $u_1$ est la réduction modulo $\pi_0$ composée avec l'inclusion $\mathcal{M}''/\pi_0\mathcal{M}''\subset M$, $u_2$ et $u_3$ sont la réduction modulo $\pi_0$, $ u_2$ est surjectif ($u_3$ aussi), et l'application naturelle $\Ker(u_2)=\pi_0 M\otimes_WS_0 \to \Ker(u_3)=\pi_0 M\otimes_WS_0$ est surjective. Donc, nous avons $M\otimes_WS_0=\mathcal{M}''+\pi_0M\otimes_WS_0$, l'idéal engendré par $\pi_0$ est inclus dans le radical de Jacobson de $S_0$, $M\otimes_WS_0$ est de type fini sur $S_0$, donc par le lemme de Nakayama, nous avons $M\otimes_WS_0=\mathcal{M}''$, donc $f=\bar f\otimes\Id $. Par conséquent le foncteur est pleinement fidèle.

La {\it quatrième étape} est l'étude du foncteur restreint à $\mathop{\bf MF_{W,tf}^h}$. Par construction, nous avons $i^*\mathop{\rm F}(N)=N$ pour tout objet $N$ de $\mathop{\bf MF_{W,tf}^h}$. Montrons :
\begin{lem}
Pour tout objet $\mathcal{N}$ de $\Gam$, libre comme $S_0$-module, si $N=i^*(\mathcal{N})$, il existe un unique isomorphisme de $(\p,\Gamma_0)$-module $\mathop{\rm F}(N)\to \mathcal{N}$ (qui se réduit modulo $\pi_0$ sur l'égalité $N=i^*(\mathcal{N})$).
\end{lem}

\begin{proof}Présentons ici une démonstration de ce fait due à N. Wach. Pour cela, considérons une base $(e_i)_{1\leq\,i\leq \,d}$ de $N$, adaptée à la graduation, et $(a_{i,j})$ la matrice des applications $\phi^r$ dans cette base (donc l'action de $\p$ est donnée sur $\mathop{\rm F}(N)$ par $\phi(e_j)=q^{r_j}\sum\limits_{1\leq\,i\leq \,d} a_{i,j}e_i$). Il faut alors prouver l'existence et l'unicité d'une base $(f_i)$ dans $\mathcal{N}$ vérifiant $\phi(f_j)=q^{r_j}\sum\limits_{1\leq\,i\leq \,d} a_{i,j}f_i$ avec $e_i=f_i$ modulo $\pi_0$. Ce sera suffisant car en posant $h(e_i)=f_i$, nous aurons un morphisme de $\phi$-module, qui fera commuter l'action de $\Gamma_0$ par unicité de celle-ci, et qui modulo $\pi_0$ redonnera l'identité.

Par construction du $\Phi$-module filtré $N$, la base $(e_i)$ se relève en une famille $(\hat e_i)$ de $\mathcal{N}$ avec $\p(\hat e_i)\in q^{r_i}\mathcal{N}$. De plus, $\mathcal{N}$ est complet pour la topologie $\pi_0$-adique (car $S_0$ l'est), et modulo $\pi_0$, $(e_i)$ est une base, donc $\mathcal{N}$ étant sans torsion, $(\hat e_i)$ est une base de $\mathcal{N}$ (nous pourrions aussi invoquer le lemme de Nakayama). Donc, il existe $\hat a_{i,j}\in S_0$ tels que :
$$\p(\hat e_j)=q^{r_j}\sum\limits_{1\leq i\leq d}\hat a_{i,j}\hat e_i$$
et $\hat a_{i,j}=a_{i,j}$ modulo $\pi_0$. Posons $\a_{i,j}\in S_0$ l'unique élément tel que $\hat a_{i,j}=a_{i,j}+\pi_0\a_{i,j}$. Nous cherchons à modifier la base $(\hat e_i)$ pour obtenir la base $(f_i)$. Cherchons $f_j$ sous la forme $f_j=\hat e_j+\pi_0c_j$, et posons $b_j=\sum\limits_{1\leq i\leq d} \a_{i,j} \hat e_i$. Alors, puisque $\phi(\pi_0)=uq^{p-1}\pi_0$, 


$\p(\hat e_j+\pi_0c_j)=\p(\hat e_j)+u\pi_0q^{p-1}\p(c_j)$, et en faisant appara\^itre $\sum\limits_{i=1}^d q^{r_j}a_{i,j}\pi_0c_i$, nous obtenons


$\p(\hat e_j+\pi_0c_j)=\sum\limits_{i=1}^d  q^{r_j}a_{i,j}(\hat e_i+\pi_0c_i)+\pi_0q^{r_j}b_j +u\pi_0q^{p-1}\p(c_j)-\sum\limits_{i=1}^d q^{r_j}a_{i,j}\pi_0c_i$

autrement dit, nous cherchons les $c_j\in \mathcal{N}$ tels que 
$$b_j +uq^{p-1-r_j}\p(c_j)-\sum\limits_{1\leq i\leq d} a_{i,j}c_i=0$$
Nous résolvons ce système de manière unique par récurrence modulo $\pi_0^n$. A chaque étape, le système se résout en faisant une récurrence modulo $p^k$, en utilisant que $p-1-r_j\geq 1$ (par hypothèse), donc que $q^{p-1-r_j}=0$ modulo $(p,\pi_0)$, et que la matrice $(a_{i,j})$ est inversible modulo $p$.\end{proof}

Pour terminer la démonstration du théorème (c'est à dire prouver le lemme précédent sans l'hypothèse sur la liberté de $N$), nous aurons besoin de résultats sur les modules de Wach, qui apparaitront plus loin dans l'article.
\end{proof}

Pour la suite, nous aurons besoin de faire intervenir un foncteur légèrement différent. Si $N$ est un objet de $\mathop{\bf MF_W^{-h}}$, son dual $N^*=\Hom_{\Z_p}(N,\Z_p)$ est un objet de $\mathop{\bf MF_W^h}$, donc $\mathop{\rm F}(N^*)$ est bien défini.
\begin{defi}
Le foncteur $\mathop{\rm F}^-$\index{$\mathop{\rm F}^-$} est défini sur $\mathop{\bf MF_W^{-h}}$ pour $h\leq p-2$ par :
$$\begin{array}{c}\mathop{\rm F}^-(N)=\big(\mathop{\rm F}(N^*)\otimes_{S_0}\O_{\E}\big)^*\end{array}=N\otimes_W\O_{\E}$$
pour tout objet $N$ de $\mathop{\bf MF_W^{-h}}$. Il donne bien un $(\phi,\Gamma)$-module étale sur $\O_{\E}$. Il s'étend de même à $ \mathop{\bf MF_W<-h>}$.
\end{defi}
Le foncteur $\mathop{\rm F}$ consiste à munir $N\otimes_W S_0$ (pour $N$ objet de $\mathop{\bf MF_W<h>}$) d'une structure de $(\phi,\Gamma_0)$-module, et comme nous voulons un foncteur défini sur $\MF$, nous prenons le dual. Ceci pose néanmoins un problème de définition, car le dual d'un $\phi$-module sur $S_0$ n'est pas un $\phi$-module, c'est pourquoi nous étendons d'abord les scalaires à $\O_{\E}$, puis nous prenons le dual.

Remarquons que $\mathop{\rm F}^-$ peut être défini sur $\mathop{\bf MF_{W,tf}^{-h}}$ (puis sur $\mathop{\bf MF_{W,tf}<-h>}$) en prenant pour un module de $p$-torsion le dual de Pontriaguine, et en passant à la limite projective pour le cas général.

\begin{rmq}
Nous pouvons définir $\mathop{\rm \widetilde F}$ sur $\mathop{\bf MF_W^{\pm h}}$ pour $h\leq \frac{p-2}{2}$ en posant $\mathop{\rm \widetilde F}(N)=\mathop{\rm F}(N\otimes_W W[h])\otimes_{S_0}\O_{\E}[-h]$ avec $\O_{\E}[-h]=\mathop{\rm F}(W[h])^*$ et $W[h]$ l'objet de $\MF$ dont le $W$-module sous-jacent est $W$, avec $\Fil^i(W[-h])=\left\{\begin{array}{c} W \text{ si }i\leq -h\\ 0\text{ si }i>-h\end{array}\right.$ et $\p^{-h}(x)=\s(x)$. Alors, $\mathop{\rm \widetilde F}(N^*$ est canoniquement isomorphe à $\mathop{\rm \widetilde F}(N)^*$ (cela se voit à l'aide de l'unicité de l'action de $\Gamma_0$ agissant trivialement modulo $\pi_0$, et commutant à $\phi$, d'après le lemme 3.1.6 de \cite{Wa}), et $\mathop{\rm \widetilde F}$ s'étend alors en un foncteur sur $\mathop{\bf MF_W<\pm h>}$ qui a des propriétés similaires à celles de $\mathop{\rm F}$, et qui préserve le dual.
\end{rmq}

\section{Lien entre le foncteur et les modules de Wach}
\subsection{Fonctorialité de $g_N$}\label{construction}
Rappelons le Théorème 1' de N. Wach (cf. \cite{Wa}) :
\begin{theore}
Si $\mathcal{N}$ est un objet de $\Gam$ avec $0\leq h\leq p-2$, alors $\mathop{\rm Hom_{\bf MF_W}}(i^*(\mathcal{N}), A_{cris})$ est isomorphe (en tant que représentation galoisienne) à $\mathop{\rm Hom_{\bf \Phi M_{S_0}}} (\mathcal{N},\O_{\widehat\E_{nr}})$.
\end{theore}
Enoncé dans le cadre (et avec les notations) qui nous intéresse, il devient :
\begin{theore}\label{theowach}
Si $N$ est un objet de $\mathop{\bf MF_{W,tf}^{-h}}$ avec $0\leq h\leq p-2$, alors il existe un isomorphisme $g_N : \mathop{\bf V_{\O_{\E}}}(\mathop{\rm F}^-(N))\to \V(N)$ de représentations galoisiennes. Si en plus $N$ est libre, en passant au dual, cela donne un isomorphisme de représentations galoisiennes ${}^tg_{N}^{-1}:\mathop{\bf V_{\O_{\E}}}\big(\mathop{\rm F}(N^*)\otimes_{S_0}\O_{\E}\big)\to \V(N)^*$.
\end{theore}
Nous allons vérifier que cet isomorphisme est fonctoriel :
\begin{theo}\label{theopartiel}
Pour tout objet $N$ de $\mathop{\bf MF_{W,tf}^{-h}}$ avec $0\leq h\leq p-2$, l'application $g_N$ construite par N. Wach vérifie les propriétés de fonctorialité suivante :
\begin{enumerate}
\item\label{morphisme} pour tout morphisme $f:N\to N'$ entre deux objets $N$ et $N'$ de $\mathop{\bf MF_{W,tf}^{-h}}$, nous avons $\V(f)\circ g_N=g_{N'}\circ \mathop{\bf V_{\O_{\E}}}(\mathop{\rm F}^-(f))$ (cela s'applique en particulier pour l'injection d'un sous-objet, ou pour la projection sur un objet quotient).
\item pour tout objet $N$ et $N'$ de $\mathop{\bf MF_{W,tf}^{-h}}$, $g_{N\oplus N'}=g_N\oplus g_{N'}$ ;
\item\label{produittensoriel} pour tout objet $N$ et $N'$ de $\mathop{\bf MF_{W,tf}^{-h}}$, pour tout sous-objet $L$ de $N\otimes N'$ tel que $L$ soit un objet de $\mathop{\bf MF_{W,tf}^{-h}}$, l'application $g_{N}\otimes g_{N'}$ restreinte à $ \mathop{\bf V_{\O_{\E}}}(\mathop{\rm F}^-(L))$ est égale à $g_{L}$. En particulier, si $N\otimes N'$ est un objet de $\mathop{\bf MF_{W,tf}^{-h}}$, alors $g_{N\otimes N'}=g_N\otimes g_{N'}$ ;
\end{enumerate}
\end{theo}

\begin{rmq}\label{pcp}
Le point (\ref{produittensoriel}) montre en particulier que $\V(N\otimes N')$ est égal à $\V(N)\otimes_{\Z_p}\V(N')$ dès que $N$, $N'$ et $N\otimes N'$ sont des objets de $\mathop{\bf MF_{W,tf}^{-h}}$ avec $0\leq h\leq p-2$.
\end{rmq}

Rappelons la construction de $g_N$ : N. Wach construit l'isomorphisme modulo $p^n$ pour tout $n$ à partir des morphismes d'anneaux (avec $A_S^+=W(R)\cap \O_{\widehat{\E}_{nr}}$) : $A_S^{+}/p^n\to W_n(R)/\pi_0$ et $A_{cris}/p^n\to W_n(R)/\pi_0$. Notons $\mathcal{N}:=\mathop{\rm F}^-(N)$. Nous avons la bijection $ \mathop{\bf V_{A_S^+}}(\mathcal{N}/p^n):={ (\mathcal{N}\otimes_{S_0} A_S^{+}/p^n)^{\p}}  \to {(\mathcal{N}\otimes_{S_0} \O_{\E_{nr}/p^n})^{\p}}$ (cf \cite{F91}, p.296, où c'est exprimé pour le foncteur contravariant). Or, N. Wach a montré que pour $N$ objet de $\mathop{\bf MF_W^{-h}}$ avec $0\leq h\leq p-2$, le schéma suivant 
\[\xymatrix{
N/p^n\otimes_WA_{cris}/p^n \ar[r]^-{k_N}& N/p^n\otimes_WW_n(R)/\pi_0& \mathcal{N}/p^n\otimes_{S_0} A_S^{+}/p^n \ar[l]_-{j_{\mathcal{N}}}\\ \V(N/p^n)\ar@{^{(}-{>}}[u]_k&&\mathop{\bf V_{A_S^+}}(\mathcal{N}/p^n)\ar@{^{(}-{>}}[u]_j
}\]
induit un isomorphisme de représentations galoisiennes de $\mathop{\bf V_{A_S^+}}(\mathcal{N}/p^n)$ sur $\V(N/p^n)$, c'est-à-dire que $K_{N}=k_{N}\circ k$ et $J_{\mathcal{N}}=j_{\mathcal{N}}\circ j$ sont toutes les deux injectives, et ont même image dans $ N/p^n\otimes_WW_n(R)/\pi_0$.

Tout ceci passe à la limite projective, et nous obtenons l'application $g_N$ bijective :
\[\xymatrix{
N\otimes_WA_{cris} \ar[r]^-{k_{N}}& N\otimes_WW(R)/\pi_0& \mathcal{N}\otimes_{S_0} A_S^{+} \ar[l]_-{j_{\mathcal{N}}}\\ \V(N)\ar@{^{(}-{>}}[u]_k\ar[ru]_{K_{N}}&&\mathop{\bf V_{A_S^+}}(\mathcal{N})\ar@{^{(}-{>}}[u]_j\ar[lu]^{J_{\mathcal{N}}}\ar[ll]^{g_N}
}\]
où $\mathop{\bf V_{A_S^+}}(\mathcal{N})=\varprojlim_{n\in \N} \mathop{\bf V_{A_S^+}}(\mathcal{N}/p^n)=(\mathcal{N}\otimes_{S_0}A_S^+)^{\p=1}=\mathop{\bf V_{\O_{\E}}}(\mathcal{N}\otimes_{S_0}\O_{\E})$.

\begin{proof}[Démonstration du théorème \ref{theopartiel}] Pour la fonctorialité au niveau des flèches, il suffit de remarquer que le diagramme suivant est commutatif (car $\mathop{\bf V_{\O_{\E}}}(\mathop{\rm F}^-(f))$ est juste $f\otimes\Id$) :
\[\xymatrix{
\mathop{\bf V_{\O_{\E}}}(\mathop{\rm F}^-(N))\ar[r]^-{\mathop{\bf V_{\O_{\E}}}(\mathop{\rm F}^-(f))}\ar@{^{(}->}[d]& \mathop{\bf V_{\O_{\E}}}(\mathop{\rm F}^*(N'))\ar@{^{(}->}[d]\\
N\otimes A_S^+ \ar[d]_-{j_{\mathcal{L}}}\ar[r]^-{f\otimes\Id}&N'\otimes A_S^+\ar[d]^-{j_{\mathcal{L}'}}\\
N\otimes W(R)/\pi_0 \ar[r]^-{f\otimes\Id}&N'\otimes W(R)/\pi_0 \\
N\otimes A_{cris}\ar[u]^-{k_L}\ar[r]^-{f\otimes\Id} &N'\otimes A_{cris}\ar[u]_-{k_{L'}}\\
\V(N)\ar@{_{(}->}[u]\ar[r]^-{\V(f)}&\V(N')\ar@{_{(}->}[u]
}\]

Le fait que $g_{N\oplus N'}=g_N\oplus g_{N'}$ se montre de la même façon. Il reste donc à voir le cas du produit tensoriel : considérons $N$ et $N'$ deux objets de $\mathop{\bf MF_{W,tf}^{-h}}$ avec $h\leq p-2$. Soit $L$ un sous-objet de $N\otimes N'$ qui est dans $\mathop{\bf MF_{W,tf}^{-h}}$, posons $\mathcal{L}=L\otimes_WS_0$. Le diagramme suivant est alors commutatif :

 \[\xymatrix{
 L\otimes_WA_{cris} \ar[d]^-{k_{L}}\ar@{^{(}-{>}}[r] &(N\otimes_WA_{cris})\otimes_{A_{cris}} (N'\otimes_WA_{cris}) \ar[d]^-{k_{N}\otimes k_{N'}}\\ L\otimes_WW(R)/\pi_0\ar@{^{(}-{>}}[r]&  \big(N\otimes_WW(R)/\pi_0\big)\otimes_{W(R)/\pi_0} \big(N'\otimes_WW(R)/\pi_0\big)\\ \mathcal{L} \otimes_{S_0} A_S^{+} \ar[u]_-{j_{\mathcal{L}}}\ar@{^{(}-{>}}[r]& (\mathcal{N}\otimes_{S_0} A_S^{+})\otimes_{A_S^+} (\mathcal{N}'\otimes_{S_0} A_S^{+})  \ar[u]_-{j_{\mathcal{N}}\otimes j_{\mathcal{N}'}}
}\]

Par conséquent, l'application $K_{N}\otimes K_{N'}$ restreinte à $\V(L)$ est égale à $K_{L}$, et l'application $J_{\mathcal{N}}\otimes J_{\mathcal{N}'}$ restreinte à $\mathop{\bf V_{A_S^+}}(\mathcal{L})$ est égale à $J_{\mathcal{L}}$. 

Le point important est que $L$ étant un objet de $\mathop{\bf MF_{W,tf}^{-h}}$ (par hypothèse), ce sont bien des bijections, et ce sont celles qui permettent de construire $g_L$. 

Donc $g_N\otimes g_{N'}$ envoie $\mathop{\bf V_{A_S^+}}(\mathcal{L})$ sur $\V(L)$ si $L$ est un sous-objet de $N\otimes N'$ qui soit dans $\mathop{\bf MF_{W,tf}^{-h}}$, et plus exactement, l'application $g_N\otimes g_{N'}$ restreinte à $\mathop{\bf V_{A_S^+}}(\mathcal{L})$ est égale à $g_L$.

Si $N\otimes N'$ est un objet de $\mathop{\bf MF_{W,tf}^{-h}}$, le résultat précédent avec $L=N\otimes N'$ nous donne $g_{N\otimes N'}=g_{N}\otimes g_{N'}$.
\end{proof}
\begin{rmq}\label{gggg}
Nous montrons de même que pour $(N_{i,j})_{1\leq j\leq n,\,1\leq i\leq n_j}$ objets de $\mathop{\bf MF_{W,tf}^{-h}}$ avec $0\leq h\leq p-2$, et pour $L$ un sous-objet de $\displaystyle\bigoplus_{j=1}^n\otimes_{i=1}^{n_j} N_{i,j}$ qui soit dans $\mathop{\bf MF_{W,tf}^{-h}}$, alors $\bigoplus\otimes g_{\mathop{\rm F}^-(N_{i,j})}$ restreinte à $\mathop{\bf V_{\O_{\E}}}(\mathop{\rm F}^-(L))$ est égale à $g_L$.
\end{rmq}

Nous pouvons traduire ces résultats en disant :
\begin{theo}\label{theosemifinal}
Soit $0\leq h\leq p-2$, et notons $\mathcal{G}$ le foncteur exact de la catégorie $\mathop{\bf MF_{W,tf}^{-h}}$ vers la catégorie des représentations continues de $\Gamma_{\K}$ sur les $ \Z_p$-modules de rang fini, défini par : si $N$ objet de $\mathop{\bf MF_{W,tf}^{-h}}$, $\mathcal{G}(N)=\mathop{\bf V_{ \O_{\E} }}(\mathop{\rm F}^-(N))$. Alors il existe $g$ un isomorphisme de foncteurs entre $\mathcal{G}$ et $\V$. De plus, nous pouvons supposer que :
\begin{itemize}
\item pour tous objet $N$ et $N'$ de $\mathop{\bf MF_{W,tf}^{-h}}$, tel que $N\otimes N'$ soit encore un objet de $\mathop{\bf MF_{W,tf}^{-h}}$, nous avons $g_{N\otimes N'}=g_{N}\otimes g_{N'}$ ;
\item pour tout uplet d'objets $(N_{i,j})_{1\leq j\leq n,\,1\leq i\leq n_j}$ de $\mathop{\bf MF_{W,tf}^{-h}}$, pour tout sous-objet $L$ (dans $\mathop{\bf MF_{W,tf}^{-h}}$) de $\displaystyle\bigoplus_{j=1}^n\otimes_{i=1}^{n_j} N_{i,j}$, l'application $\bigoplus\otimes g_{N_{i,j}}$ restreinte à $\mathop{\bf V_{ \O_{\widehat{\E}} }}(\mathop{\rm F}^-(L))$ est égale à $g_{L}$.
\end{itemize}
\end{theo}

\subsection{Lien entre $\Gam$ et $\mathop{\bf \Gamma\Phi M_S^h}$}
Avant de parler de modules de Wach (qui sont des $S$-modules), il faut comprendre l'extension des scalaires $S_0\to S$.

\begin{lem}
$S=\bigoplus\limits_{0\leq i\leq p-2}S_{i}$, où si $x\in S_{i}$ et $g\in\Gamma_f$ est $[\a]$ (le relèvement de Teichmuller de $\a\in\F_p^*$), alors $g$ agit sur $x$ par $g(x)=[\a]^ix$.
\end{lem}

\begin{proof}
L'application $p_i=\frac{1}{|\Gamma_f|}\sum\limits_{g\in \Gamma_f}\X(g)^{-i}g$ est un projecteur dont l'image est $S_i$, et les $p_i$ vérifient $\sum\limits_{0\leq i\leq p-2} p_i=\Id$.
\end{proof}

\begin{lem}\label{lemmess0}
$S$ a une base normale sur $S_0$, c'est à dire qu'il existe $e\in S$ tel que $(g(e))_{g\in\Gamma_f}$ soit une base de $S$ sur $S_0$. De plus, $p$ ne divise aucun $p_i(e)$.
\end{lem}
\begin{proof}
En effet, il suffit de le montrer modulo $p$ (et ensuite de relever une base normale de $k[[\pi]]$ sur $k[[\pi_0]]$, puisque $S_0$ est complet pour la topologie $p$-adique). Or, Fontaine a montré dans \cite{F91}, page 270, que le corps des fractions de $k[[\pi]]$, $k((\pi))$, est une extension galoisienne cyclique de degré $p-1$ (donc modérément ramifiée) de $k((\pi_0))$, dont le groupe de Galois est donné par $\Gamma_f$. Donc, par un théorème de E. Noether, il existe une base normale pour les anneaux d'entiers correspondants. Enfin, si $\overline e$ est cette base (modulo $p$), alors $p_i(\overline e)=\sum_g\frac{\X(g)^{-i}}{|\Gamma_f|}g(\overline e)$ est non nul (puisque chaque coordonnée suivant la base $(g(\overline e))$ est non nulle (même modulo $p$)), donc $p_i(e)$ sera bien non divisible par $p$ si $e$ relève $\overline e$. 
\end{proof}

En particulier, nous avons $S_i=p_i(e)S_0$ (car $p_i(e)S_0\subset S_i$, puis $e\in \oplus_i p_i(e)S_0$, $\oplus_ip_i(e)S_0$ est donc un $S_0$-module contenant $e$ et stable par $\Gamma_f$, donc $S=\oplus_i p_i(e)S_0\subset \oplus_i S_i=S$). Puis, remarquons que $p_i\circ p_j=0$ pour $i\neq j$, donc pour $\mathcal{M}$ un objet de $\mathop{\bf \Gamma\Phi M_S^h}$, nous avons les $p_i(\mathcal{M})$ en somme directe dans $\mathcal{M}$. Enfin, $p_i(e)\mathcal{M}^{\Gamma_f}\subset p_i(\mathcal{M})$, et $p_i(e)\mathcal{M}^{\Gamma_f}$ est isomorphe comme $S_0$-module à $\mathcal{M}^{\Gamma_f}$ car $\mathcal{M}$ est sans $p'$-torsion, et $p$ ne divise pas $p_i(e)$. Donc, nous avons que pour $\mathcal{M}$ un objet de $\mathop{\bf \Gamma\Phi M_S^h}$, $\mathcal{M}^{\Gamma_f}\otimes_{S_0}S=\oplus_i \mathcal{M}^{\Gamma_f}\otimes_{S_0}S_0p_i(e)$ s'injecte dans $\mathcal{M}$.

\begin{prop}
Soit $\mathcal{M}$ un objet de $\mathop{\bf \Gamma\Phi M_S^h}$, alors $\mathcal{M}^{\Gamma_f}$ est un objet de $\Gam$, et $\mathcal{M}=\mathcal{M}^{\Gamma_f}\otimes_{S_0}S$. De plus, $\mathcal{M}^{\Gamma_f}/\pi_0=\mathcal{M}/\pi$. Enfin, si $\mathcal{M}$ est $S$-libre, alors $\mathcal{M}^{\Gamma_f}$ est $S_0$-libre.
\end{prop}

\begin{proof}
Nous avons que $p_0(\mathcal{M})=\mathcal{M}^{\Gamma_f}$. Or, comme l'action de $\Gamma_f$ est triviale modulo $\pi$, nous avons que pour tout $x\in \mathcal{M}$, $x-p_0(x)\in\pi\mathcal{M}$, donc si $\mathcal{N}$ est le $S$-module engendré par $\mathcal{M}^{\Gamma_f}$ (c'est à dire que $\mathcal{N}=\mathcal{M}^{\Gamma_f}\otimes_{S_0}S$ d'après la remarque précédent la proposition), alors $\mathcal{M}=\mathcal{N}+\pi\mathcal{M}$, donc comme $\mathcal{M}$ est de type fini sur $S$, et que l'idéal engendré par $\pi$ est dans le radical de $S$, le lemme de Nakayama nous donne que $\mathcal{M}=\mathcal{N}$.

Puis, $\mathcal{M}$ est de type fini sur $S$, donc sur $S_0$ (car $S$ est un $S_0$-module libre de rang fini par le lemme \ref{lemmess0}), donc engendré sur $S_0$ par exemple par la famille finie $(m_i)$. Alors, $p_0(\mathcal{M})=\mathcal{M}^{\Gamma_f}$ est engendré par la famille $(p_0(m_i))$ (car $p_0$ est un morphisme de $S_0$-modules), donc est de type fini. De plus, $\mathcal{M}^{\Gamma_f}$ étant inclus dans $\mathcal{M}$, il est sans $p'$-torsion.

Ensuite, nous avons que $\pi\mathcal{M}\cap \mathcal{M}^{\Gamma_f}=\pi_0\mathcal{M}^{\Gamma_f}$ : pour $S$, l'égalité $\pi S\cap S_0=\pi_0S_0$ provient juste de ce que $\pi_0$ est un multiple de $\pi$, donc $\pi_0S_0\subset \pi S\cap S_0$, et pour la réciproque, que $S_0=W[[\pi_0]]$. Cela se traduit par la suite exacte de $S_0$-modules \[\xymatrix{
0 \ar[r] & \pi_0 S_0\ar[r] & S\ar[r]& S/\pi S\oplus S/S_0\ar[r]&0
}\]
(la surjectivité vient juste de ce que $S/\pi=W$, et que $W\subset S_0$), et en tensorisant par $\mathcal{M}^{\Gamma_f}$ au dessus de $S_0$, nous avons la suite exacte de $S_0$-modules \[\xymatrix{
0 \ar[r] & \pi_0\mathcal{M}^{\Gamma_f}\ar[r] & \mathcal{M}\ar[r]& \mathcal{M}/\pi \oplus \mathcal{M}/\mathcal{M}^{\Gamma_f}\ar[r]&0
}\]
ce qui traduit bien $\pi\mathcal{M}\cap \mathcal{M}^{\Gamma_f}=\pi_0\mathcal{M}^{\Gamma_f}$. Par conséquent, $\mathcal{M}^{\Gamma_f}/\pi_0$ s'injecte dans $\mathcal{M}/\pi$, et l'action de $\Gamma_0$ provient de celle sur $\mathcal{M}/\pi$, qui est triviale par définition. De plus, nous avons vu que pour tout $x\in \mathcal{M}$, $x-p_0(x)\in\pi\mathcal{M}$, donc comme $p_0(x)\in \mathcal{M}^{\Gamma_f}$, l'application naturelle $\mathcal{M}^{\Gamma_f}/\pi_0\to \mathcal{M}/\pi$ (dont nous avons vu l'injectivité) est surjective. Par conséquent, si $\mathcal{M}$ est $S$-libre, $\mathcal{M}^{\Gamma_f}/\pi_0$ est $W$-libre et $\mathcal{M}^{\Gamma_f}$ est sans $\pi_0$-torsion, et donc $S_0$ étant complet pour la topologie $\pi_0$-adique, une $W$-base de $\mathcal{M}^{\Gamma_f}/\pi_0$ se relève en une $S_0$-base de $\mathcal{M}^{\Gamma_f}$.

Enfin, $\phi$ commute à $\Gamma$, donc laisse stable $\mathcal{M}^{\Gamma_f}$, donc induit un morphisme $\Phi_0:\mathcal{M}^{\Gamma_f}\otimes_{\s}S_0\to \mathcal{M}^{\Gamma_f}$. Pour étudier le conoyau, remarquons d'abord que $x\otimes y\in S_0\otimes_{\s}S_0\mapsto \phi(x)y\in S_0$ et $x\otimes y\in S\otimes_{\s}S\mapsto \phi(x)y\in S$ sont des isomorphismes (préservant l'action naturelle de $\Gamma_f$), donc $S_0\otimes_{\s(S_0)}S_0[\frac{1}{q}]\simeq S_0[\frac{1}{q}]$ et $S\otimes_{\s(S)}S[\frac{1}{q}]\simeq S[\frac{1}{q}]$ (puisque $S[\frac{1}{q}]$ est plat sur $S$ et $S_0[\frac{1}{q}]$ est plat sur $S_0$). Par conséquent, $S_0\otimes_{\s(S_0)}S\simeq S\otimes_{\s(S)}S$ et $S_0\otimes_{\s(S_0)}S[\frac{1}{q}]\simeq S\otimes_{\s(S)}S[\frac{1}{q}]$ ; plus présicément, si $y_i\in S\otimes_{\s(S)}S$ s'envoye dans $S$ sur $p_i(e)$ (nous pouvons supposer que $y_0=1$ car $p_0(e)$ est inversible dans $S_0$), alors $ S\otimes_{\s(S)}S=\oplus_i S_0\otimes_{\s(S_0)}S_0y_i$ et $ S\otimes_{\s(S)}S[\frac{1}{q}]=\oplus_i S_0\otimes_{\s(S_0)}S_0[\frac{1}{q}]y_i$ (c'est bien le même $y_i$, car $ S\otimes_{\s(S)}S$ s'injecte dans $S\otimes_{\s(S)}S[\frac{1}{q}]$, puisque $ S\otimes_{\s(S)}S$ est sans $q$-torsion). Et l'action naturelle de $\Gamma_f$ sur $ S\otimes_{\s(S)}S[\frac{1}{q}]$ revient à dire que $g(y_i)=\X(g)^iy_i$ pour $g\in\Gamma_f$. Puisque $\mathcal{M}=\mathcal{M}^{\Gamma_f}\otimes_{S_0}S$, nous avons que $\mathcal{M}\otimes_{\s(S)}S[\frac{1}{q}]=\mathcal{M}^{\Gamma_f}\otimes_{\s(S_0)}S[\frac{1}{q}]=\oplus_i \mathcal{M}^{\Gamma_f}\otimes_{\s(S_0)}S_0[\frac{1}{q}]y_i$. Donc, $\mathcal{M}^{\Gamma_f}\otimes_{\s(S_0)}S_0[\frac{1}{q}]$ s'injecte naturellement dans $\mathcal{M}\otimes_{\s(S)}S[\frac{1}{q}]$, et $(\mathcal{M}\otimes_{\s(S)}S[\frac{1}{q}])^{\Gamma_f}=\mathcal{M}^{\Gamma_f}\otimes_{\s(S_0)}S_0[\frac{1}{q}]$.

Ensuite, $\Phi: \mathcal{M}\otimes_{\s}S\to \mathcal{M}$ est injective, de conoyau tué par $q^h$ (par définition), donc comme $S[\frac{1}{q}]$ est plat sur $S$, $\Phi$ induit une bijection $\Phi:\mathcal{M}\otimes_{\s(S)}S[\frac{1}{q}]\to \mathcal{M}\otimes_{S}S[\frac{1}{q}]$. Puis, $S[\frac{1}{q}]=\oplus_i S_0[\frac{1}{q}]p_i(e)$, donc $\mathcal{M}\otimes_{S}S[\frac{1}{q}]=\mathcal{M}^{\Gamma_f}\otimes_{S_0}S[\frac{1}{q}]=\oplus_i \mathcal{M}^{\Gamma_f}\otimes_{S_0}S_0[\frac{1}{q}]p_i(e)$, donc $\mathcal{M}^{\Gamma_f}\otimes_{S_0}S_0[\frac{1}{q}]$ s'injecte dans $\mathcal{M}\otimes_{S}S[\frac{1}{q}]$ et $\mathcal{M}^{\Gamma_f}\otimes_{S_0}S[\frac{1}{q}]=(\mathcal{M}\otimes_{S}S[\frac{1}{q}])^{\Gamma_f}$. Par conséquent, le diagramme

\[\xymatrix{ \mathcal{M}\otimes_{\s(S)}S[\frac{1}{q}]\ar[r]^-{\Phi}&\mathcal{M}\otimes_{S}S[\frac{1}{q}]\\
\mathcal{M}^{\Gamma_f}\otimes_{\s(S_0)}S_0[\frac{1}{q}]\ar[r]^-{\Phi_0}\ar@{^{(}-{>}}[u]^-{i}&\mathcal{M}^{\Gamma_f}\otimes_{S_0}S_0[\frac{1}{q}]\ar@{^{(}-{>}}[u]^-{j}
}\]

est commutatif, avec $\Phi$ bijective, $i$ et $j$ injective, et $\Phi$ (qui commute à l'action de $\Gamma_f$) qui identifie $(\mathcal{M}\otimes_{\s(S)}S[\frac{1}{q}])^{\Gamma_f}$ à $(\mathcal{M}\otimes_{S}S[\frac{1}{q}])^{\Gamma_f}$, donc $\Phi_0$ est bijective (donc $\mathcal{M}^{\Gamma_f}/\Phi_0(\mathcal{M}^{\Gamma_f}\otimes_{\s}S_0)$ est de $q$-torsion, donc tué par une puissance de $q$ car $\mathcal{M}^{\Gamma_f}$ est de type fini sur $S_0$).

Soit alors $x\in\mathcal{M}^{\Gamma_f}$. Par définition, il existe $y\in  \mathcal{M}\otimes_{\s(S)}S=\oplus_i\mathcal{M}^{\Gamma_f}\otimes_{\s(S_0)}S_0y_i$ tel que $\Phi(y)=q^hx$. La commutativité du diagramme et la bijectivité de $\Phi_0$ nous donne que $y\in \mathcal{M}^{\Gamma_f}\otimes_{\s(S_0)}S_0[\frac{1}{q}]$. Donc nous avons $y\in \big(\mathcal{M}^{\Gamma_f}\otimes_{\s(S_0)}S_0[\frac{1}{q}]\big)\cap\big(\oplus_i\mathcal{M}^{\Gamma_f}\otimes_{\s(S_0)}S_0y_i\big)=  \big(\mathcal{M}^{\Gamma_f}\otimes_{\s(S_0)}S_0[\frac{1}{q}]\big)\cap\big(\mathcal{M}^{\Gamma_f}\otimes_{\s(S_0)}S_0\big)=\mathcal{M}^{\Gamma_f}\otimes_{\s(S_0)}S_0$. En définitive, nous avons bien que $\mathcal{M}^{\Gamma_f}/\Phi_0(\mathcal{M}^{\Gamma_f}\otimes_{\s}S_0)$ est tué par $q^h$.

Finalement, nous avons bien que $\mathcal{M}^{\Gamma_f}$ est un objet de $\Gam$.
\end{proof}
\begin{rmq}
De la même façon que pour $S_i$, nous montrons pour $\mathcal{M}$ un objet de $\mathop{\bf \Gamma\Phi M_S^h}$ que $p_i(\mathcal{M})=\mathcal{M}^{\Gamma_f}\otimes_{S_0}S_i=p_i(e)\mathcal{M}^{\Gamma_f}$.
\end{rmq}

\begin{theo}\label{theoextensionscalaire}
L'extension des scalaires de $S_0$ à $S$ induit une équivalence de catégories entre $\Gam$ et $\mathop{\bf \Gamma\Phi M_S^h}$, préservant suites exactes et produit tensoriel (si ce dernier est encore dans la catégorie). Un quasi-inverse est donné par les points fixes par $\Gamma_f$.
\end{theo}
\begin{proof}
L'essentielle surjectivité se prouve en remarquant que si $f:\mathcal{M}\to \mathcal{N}$ est un morphisme de $\mathop{\bf \Gamma\Phi M_S^h}$, alors comme il commute à l'action de $\Gamma_f$, $f$ induit bien un morphisme de $(\phi,\Gamma_0)$-modules entre $\mathcal{M}^{\Gamma_f}$ et $\mathcal{N}^{\Gamma_f}$ (qui redonne $f$ en étendant les scalaires de $S_0$ à $S$). Le reste est immédiat à partir des résultats précédents.
\end{proof}

\subsection{Modules de Wach}\label{modulesdeWach}
L. Berger a défini dans \cite{Berger} le module de Wach $\mathop{\bf N}(T)$ d'un réseau $T$ d'une $\Q_p$-représentation cristalline $V$ à poids de Hodge-Tate négatifs comme l'unique $S$-sous-module de $\mathop{\bf D^+}(T):=(A_S^+\otimes_{\Z_p}T)^{H_{\K}}$ (avec $A_S^+=W(R)\cap \O_{\widehat{\E}_{nr}}$) vérifiant :
\begin{itemize}
\item $\mathop{\bf N}(T)$ est un $S$-module libre de rang la dimension de $V$ ;
\item l'action de $\Gamma$ préserve $\mathop{\bf N}(T)$ et est triviale sur $\mathop{\bf N}(T)/\pi \mathop{\bf N}(T)$ ;
\item il existe un entier $r\geq 0$ tel que $\pi^r \mathop{\bf D^+}(T)\subset \mathop{\bf N}(T)$.
\end{itemize}
Il définit de même le module de Wach $\mathop{\bf N}(V)$ d'une représentation cristalline $V$ à poids de Hodge-Tate négatifs. L'unicité donne en particulier que $\mathop{\bf N}$ va préserver somme directe et produit tensoriel, ce qui nous intéressera tout particulièrement.

Rappelons le Théorème 1' de \cite{Wa} : 

\begin{theore}
Si $N$ est un objet de $\MF$ avec $0\leq h\leq p-2$, alors $\V(N)$ est isomorphe (via l'application $g_N$) comme représentation galoisienne à $\mathop{\bf V_{\O_{\E}}}(\mathop{\rm F}^-(N))$. 
\end{theore}
Alors, nous avons
\begin{prop}\label{Wach}
Si $N$ est un objet de $\MF$ avec $0\leq h\leq p-2$, $\mathop{\bf D_{\O_{\E}}}(g_N)$ (qui identifie $\mathop{\rm F}^-(N)=N\otimes_W\O_{\E}$ à $\mathop{\bf D_{\O_{\E}}}(\V(N))$) envoie $N\otimes_W S$ sur $\mathop{\bf N}(\V(N))$ (le module de Wach associé à $\V(N)$).
\end{prop}

\begin{proof} En passant au dual, cela revient à dire que $\mathop{\rm F}(N^*)\otimes_{S_0}S$ est isomorphe à $\mathop{\bf N}(\V(N)^*)$ par fonctorialité du module de Wach envers le dual. Appelons $T=\V(N)^*$ et $r\leq p-2$ l'entier tel que $\Fil^r(N^*)\neq \{0\}$, $\Fil^{r+1}(N^*)=\{0\}$. Remarquons que la structure de $(\phi,\Gamma_0)$-module de $\mathop{\rm F}(N^*)$ induit une structure de $(\phi,\Gamma)$-module sur $N^*\otimes_W S$, et que $N^*\otimes_W \frac{1}{\pi^r}S$ est le dual (au sens généralisé des modules de Wach) d'un $(\phi,\Gamma)$-module de hauteur finie (puisque égale à $ r$) sur $S$, donc par le résultat de J.-M. Fontaine (cf \cite{F91}, p.296), les périodes de $N^*\otimes_W \frac{1}{\pi^r}S$ sont dans $A_S^+$. Par conséquent, $\mathop{\bf V_{\O_{\E}}}\big((N^*\otimes_W \frac{1}{\pi^r}S)\otimes_{S}\O_{\E}\big)=\mathop{\bf V_{\O_{\E}}}(\mathop{\rm F}(N^*)\otimes_{S_0}\O_{\E})=T=((N^*\otimes_W \frac{1}{\pi^r}S)\otimes_SA_S^+)^{\phi}\subset N^*\otimes_W \frac{1}{\pi^r}A_S^+$.

Puis, l'identification de $\mathcal{N}$ avec $\mathop{\bf D_{\O_{\E}}}(\mathop{\bf V_{\O_{\E}}}(\mathcal{N}))$ pour $\mathcal{N}$ un $(\phi,\Gamma)$-module sur $\O_{\E}$ est induite par la multiplication dans $\O_{\widehat{\E}_{nr}}$. Donc, comme $T\subset N^*\otimes_W \frac{1}{\pi^r}A_S^+$, nous avons $\mathop{\bf D^+}(T)\subset ((N^*\otimes_W \frac{1}{\pi^r}A_S^+)\otimes A_S^+)^{H_{\K}}$ qui est identifié à $(N^*\otimes_W \frac{1}{\pi^r}A_S^+)^{H_{\K}}=N^*\otimes_W \frac{1}{\pi^r}S$. Donc la dernière condition de la définition d'un module de Wach, $\pi^r \mathop{\bf D^+}(T)\subset N^*\otimes_W S$, est vérifiée.

\begin{lem}\label{uniciteberger}
Sous les notations précédentes, nous avons l'inclusion $\mathop{\bf N}(T)\subset N^*\otimes_W S$. 
\end{lem}
\begin{rmq}
La démonstration donnée ci-dessous est exactement l'idée principale de la démonstration de l'unicité du module de Wach (cf. proposition II.1.1 de \cite{Berger})
\end{rmq}
\begin{proof} Notons $\mathcal{N}_1=\mathop{\bf N}(T)$ et $\mathcal{N}_2=N^*\otimes_W S$. $\mathcal{N}_1\subset \mathop{\bf D^+}(T)$ par définition, donc nous avons l'inclusion $\pi^r\mathcal{N}_1\subset \mathcal{N}_2$. Soit $x\in \mathcal{N}_1$ et $s$ l'entier tel que $\pi^sx\in \mathcal{N}_2$, mais $\pi^sx\notin \pi \mathcal{N}_2$. Choisisons $x\notin \pi \mathcal{N}_1$ tel qu'en plus $s$ soit maximal, ce qui fait que $\pi^s\mathcal{N}_1\subset \mathcal{N}_2$. Comme $\pi^sx\in \mathcal{N}_2$ et que $\Gamma$ agit trivialement sur $\mathcal{N}_2/\pi\mathcal{N}_2$, nous avons pour tout $g\in\Gamma$ que $(g-1)(\pi^sx)\in \pi \mathcal{N}_2$, et nous pouvons écrire $(g-1)(\pi^sx)=g(\pi^s)(g(x)-x)+(g(\pi^s)-\pi^s)x$. Comme $\Gamma$ agit trivialement sur $\mathcal{N}_1/\pi\mathcal{N}_1$, et que $\pi^s\mathcal{N}_1\subset \mathcal{N}_2$, nous avons que $g(\pi^s)(g(x)-x)\in \pi \mathcal{N}_2$, et donc que $(g(\pi^s)-\pi^s)x\in \pi \mathcal{N}_2$, ce qui est une contradiction si $s\geq 1$, parce qu'alors $g(\pi^s)-\pi^s=(\X(g)^s-1)\pi^s+\cdots$. Donc nous avons bien $\mathcal{N}_1\subset \mathcal{N}_2$, autrement dit $\mathop{\bf N}(T)\subset N^*\otimes_W S$.
\end{proof}

L'étude du paragraphe précédent nous donne que le $S_0$-module $\mathcal{N}=\mathop{\bf N}(T)^{\Gamma_f}$ est libre et $\mathop{\bf N}(T)=\mathop{\bf N}(T)^{\Gamma_f}\otimes_{S_0}S$. Utilisons alors le fait que le foncteur $\mathop{\rm F}$ est essentiellement surjectif (à cause de l'hypothèse sur $h$) pour dire que $ \mathcal{N}$ est isomorphe en tant que $(\phi,\Gamma_0)$-module à $\mathop{\rm F}(N^*)$, donc $\mathop{\bf N}(T)$ est isomorphe au $(\phi,\Gamma)$-module $N^*\otimes_W S$. Notons $i$ cet isomorphisme.

Remarquons que $\mathop{\bf N}(T)\otimes_S\O_{\E}=\mathop{\bf D_{\O_{\E}}}(T)=N^*\otimes_W \O_{\E}$, car une représentation cristalline est de hauteur finie. Par conséquent, $i$ induit un isomorphisme de $\mathop{\bf D_{\O_{\E}}}(T)$ qui envoye $\mathop{\bf N}(T)$ sur $N^*\otimes_W S$, et comme il préserve $D^+(T)$, nous obtenons bien $N^*\otimes_WS\subset D^+(T)$, donc $N^*\otimes_WS=\mathop{\bf N}(T)$ car il vérifie toutes les conditions de la définition du module de Wach.
\end{proof}

Nous pouvons alors en déduire la proposition qui nous intéresse :

\begin{prop}
Soit $N_{i,j}$ des objets de $\mathop{\bf MF_W^{h}}$ avec $0\leq h\leq p-2$, $L$ un sous-objet (dans $\mathop{\bf MF_W^+}$) de $M=\displaystyle \bigoplus_i\otimes_j N_{i,j}$. Alors les isomorphismes de modules de Wach $${}^t\mathop{\bf D_{\O_{\E}}}(g_{N^*_{i,j}}): \mathop{\bf N}(\V(N_{i,j}^*)^*) \to N_{i,j}\otimes_{W} S $$ identifient $L\otimes_W S$ à un module de Wach.
\end{prop}

\begin{proof} Les isomorphismes ${}^t\mathop{\bf D_{\O_{\E}}}(g_{N^*_{i,j}})$ induisent un isomorphisme $$\displaystyle\bigoplus_{i=1}^n\otimes_{j=1}^{m_i}{}^t\mathop{\bf D_{\O_{\E}}}(g_{N^*_{i,j}}) : \mathop{\bf N}(\displaystyle\bigoplus_{i=1}^n\otimes_{j=1}^{m_i}\V(N^*_{i,j})^*)\to M\otimes_{W}S$$ (puisque le module de Wach préserve le produit tensoriel). Nous utiliserons cet isomorphisme pour identifier ces deux espaces.

Notons $(e_i)$ une base de $M$ telle que $(p^{\a_i}e_i)$ soit une base de $L$, avec $\a_i\in \N\cup\{+\infty\}$. Notons aussi $n=\rg_W(M)$.

La proposition \ref{propositiondur} affirme que $L\otimes_W S$ est stable par l'action de $\Gamma$. Considérons alors la sous-représentation galoisienne $T$ de $U_M:=\displaystyle\bigoplus_{i=1}^n\otimes_{j=1}^{m_i}\V(N^*_{i,j})^*$ définie par $T= \mathop{\bf V_{\O_{\E}}}(L\otimes_W \O_{\E})$. Montrons que $\mathop{\bf N}(T)=L\otimes_W S$, c'est-à-dire vérifions les conditions qui caractérisent un module de Wach :
\begin{itemize}
\item $L\otimes_{W}S\subset T\otimes_{\Z_p} \O_{\widehat\E_{nr}}\cap (U_M\otimes_{\Z_p}A_S^+)^{H_{\K}}=\mathop{\bf D^+}(T)$ : l'inclusion provient de ce que $T\otimes_{\Z_p} \O_{\widehat\E_{nr}}=L\otimes_{W} \O_{\widehat\E_{nr}}$ et $\mathop{\bf N}(U_M)=M\otimes_{W}S$ via l'isomorphisme (et donc $M\otimes_{W}S\subset D^+(U_M)=(U_M\otimes_{\Z_p}A_S^+)^{H_{\K}}$) ; l'égalité se montre en considérant les coordonnées suivant la base $(e_i)$, car si $x\in T\otimes_{\Z_p} \O_{\widehat\E_{nr}}\cap (U_M\otimes_{\Z_p}A_S^+)^{H_{\K}}$, alors il existe $(\b_i)\in (A_S^+)^n$ et $(\delta_i)\in \O_{\widehat\E_{nr}}^n$ avec $x=\sum_i \b_i e_i=\sum_i p^{\a_i}\delta_i e_i$, donc $\b_i=p^{\a_i}\delta_i$ pour tout $i$, donc $\b_i\in p^{\a_i}A_S^+$ pour tout $i$, ce qui donne $T\otimes_{\Z_p} \O_{\widehat\E_{nr}}\cap (U_M\otimes_{\Z_p}A_S^+)^{H_{\K}}\subset \mathop{\bf D^+}(T)$ (l'inclusion réciproque étant immédiate);
\item $L\otimes_{W}S$ est un $S$-module libre de rang égal à celui de $T$ sur $\Z_p$ (qui est celui de $L\otimes_W \O_{\E}$ sur $\O_{\E}$, donc celui de $L$ sur $W$) ;
\item l'action de $\Gamma$ laisse stable $L\otimes_{W}S$ (c'est la proposition \ref{propositiondur}) et est triviale modulo $\pi$ : l'action de $\Gamma$ sur $M\otimes_W S$ étant triviale modulo $\pi$ par construction, pour $\gamma\in \Gamma$, pour $i$ fixé, il existe $(x_j)\in S^{n-1}$ et $(y_j)\in S^n$ tels que $\gamma(e_i)=e_i+\pi\sum_{j\neq i}x_je_j$ et $\gamma(p^{\a_i}e_i)=\sum_jy_jp^{\a_j}e_j$ ; donc $y_i=1$ et $p^{\a_j}y_j=\pi x_jp^{\a_i}$ pour $j\neq i$, donc $\pi$ divise $y_j$ dans $S$ pour $j\neq i$. 
\item il existe $r$ un entier positif tel que $\pi^r \mathop{\bf D^+}(U_M)\subset M\otimes_{W}S$, donc ce $r$ donne $\pi^r \mathop{\bf D^+}(T)\subset  M\otimes_{W}S\cap L\otimes_W \O_{\widehat\E_{nr}}= L\otimes_{W}S$. En effet, si $x\in M\otimes_{W}S\cap L\otimes_W \O_{\widehat\E_{nr}}$, alors il existe $(\b_i)\in S^n$ et $(\delta_i)\in \O_{\widehat\E_{nr}}^n$ avec $x=\sum_i \b_i e_i=\sum_i p^{\a_i}\delta_i e_i$, donc $\b_i=p^{\a_i}\delta_i$ pour tout $i$, donc $\b_i\in p^{\a_i}S$ pour tout $i$, ce qui donne $M\otimes_{W}S\cap L\otimes_W \O_{\widehat\E_{nr}}\subset L\otimes_{W}S$ (l'inclusion réciproque étant immédiate).
\end{itemize}
D'où, nous avons bien $\mathop{\bf N}(T)=L\otimes_{W}S$.
\end{proof}

Ce qui nous intéressera tout particulièrement, c'est le corollaire suivant :
\begin{prop}\label{propositionmodulewach}
Soit $N_{i,j}$ des objets de $\mathop{\bf MF_W^{-h}}$ avec $0\leq h\leq p-2$, $L$ un sous-objet (dans $\mathop{\bf MF_W^-}$) facteur direct (comme $W$-module) de $M=\displaystyle \bigoplus_i\otimes_j N_{i,j}$. Alors les isomorphismes de modules de Wach $$\mathop{\bf D_{\O_{\E}}}(g_{N_{i,j}}) :N_{i,j}\otimes_{W} S\to \mathop{\bf N}(\V(N_{i,j})) $$ induisent un isomorphisme de module de Wach $$ L\otimes_{W}S\to \mathop{\bf N}(\mathop{\bf \overline{V}_{cris}}(L))$$ où $\mathop{\bf \overline{V}_{cris}}(L)$\index{$\mathop{\bf \overline{V}_{cris}}$}$:=\mathop{\bf V_{cris,p}}(D_L)\cap \displaystyle\bigoplus_{i=1}^n\otimes_{j=1}^{m_i} \V(N_{i,j})$.
\end{prop}
\begin{proof} Les isomorphismes $\mathop{\bf D_{\O_{\E}}}(g_{N_{i,j}})$ induisent un isomorphisme $$\mathop{\bf D_{\O_{\E}}}(g_M):=\displaystyle\bigoplus_{i=1}^n\otimes_{j=1}^{m_i}\mathop{\bf D_{\O_{\E}}}(g_{N_{i,j}}) : M\otimes_{W}S\to \mathop{\bf N}(\mathop{\bf \overline{V}_{cris}}(M))$$ Par dualité, il suffit de voir que si $L_0=M/L$, alors $\mathop{\bf D_{\O_{\E}}}(g_M)$ induit un isomorphisme de $\mathop{\bf \overline{V}_{cris}}(L_0^*)$ sur $L_0^*\otimes_WS$. Posons $T= \mathop{\bf V_{\O_{\E}}}(L_0^*\otimes_W \O_{\E})$. La proposition précédente nous donne bien que $\mathop{\bf N}(T)=L_0^*\otimes_{W}S$ (via $\mathop{\bf D_{\O_{\E}}}(g_M)$).

Puis, une propriété du module de Wach nous permet de conclure : $\mathop{\bf N}(T[\frac{1}{p}])/\pi$ s'identifie à $\mathop{\bf D_{cris,p}}(T\otimes_{\Z_p}\Q_p)$ (par le théorème III.4.4 de \cite{Berger}), et l'application $g_M$ envoie $\mathop{\bf N}(T)/\pi$ sur $L_0^*$, donc nous avons bien que $\mathop{\bf D_{cris,p}}(T\otimes_{\Z_p}\Q_p)=L_0^*\otimes_{W}\K$, donc que $T=\mathop{\bf \overline{V}_{cris}}(L_0^*)$. 
\end{proof}

\begin{rmq}\label{rmqquotient}
Si $M'$ est le quotient du $L$ considéré dans la proposition \ref{propositionmodulewach} par le sous-objet $L'$ (facteur direct comme $W$-module), alors $\mathop{\bf D_{\O_{\E}}}(g_M):L\otimes_W S\to \mathop{\bf N}(\mathop{\bf \overline{V}_{cris}}(L))$ (qui induit aussi un isomorphisme $L'\otimes_W S\to \mathop{\bf N}(\mathop{\bf \overline{V}_{cris}}(L'))$) induit par passage au quotient un isomorphisme $M' \otimes_W S\to \mathop{\bf N}(\mathop{\bf \overline{V}_{cris}}(M'))$. 
\end{rmq}

\begin{cor}\label{corollaireZink}
Soit $0\leq h\leq p-2$. Soit $M$ et $M'$ deux objets de $\mathop{\bf MF_W<h>}$, et $f:M\to M'$ un morphisme $\phi$-modules filtrés. Alors $\mathop{\bf V_{\E}}(\mathop{\rm F}(f)\otimes \Id_{\E})=\mathop{\bf V_{cris,p}}(f)$.
\end{cor}

\begin{proof}
Soient $N_{i,j}$ et $N'_{i,j}$ des objets de $\mathop{\bf MF_W^h}$, $L$ un sous-objet de $\oplus\otimes N_{i,j}$ et $L_0$ un sous-objet facteur direct de $L$, tel que $M=L/L_0$, $L'$ un sous-objet de $\oplus\otimes N'_{i,j}$ et $L'_0$ un sous-objet facteur direct de $L'$, tel que $M'=L'/L'_0$. Nous allons montrer que les isomorphismes $\oplus\otimes\mathop{\bf D_{\O_{\E}}}(g_{N_{i,j}})$ et $\oplus\otimes\mathop{\bf D_{\O_{\E}}}(g_{N'_{i,j}})$ identifient $f\otimes\Id$ à $\mathop{\bf D_{\E}}(\mathop{\bf V_{cris,p}}(f))$.

Pour cela, il suffit d'utiliser la fidélité et la pleine fidélité de $\mathop{\rm F}$ combiné au théorème \ref{theoextensionscalaire} (pour pouvoir dire que la réduction modulo $\pi$ est injective sur les morphismes de $(\phi,\Gamma)$-module entre $M\otimes_WS$ et $M'\otimes_WS$), plus le fait que $\mathop{\bf D_{\E}}(\mathop{\bf V_{cris,p}}(f))$ modulo $\pi$ redonne $f$ (d'après les résultats de Berger).
\end{proof}

\begin{cor}
Soit $M$ et $M'$ deux objets de $\mathop{\bf MF_W<h>}$, et $f:M\to M'$ un morphisme $\phi$-modules filtrés. Alors $\mathop{\bf V_{cris,p}}(f)$ envoye $\mathop{\bf V_{\O_{\E}}}(M\otimes_W\O_{\E})$ dans $\mathop{\bf V_{\O_{\E}}}(M'\otimes_W\O_{\E})$. En particulier, $\mathop{\bf \overline{V}_{cris}}$ devient un foncteur en posant $\mathop{\bf \overline{V}_{cris}}(f)=\mathop{\bf V_{cris,p}}(f)$.
\end{cor}
\begin{proof}
C'est une conséquence immédiate du corollaire précédent, et de ce que si $T$ et $T'$ sont deux $\Z_p$-représentations cristallines, alors un morphisme de $(\phi,\Gamma)$-modules $g:\mathop{\bf N}(T)\to\mathop{\bf N}(T')$ induit $\mathop{\bf V_{\O_{\E}}}(g):T=\mathop{\bf V_{\O_{\E}}}(\mathop{\bf N}(T)\otimes_S\O_{\E})\to T'=\mathop{\bf V_{\O_{\E}}}(\mathop{\bf N}(T')\otimes_S\O_{\E})$.
\end{proof}

\section{Fin de la démonstration du théorème \ref{propositionbesoin2}}
\begin{theo}\label{theoessentielsurj}
Pour $0\leq h\leq p-2$, le foncteur $\mathop{\rm F}$ de $\mathop{\bf MF_{W,tf}^h}$ vers $\Gam$ a pour image essentielle $\Gam$.
\end{theo}

Pour montrer ce résultat, nous allons utiliser le théorème \ref{Wach} qui nous dit que pour $N$ objet de $\mathop{\bf MF_{W}^h}$, $\mathop{\rm F}(N)\otimes_{S_0} S$ est le module de Wach de $\V(N^*)^*$. Commençons par montrer :
\begin{prop}
Soit $\mathcal{M}$ un objet de $\Gam$ (avec $0\leq h\leq p-2$) de $p$-torsion, et $T'=\mathop{\bf V_{\O_{\E}}}(\mathcal{M}\otimes_{S_0}\O_{\E})$ la $\Z_p$-représentation galoisienne correspondant au $(\phi,\Gamma)$-module sur $\O_{\E}$ obtenu à partir de $\mathcal{M}$. Alors il existe $T''\subset T$ deux $\Z_p$-représentations galoisiennes cristallines (c'est à dire que le module sous-jacent est libre sur $\Z_p$, et en rendant $p$ inversible nous avons une représentation cristalline) à poids de Hodge-Tate dans $[\![-h,0]\!]$ telles que $T'$ s'identifie au quotient de $T$ par $T''$.
\end{prop}
\begin{proof}
Le Théorème 1' de \cite{Wa} (et la proposition \ref{Wach}) donne que $T'=\Hom_{\Z_p}\big(\V\big(\Hom_{W}(i^*(\mathcal{M}),\varinjlim{W/p^n})\big),\varinjlim{\Z_p/p^n}\big)$. En notant $X^*$ le dual de Pontriaguine d'un module de torsion $X$, cela s'écrit plus simplement en $T'=\V\big((\mathcal{M}/\pi_0)^*\big)^*$. Puis, puisque $(\mathcal{M}/\pi_0)^*$ est un objet de la catégorie $\mathop{\bf MF_{W,tf}^{-h}}$, la proposition 1.6.3 de \cite{Wi} nous donne qu'il existe $M_1\in \MF$ et un épimorphisme $M_1\to (\mathcal{M}/\pi_0)^*$. Le foncteur $\V$ étant exact, il existe donc une $\Z_p$-repésentation cristalline $T_1$ ($T_1=\V(M_1)$) dont les poids de Hodge-Tate sont dans $[\![0,h]\!]$ et un épimorphisme $T_1\to \V\big((\mathcal{M}/\pi_0)^*\big)$. 

Comme $\mathcal{M}$ est supposé de $p$-torsion, $ \V\big((\mathcal{M}/\pi_0)^*\big)$ est de $p$-torsion et de type fini, donc il existe un entier $n$ tel que $p^n \V\big((\mathcal{M}/\pi_0)^*\big)=\{0\}$. Alors $T_1/p^n$ se surjecte toujours sur $\V\big((\mathcal{M}/\pi_0)^*\big)$, et en passant au dual de Pontriaguine, $T'$ s'injecte dans $\Hom_{\Z_p}\big(T_1/p^n,\varinjlim{\Z_p/p^n}\big)=\Hom_{\Z_p}(T_1,\Z_p)/p^n$ (car $T_1$ est un $\Z_p$-module libre). Si $f$ est la projection canonique $\Hom_{\Z_p}\big(T_1,\Z_p)\to \Hom_{\Z_p}\big(T_1,\Z_p)/p^n$, alors $T=f^{-1}(T')$ convient (et il suffit de prendre $T''$ égal au noyau de la projection $f|_{T}$).
\end{proof}

\begin{prop}\label{propquotient}
Soit $\mathcal{M}$ un objet de $\Gam$ (avec $0\leq h\leq p-2$) de $p$-torsion, $T'=\mathop{\bf V_{\O_{\E}}}(\mathcal{M}\otimes_{S_0}\O_{\E})$ et $T''\subset T$ les représentations données par la proposition ci-dessus. Alors $\mathcal{M}\otimes_{S_0}S$ s'identifie à $\mathop{\bf N}(T)/\mathop{\bf N}(T'')$.
\end{prop}
\begin{proof}
Notons $\mathcal{M}_1= \mathcal{M}\otimes_{S_0}S$ et $\mathcal{M}_2= \mathop{\bf N}(T)/\mathop{\bf N}(T'')$ (tous les deux vus dans le $(\phi,\Gamma)$-module $\mathcal{M}\otimes_{S_0}\O_{\E}$, car $\mathop{\bf N}(T)\cap\mathop{\bf D_{\O_{\E}}}(T")=\mathop{\bf N}(T'')$ :  en effet, notons $\mathcal{N}=\mathop{\bf N}(T)\cap\mathop{\bf D_{\O_{\E}}}(T")=\mathop{\bf N}(T)\cap T''\otimes_{\Z_p}\O_{\hat{\E}_{nr}}$ qui est stable par l'action de $\Gamma$, nous avons que $\mathcal{N}\cap \pi\mathop{\bf N}(T)=\pi\mathcal{N}$ puisque $\pi$ est inversible dans $\O_{\hat{\E}_{nr}}$, donc $\mathcal{N}/\pi$ s'injecte dans $\mathop{\bf N}(T)/\pi$, donc $\Gamma$ agit bien trivialement sur $\mathcal{N}$. Puis, $T''\otimes_{\Z_p}\O_{\hat{\E}_{nr}}\cap (T\otimes_{\Z_p}A_S^+)^{H_{\K}}=\mathop{\bf D^+}(T'')$, car si $(e_i)$ est une base de $T$ telle que $(p^{\a_i}e_i)$ est une base de $T''$ (avec $\a_i\in\N\cup\{+\infty\}$), alors un élément $x$ de l'intersection s'écrit $x=\sum_ix_ie_i=\sum_ip^{\a_i}y_ie_i$ avec $x_i\in A_S^+$ et $y_i\in\O_{\hat{\E}_{nr}}$ ; donc $y_i\in p^{-a_i}A_S^+\cap\O_{\hat{\E}_{nr}}=A_S^+$ si $\a_i\neq +\infty$, et $\{0\}$ sinon, donc $x\in T''\otimes_{\Z_p}A_S^+$ et est fixé par $H_{\K}$, donc $T''\otimes_{\Z_p}\O_{\hat{\E}_{nr}}\cap (T\otimes_{\Z_p}A_S^+)^{H_{\K}}\subset \mathop{\bf D^+}(T'')$ (l'inclusion réciproque étant immédiate). Donc nous avons $\mathcal{N}\subset \mathop{\bf D^+}(T'')$ puisque $\mathop{\bf N}(T)\subset (T\otimes_{\Z_p}A_S^+)^{H_{\K}}= \mathop{\bf D^+}(T)$. Enfin, $\pi^h\mathop{\bf D^+}(T)\subset \mathop{\bf N}(T)$, donc $\pi^h\mathop{\bf D^+}(T'')\subset  \mathop{\bf N}(T)$, et comme $\pi^h\mathop{\bf D^+}(T'')\subset T''\otimes_{\Z_p} \O_{\hat{\E}_{nr}}$, nous avons bien que $\pi^h\mathop{\bf D^+}(T'')\subset \mathcal{N}$. Ces conditions caractérisent le module de Wach de $T''$, donc $\mathop{\bf N}(T)\cap\mathop{\bf D_{\O_{\E}}}(T")=\mathop{\bf N}(T'')$).

D'après les résultats p.296 de \cite{F91} (l'égalité entre $\mathop{\bf D_S^*}$ et $j_*\circ \mathop{\bf D_{\E}^*}$) (ou bien le lemme III.5 de \cite{Colmez}), nous avons $\mathcal{M}_1\subset \mathop{\bf D^+}(T')$ et $\mathcal{M}_2\subset \mathop{\bf D^+}(T')$, puisque tout deux sont des $S$-modules de type fini stables par $\phi$ et $p$-étales (puisque de $q$-hauteur finie). Puis, l'action de $\Gamma$ est triviale modulo $\pi$ dans les deux cas (puisque c'est le cas par définition sur $\mathop{\bf N}(T)$, et que l'action de $\Gamma_0$ est triviale modulo $\pi_0$ sur $\mathcal{M}$).

D'après le Théorème III.3.1 de \cite{Berger}, nous avons l'inclusion $\pi^h T\otimes_{\Z_p}A_{S}^+\subset  \mathop{\bf N}(T)\otimes_SA_S^+$. Par conséquent, en projetant nous obtenons que $\pi^h T'\otimes_{\Z_p}A_{S}^+\subset  \mathcal{M}_2\otimes_SA_S^+$. Par définition, nous avons que $\mathop{\bf D^+}(T')\subset  \mathop{\bf D^+}(T')\otimes_SA_S^+\subset T'\otimes_{\Z_p}A_S^+$, donc en prenant les points fixes sous l'action de $H_{\K}$, nous avons $\mathop{\bf D^+}(T')\subset \big(\mathop{\bf D^+}(T')\otimes_SA_S^+\big)^{H_{\K}}\subset \big(T'\otimes_{\Z_p}A_S^+\big)^{H_{\K}}=\mathop{\bf D^+}(T')$. Donc, en prenant les points fixes sous $H_{\K}$ dans l'inclusion $\pi^h T'\otimes_{\Z_p}A_{S}^+\subset  \mathcal{M}_2\otimes_SA_S^+$, nous obtenons que $\pi^h  \mathop{\bf D^+}(T')\subset  \big(\mathcal{M}_2\otimes_SA_S^+\big)^{H_{\K}}$. Donc nous avons $\pi^h  \mathop{\bf D^+}(T')\subset  \mathcal{M}_2$ en vertu du lemme :
\begin{lem}
Soit $\mathcal{N}$ un $S$-module de type fini sans $p'$-torsion, alors $\big(\mathcal{N}\otimes_SA_S^+\big)^{H_{\K}}=\mathcal{N}$.
\end{lem}
\begin{proof}
C'est une conséquence de la proposition 1.2.7 de \cite{F91}, qui nous donne (sous les hypothèses du lemme) une filtration décroissante $\mathcal{N}_i$ de $\mathcal{N}$, telle que $\mathcal{N}_i/\mathcal{N}_{i+1}$ est soit $S/p$-libre, soit $S$-libre. La propriété cherchée est stable par suite exacte, c'est à dire vérifie que si $0\to \mathcal{N''}\to \mathcal{N}\to \mathcal{N'}\to 0$ est une suite exacte de $S$-modules, et que $\big(\mathcal{N''}\otimes_SA_S^+\big)^{H_{\K}}=\mathcal{N''}$, $\big(\mathcal{N'}\otimes_SA_S^+\big)^{H_{\K}}=\mathcal{N'}$, alors $\big(\mathcal{N}\otimes_SA_S^+\big)^{H_{\K}}=\mathcal{N}$. Donc il suffit de montrer le lemme pour $\mathcal{N}$ qui est $S$-libre ou $S/p$-libre, ce qui provient de ce que $(A_S^+)^{H_{\K}}=S$ et $(A_S^+/p)^{H_{\K}}=S/p$.
\end{proof}

Puis $\frac{1}{\pi^h}\mathcal{M}_1$ est le dual (de Pontriaguine) d'un $(\phi,\Gamma)$-module sur $S$ de hauteur inférieure ou égale à $h$, sans $p'$-torsion, donc $T'=\mathop{\bf V_{\O_{\E}}}(\frac{1}{\pi^h}\mathcal{M}\otimes_{S_0}\O_{\E})$ vérifie $T'=(\frac{1}{\pi^h}\mathcal{M}\otimes_{S_0}A_S^+)^{\phi}$ (cf \cite{F91}, p.296) puisque $0\leq h$. Donc $T'\otimes_{\Z_p}A_S^+\subset \frac{1}{\pi^h}\mathcal{M}\otimes_{S_0}A_S^+$, et en prenant les points fixes sous $H_{\K} $ (et par le lemme précédent), nous obtenons $ \mathop{\bf D^+}(T')\subset  \frac{1}{\pi^h}\mathcal{M}_1$, donc $\pi^h \mathop{\bf D^+}(T')\subset\mathcal{M}_1$.

Ces conditions impliquent que la démonstration du lemme \ref{uniciteberger} s'applique ici (car $h\leq p-2$, pour que nous ayons si $0\leq s\leq h$, $\X(g)^s-1$ inversible dans $\Z_p$ (c'est à dire $\X(g)^s-1\neq 0$ modulo $p$) pour un $g\in\Gamma$), et donc $\mathcal{M}_1=\mathcal{M}_2$. 
\end{proof}
\begin{rmq}
L'unicité d'un tel module n'est plus vrai en général :  dans $S/pS$, $S/pS$ et $\pi^{p-1}S/pS$ sont deux $S$-modules de type fini, avec action de $\Gamma$ triviale modulo $\pi$, et si $T=\mathop{\bf V_{O_{\E}}}(\O_{\E}/p)$ (c'est à dire $\F_p$ avec l'action triviale), alors $\mathop{\bf D^+}(T)=S/pS$, donc la dernière condition est aussi vérifiée.
\end{rmq}

Il ne reste donc plus qu'à passer d'un module sur $S$ à un module sur $S_0$, ce qui est donné par le lemme suivant (qui est une conséquence immédiate de l'égalité $S=\bigoplus\limits_{0\leq i\leq p-2}S_{i}$) :
\begin{lem}\label{lemmepi}
Soit $\mathcal{M}$ un $S_0$-module, alors $\big(\mathcal{M}\otimes_{S_0}S\big)^{\Gamma_f}=\mathcal{M}$.
\end{lem}

Ces propositions et ces lemmes mis bout à bout nous donnent le théorème dans le cas d'un objet de $\Gam$ de $p$-torsion. C'est à dire que si $\mathcal{M}$ est un objet de $\Gam$ (avec $0\leq h\leq p-2$) de $p$-torsion, alors il existe $M$ un objet de $\mathop{\bf MF_{W,tf}^{h}}$ tel que $\mathcal{M}=\mathop{\rm F}(M)$. Et plus précisément, nous avons $M=i^*(\mathcal{M})=\mathcal{M}/\pi$. Donc, dans le cas où $\mathcal{M}$ n'est pas supposé de $p$-torsion, nous avons que $\mathcal{M}/p^n=\mathop{\rm F}(i^*(\mathcal{M}/p^n))$ pour tout $n$, donc en passant à la limite projective, nous obtenons bien que $\mathcal{M}=\mathop{\rm F}(i^*(\mathcal{M}))$, ce qui donne bien l'essentielle surjectivité de $\mathop{\rm F}$, et donc termine la démonstration du théorème \ref{theoessentielsurj}.

\section{Le point du torseur}
\subsection{Conséquence des théorèmes précédents}
Pour tout objet $N$ de $\MF$ avec $0\leq h\leq p-2$, construisons $\f_N$ comme la composée :

\[\xymatrix{ {\V(N)\otimes_{\Z_p} \O_{\widehat{\E}_{nr}}}  \ar[r]^-{g_N^{-1}}&{\mathop{\bf V_{ \O_{\widehat{\E}}  }}(\mathop{\rm F^-}(N))\otimes_{\Z_p} \O_{\widehat{\E}_{nr}}}  \ar[r]^-{\psi_{\mathop{\rm F^-}(N)}}& {\mathop{F^-}(N)\otimes_{\O_{\E}}  \O_{\widehat{\E}_{nr}}} \ar@{=}[d]\\&& N\otimes_{W}  \O_{\widehat{\E}_{nr}} 
}\]
où $\psi$ est l'isomorphisme de Fontaine (cf. paragaraphe \ref{isomfontaine}), $g_N$ l'isomorphisme de N. Wach (cf. paragraphe \ref{construction}), et $\mathop{\rm F^-}$ est le foncteur construit à la fin de la partie \ref{fonc}.

De la proposition \ref{propositionmodulewach} nous déduisons (toujours à $0\leq h\leq p-2$ fixé) :
\begin{prop}
Pour tout uplet d'objets $(N_{i,j})_{1\leq j\leq n,\,1\leq i\leq n_j}$ de $\MF$, pour tout sous-$\Phi$-module filtré $L$ facteur direct (comme $W$-module) de $\displaystyle\bigoplus_{j=1}^n\otimes_{i=1}^{n_j} N_{i,j} $, l'application $\bigoplus\otimes \f_{N_{i,j}}$ envoie $\mathop{\bf \overline{V}_{cris}}(L)\otimes_{\Z_p}\O_{\widehat{\E}_{nr}}$ bijectivement sur $L\otimes_{W}  \O_{\widehat{\E}_{nr}}$.
\end{prop}
\begin{proof}
Rappelons que $\mathop{\bf \overline{V}_{cris}}(L)= \mathop{\bf V_{cris,p}}(D_L)\cap \displaystyle\bigoplus_{j=1}^n\otimes_{i=1}^{n_j} \V(N_{i,j})$. Comme corollaire de la proposition \ref{propositionmodulewach}, l'inverse de la fonction $\psi_{\mathop{\bf D_{\O_{\E}}}(\V(N))}^{-1}\circ \big(\mathop{\bf D_{\O_{\E}}}(g_N)\otimes\Id\big)$ vérifie la propriété recherchée, car $\mathop{\bf D_{\O_{\E}}}(g_N)$ envoie $L\otimes_W S$ sur $ \mathop{\bf N}(\mathop{\bf \overline{V}_{cris}}(L))$, donc $L\otimes_W\O_{\E}$ sur $\mathop{\bf D_{\O_{\E}}}(\mathop{\bf \overline{V}_{cris}}(L))$. Il suffit alors de remarquer que $\f_N=\psi_{\mathop{\rm F}^-(N)}\circ \big(g_N^{-1}\otimes \Id\big)=\big(\mathop{\bf D_{\O_{\E}}}(g_N^{-1})\otimes\Id\big)\circ \psi_{\mathop{\bf D_{\O_{\E}}}(\V(N))}$ par commutativité du diagramme suivant :
\[\xymatrix{
\V(N)\otimes_{\Z_p}\O_{\widehat{\E}_{nr}}\ar[rr]^-{g_N^{-1}\otimes\Id} \ar[d]^-{\psi_{\mathop{\bf D_{\O_{\E}}}(\V(N))}}&&\mathop{\bf V_{\O_{\E}}}(\mathop{\rm F}^-(N))\otimes_{\Z_p}\O_{\widehat{\E}_{nr}}\ar[d]^-{\psi_{\mathop{\rm F}^-(N)}}\\\mathop{\bf D_{\O_{\E}}}(\V(N))\otimes_{\O_{\E}}\O_{\widehat{\E}_{nr}}\ar[rr]_-{\mathop{\bf D_{\O_{\E}}}(g_N^{-1})\otimes\Id} &&\mathop{\rm F}^-(N)\otimes_{\O_{\E}}\O_{\widehat{\E}_{nr}}
}\]
car $\mathop{\bf D_{\O_{\E}}}$ et $\mathop{\bf V_{\O_{\E}}}$ sont des foncteurs quasi-inverses l'un de l'autre.
\end{proof}

En combinant ce résultat et celui de la remarque \ref{gggg}, nous obtenons
\begin{prop}
Si $L$ est un sous-objet dans $\MF$ de $\displaystyle \bigoplus_i\otimes_j N_{i,j}$ avec $N_{i,j}$ des objets de $\mathop{\bf MF_W^{-h}}$ avec $0\leq h\leq p-2$, alors $\V(L)\subset \displaystyle \bigoplus_i\otimes_j \V(N)_{i,j}$.
\end{prop}
\begin{rmq}
Cette propriété peut être montré directement, en utilisant les propriétés des périodes des Lubin-Tate (qui donnent par produit les périodes des modules élémentaires) et le fait qu'un $\Phi$-module filtré simple est élémentaire, donc que (par Jordan-Hölder) tout $\Phi$-module filtré tué par $p$ a une filtration dont le gradué associé est somme directe de modules élémentaires. 
\end{rmq}

Combiné avec le théorème \ref{theowach}, et en introduisant $\mathcal{F}_1$ et $\mathcal{F}_2$ les foncteurs exacts de la catégorie $\mathop{\bf MF_W^{-h}}$ vers la catégorie des $ \O_{\widehat{\E}_{nr}}$-modules libres de rang fini, définis par : si $M$ objet de $\mathop{\bf MF_W^{-h}}$, $\mathcal{F}_1(M)=\V(M)\otimes_{\Z_p} \O_{\widehat{\E}_{nr}} $ et $\mathcal{F}_2(M)=M\otimes_{W}  \O_{\widehat{\E}_{nr}}  $, nous obtenons :
\begin{theo}\label{theoxxx}
Pour $0\leq h\leq p-2$ fixé, il existe $\f$ un isomorphisme de foncteur entre $\mathcal{F}_1$ et $\mathcal{F}_2$. 
De plus, vis à vis du produit tensoriel, l'isomorphisme peut être choisi de telle sorte que :
\begin{itemize}
\item pour tous objets $M$ et $N$ de $\MF$ tels que $M\otimes N$ est encore un objet de $\MF$, alors  le diagramme suivant est commutatif :\\
\end{itemize}
\[\xymatrix{
 {\V(N\otimes M)\otimes \O_{\widehat\E_{nr}}}\ar[d]  \ar[r]^-{\f_{N\otimes M}}   & (N\otimes M)\otimes \O_{\widehat\E_{nr}} \ar[d]\\
{(\V(M)\otimes \O_{\widehat\E_{nr}})\otimes (\V(N)\otimes\O_{\widehat\E_{nr}}}) \ar[r]_-{\f_M\otimes \f_N}   & (N\otimes \O_{\widehat\E_{nr}}) \otimes (M\otimes \O_{\widehat\E_{nr}}) 
}\]
\begin{itemize}
\item pour tout uplet d'objets $(N_{i,j})_{1\leq j\leq n,\,1\leq i\leq n_j}$ de $\MF$, pour tout sous-objet $L$ de $\displaystyle\bigoplus_{j=1}^n\otimes_{i=1}^{n_j} N_{i,j} $ dans $\MF$, l'application $\bigoplus\otimes \f_{N_{i,j}}$ restreinte à $\V(L)$ est égale à $\f_L$ ;
\item pour tout uplet d'objets $(N_{i,j})_{1\leq j\leq n,\,1\leq i\leq n_j}$ de $\MF$, pour tout sous-$\Phi$-module filtré $L$ facteur direct (comme $W$-module) de $\displaystyle\bigoplus_{j=1}^n\otimes_{i=1}^{n_j} N_{i,j} $, l'application $\bigoplus\otimes \f_{N_{i,j}}$ envoie $\big(\mathop{\bf V_{cris,p}}(D_L)\cap \displaystyle\bigoplus_{j=1}^n\otimes_{i=1}^{n_j} \V(N_{i,j})\big)\otimes_{\Z_p}\O_{\widehat{\E}_{nr}} $ bijectivement sur $L\otimes_{W}  \O_{\widehat{\E}_{nr}}$.
\end{itemize}
\end{theo}

Il faut juste regarder le comportement de $\f$ vis à vis du dual pour pouvoir déduire du théorème \ref{theoxxx} le théorème \ref{theotorseur}

\subsection{$\f$ et le dual}

Pour la suite, nous aurons besoin de définir $\f_{N}$ pour $N$ ayant des poids à la fois positifs et négatifs (par exemple si $N=\End(M)$ pour $M$ un objet de $\MF$). 

Appelons $W[-h]$ l'objet de $\MF$ dont le $W$-module sous-jacent est $W$, avec $\Fil^i(W[-h])=\left\{\begin{array}{c} W \text{ si }i\leq -h\\ 0\text{ si }i>-h\end{array}\right.$ et $\p^{-h}(x)=\s(x)$. Rappelons que $\Z_p(h)$ est la représentation galoisienne $\Z_p(1)^{\otimes h}$ pour $h\geq 0$ (et $\Z_p(h)=\Z_p(-h)^*$ si $h\leq 0$), et que $ \O_{\widehat\E_{nr}}(h)= \O_{\widehat\E_{nr}}\otimes_{\Z_p}\Z_p(h)$. Nous pouvons alors définir :

\begin{defi}
Supposons $0\leq h\leq \frac{p-2}{2}$. Posons $\mathop{\bf \widetilde V_{cris} }(N)=\V(N\otimes W[-h])\otimes_{\Z_p}\Z_p(-h) $\index{$\mathop{\bf \widetilde V_{cris} }$} pour $N$ objet de $\Mf$, et $\mathop{\bf \widetilde V_{cris} }(f)=\mathop{\bf V_{cris,p} }(f)$ restreinte à $\mathop{\bf \widetilde V_{cris} }(N)$ pour $f:N\to N'$ flèche de $\Mf$. Pour tout objet $N$ de $\Mf$ nous pouvons définir $\tilde\f_N$\index{$\tilde\f_N$} de la façon suivante : remarquons que $\mathop{\bf \widetilde V_{cris} }(N)\otimes \O_{\widehat\E_{nr}}= \big( \V(N\otimes W[-h])\otimes \Z_p(-h)\big)\otimes_{\Z_p} \O_{\widehat\E_{nr}}=  \big( \V(N\otimes W[-h])\otimes  \O_{\widehat\E_{nr}}\big)\otimes_{ \O_{\widehat\E_{nr}}}  \O_{\widehat\E_{nr}}(-h)$, et notons $\f_{h}={}^t\f_{ \O_{\widehat\E_{nr}}[-h]}^{-1}$, alors  $\tilde\f_N$ est l'isomorphisme
\[\xymatrix{  {\mathop{\bf \widetilde V_{cris} }(N)\otimes \O_{\widehat\E_{nr}}} \ar[rrr]^-{\f_{N\otimes W[-h]}\otimes \f_{h}} &&& \big((N\otimes W[-h])\otimes \O_{\widehat\E_{nr}}\big) \otimes  \O_{\widehat\E_{nr}}[h]  \ar[d]\\&&&  {N\otimes_W \O_{\widehat\E_{nr}} }}\] 
\end{defi}

Nous avons bien que $\mathop{\bf \widetilde V_{cris} }(N)\simeq \V(N)\otimes_{\Z_p} \V(W[-h])\otimes_{\Z_p}\Z_p(-h)\simeq \V(N)$ de manière naturelle pour $N$ un objet de $\MF$, car $N\otimes W[-h]$ est un objet de $\mathop{\bf MF_W^{-2h}}$, et comme $2h\leq p-2$, nous pouvons appliquer la remarque \ref{pcp}.

Un quasi-inverse de $\mathop{\bf \widetilde V_{cris} }$ est donné par $\mathop{\bf \widetilde D_{cris} }(V)=\D(V\otimes \Z_p(h))\otimes_W W[h]$. Une façon naturelle de voir $\mathop{\bf \widetilde V_{cris} }(N)$ dans $N\otimes_W A_{cris}$ est de dire que $$\mathop{\bf \widetilde V_{cris} }(N)=\Fil^0(N\otimes_W t^{-h}A_{cris})^{\Phi}t^{-h}$$

Avant de continuer de regarder les propriétés de $\mathop{\bf \widetilde V_{cris} }$, introduisons $\tilde g_N$ pour $N$ un objet de $\Mf$ si $0\leq h\leq \frac{p-2}{2}$ de la même façon que précédemment. De la proposition \ref{propositionmodulewach} et de la remarque \ref{rmqquotient}, nous déduisons :

\begin{cor}\label{co}
Pour tout uplet d'objets $(N_{i,j})_{1\leq j\leq n,\,1\leq i\leq m}$ de $\Mf$ avec $0\leq h\leq \frac{p-2}{2}$, pour tout sous-objet $L$ facteur direct (comme $W$-module) de $\displaystyle\bigoplus_{j=1}^n\otimes_{i=1}^{m} N_{i,j}$ et pour tout quotient $M$ de $L$, les applications $\tilde g_{N_{i,j}}$ induisent un isomorphisme de représentations de $\Gamma_{\K}$, de $\mathop{\bf V_{ \O_{\widehat{\E}}  }}(\mathop{\rm F}^-(M\otimes W[-mh]))\otimes \Z_p(-mh)$ sur un réseau de $\mathop{\bf V_{cris,p}}(D_M)$ (qui est l'image de $\mathop{\bf V_{cris,p}}(D_L)\cap \displaystyle\bigoplus_{j=1}^n\otimes_{i=1}^{n_j} \mathop{\bf \widetilde V_{cris}}(N_{i,j})$ par l'application projection).
\end{cor} 

Puis, pour l'étude vis à vis du dual, montrons d'abord le lemme suivant :
\begin{lem}\label{lemmedual}
Pour tout objet $N$ de $\Mf$ avec $0\leq h\leq \frac{p-2}{2}$, l'application surjective naturelle 
\[\xymatrix{ N\otimes N^*\ar[r]^-{\pi} &W}\]
induit un isomorphisme $$\mathop{\bf \widetilde V_{cris} }(N^*)\simeq \mathop{\bf \widetilde V_{cris} }(N)^*$$ dont le crochet de dualité correspond à $\mathop{\bf V_{cris,p}}(\pi)$.
\end{lem}

\begin{proof}
Soient $l\in\N$ et $l'\in\N$ tels que $N\otimes W[-l]$ et $N^*\otimes W[-l']$ soient des objets de $\mathop{\bf MF_W^{-2h}}$. Puis, $W$ désigne l'objet trivial de $\Mf$, et donc $\pi$ induit une application $\mathop{\rm F}^-(\pi):\mathop{\rm F}^-(N\otimes W[-l]))\otimes_{\O_{\E}}\mathop{\rm F}^-(N^*\otimes W[-l'])\to  \O_{\E}[-(l+l')]$, dont l'application linéaire sous-jacente est toujours celle obtenue par le crochet de dualité (c'est juste $\pi\otimes \Id$), donc induit un isomorphisme entre le dual de $\mathop{\rm F}^-(N\otimes W[-l]))\otimes_{\O_{\E}}\O_{\E}e_{l}$ et $\mathop{\rm F}^-(N^*\otimes W[-l']))\otimes_{\O_{\E}}\O_{\E}e_{l'}$ (où $\phi(e_{r})=q^{r}e_{r}$ et $g(e_{r})=\frac{\X(g)^{-r}\pi^{-r}}{g(\pi^{-r})}$ pour $g\in\Gamma$). Cet isomorphisme de $\O_{\E}$-modules est en fait un isomorphisme de $(\phi,\Gamma)$-module, car $\mathop{\rm F}^-(\pi)$ est un morphisme de $(\phi,\Gamma)$-modules. Comme $\mathop{\bf V_{\O_{\E}}}$ préserve le dual, en notant $\mathop{\bf \widetilde{V}_{\O_{\E}}}(N)=\mathop{\bf V_{\O_{\E}}}(\mathop{\rm F}^-(N\otimes W[-l]))\otimes_{\Z_p} \Z_p(-l)=\mathop{\bf V_{\O_{\E}}}\big(\mathop{\rm F}^-(N\otimes W[-l])\otimes_{\O_{\E}}\O_{\E}e_{l}\big)$ et $\mathop{\bf \widetilde{V}_{\O_{\E}}}(N^*)=\mathop{\bf V_{\O_{\E}}}(\mathop{\rm F}^-(N\otimes W[-l']))\otimes_{\Z_p} \Z_p(-l')=\mathop{\bf V_{\O_{\E}}}\big(\mathop{\rm F}^-(N^*\otimes W[-l'])\otimes_{\O_{\E}}\O_{\E}e_{l'}\big)$, nous avons que l'application $\mathop{\bf \widetilde V_{\O_{\E}}}(\mathop{\rm F}^-(\pi))$ identifie le dual de $\mathop{\bf \widetilde{V}_{\O_{\E}}}(N)$ (comme représentation de $\Gamma_{\K}$) à $\mathop{\bf \widetilde{V}_{\O_{\E}}}(N)^*$. 

Or, $\mathop{\bf \widetilde{V}_{\O_{\E}}}(N)$ est isomorphe à $\mathop{\bf \widetilde V_{cris} }(N)$ (via $\tilde{g}_{N}$) et $\mathop{\bf \widetilde{V}_{\O_{\E}}}(N^*)$ est isomorphe à $\mathop{\bf \widetilde V_{cris} }(N^*)$ (via $\tilde{g}_{N^*}$). Pour conclure, il suffit d'invoquer la commutativité du diagramme suivant (d'après le corollaire \ref{corollaireZink}) :
\[\xymatrix{
\mathop{\bf \widetilde{V}_{\O_{\E}}}(N)\otimes\mathop{\bf \widetilde{V}_{\O_{\E}}}(N^*) \ar[d]_-{\mathop{\bf \widetilde V_{\O_{\E}}}(\mathop{\rm F}^-(\pi))}\ar[rr]_-{\tilde{g}_{N}\otimes\tilde{g}_{N^*}}&&\mathop{\bf \widetilde{V}_{cris}}(N)\otimes \mathop{\bf \widetilde{V}_{cris}}(N^*)\ar[d]_-{\mathop{\bf V_{cris,p}}(\pi)} \\
\mathop{\bf V_{\O_{\E}}}(\mathop{\rm F}^-(W))\ar[rr]_-{\tilde{g}_{W}}& &\mathop{\bf \widetilde{V}_{cris}}(W)
}\]
\end{proof}

Nous en déduisons alors :
\begin{lem}
Sous les conditions du lemme \ref{lemmedual}, l'application $\tilde\f_N$ est fonctorielle vis à vis du dual, c'est-à-dire que le diagramme suivant commute :
\[\xymatrix{
 { \mathop{\bf \widetilde V_{cris} }(N)^* \otimes_{\Z_p} \O_{\widehat\E_{nr}}}\ar[d]  \ar[r]^-{({}^t\tilde\f_{N})^{-1}}   & N^*\otimes_W \O_{\widehat\E_{nr}}  \ar[d]\\
{\mathop{\bf \widetilde V_{cris} }(N^*)\otimes_{\Z_p} \O_{\widehat\E_{nr}}}  \ar[r]^-{\tilde\f_{N^*}}   & N^*\otimes_W \O_{\widehat\E_{nr}}  
}\]
\end{lem}


D'où, en rassemblant tout ceci, nous obtenons le théorème suivant :

\begin{theo}\label{th1}
Soit $0\leq h\leq \frac{p-2}{2}$, et notons $\widetilde{\mathcal{F}}_1$ le foncteur exact de la catégorie $\mathop{\bf MF_W^{\pm h}}$ vers la catégorie des $ \O_{\widehat\E_{nr}}$-modules libres de rang fini, défini par : si $N$ objet de $\mathop{\bf MF_W^{\pm h}}$, $\widetilde{\mathcal{F}}_1(N)=\mathop{\bf \widetilde  V_{cris} }(N) \otimes_{\Z_p} \O_{\widehat\E_{nr}}$. Alors il existe $\tilde\f$ un isomorphisme de foncteurs entre $\widetilde{\mathcal{F}}_1$ et $\mathcal{F}_2$, préservant le dual. Nous pouvons supposer de plus :
\begin{itemize}
\item pour tous objet $N$ et $N'$ de $\Mf$, tel que $N\otimes N'$ soit encore un objet de $\Mf$, $\tilde\f_{N\otimes N'}=\tilde\f_{N}\otimes \tilde\f_{N'}$ ;
\item pour tout uplet d'objets $(N_{i,j})_{1\leq j\leq n,\,1 \leq i\leq n_j}$ de $\Mf$, pour tout sous-$\Phi$-module filtré $L$ facteur direct (comme $W$-module) de $\displaystyle\bigoplus_{j=1}^n\otimes_{i=1}^{n_j} N_{i,j}$, l'application $\bigoplus \otimes \tilde\f_{N_{i,j}}$ envoie $(\mathop{\bf V_{cris,p}}(D_L)\cap \displaystyle\bigoplus_{j=1}^n\otimes_{i=1}^{n_j}  \mathop{\bf \widetilde V_{cris}}(N_{i,j}))\otimes_{\Z_p}\O_{\widehat\E_{nr}}$ bijectivement sur $L\otimes_W \O_{\widehat\E_{nr}}$.
\end{itemize}
\end{theo}

\section{Position des réseaux}

Pour tout ce paragraphe, nous supposerons donnés $\r\,:\,\Gamma_{\K}\to G(\Z_p)$ une représentation cristalline à valeurs dans les point sur $\Z_p$ d'un groupe algébrique lisse sur $\Z_p$, $G$, et $U$ un $\Z_p$-module libre de rang $n$, avec $\a : G\to GL_U$ une immersion fermée.

\subsection{Description des groupes plats sur $\Z_p$}\label{para}
Notons $U_W=U\otimes_{\Z_p}W$. Identifions $G$ avec son image dans $GL_U$. Nous allons donner une définition plus exploitable de $G_W=G\times_{\Z_p}W$ dans un cas particulier :

Prenons une base de $U$ et supposons que l'immersion de $G$ dans $GL_U$ induise une immersion dans $\End_U$ (c'est-à-dire $G=\Spec(\Z_p[X_{i,j}]_{1\leq i,j\leq n}/I)$).

Notons $\Z_p[X_{i,j}]_{\leq d}$ les polyn\^omes de degré total inférieur ou égal à $d$. Le groupe $GL_{U}$ agit sur $\Z_p[X_{i,j}]$ par : pour toute $\Z_p$-algèbre $R$, si $f\in R[X_{i,j}]$, et $s\in GL_{U}(R)$, alors $\eta_s(f)$ est le polyn\^ome défini par $\eta_s(f)(y)=f(s^{-1}y)$. Cette action est linéaire et laisse stable $ \Z_p[X_{i,j}]_{\leq d}$. 

\begin{prop}\label{rr}
Soit $G=\Spec(\Z_p[X_{i,j}]_{1\leq i,j\leq n}/I)$ un groupe algébrique plat sur $\Z_p$, soient $(f_1,\cdots ,f_r)$ des générateurs de l'idéal $I$, et si $d$ est le maximum des degrés totaux des $f_i$, posons $E=I\cap \Z_p[X_{i,j}]_{\leq d}$. Alors $E$ est facteur direct, et pour toute $\Z_p$-algèbre $R$, si $E_R=E\otimes_{\Z_p}R$, $G(R)=\{g\in GL_U(R)|\eta_g(E_R)=E_R\}$.
\end{prop}

\begin{rmq}
Si $\a: G\to GL_U$ n'induit pas une immersion fermée de $G$ dans $\End_U$, il suffit de composer $\a$ avec une immersion fermée $\beta : GL_U\to GL_{U'}$ tel que $\beta$ induise une immersion fermée de $GL_U$ dans $\End_{U'}$. Par exemple, $U'=U\oplus \Z_p$ avec $\beta=\Id\oplus \frac{1}{\mathop{\rm det}}$, ou bien $U'=U\oplus U^*$ (où $U^*$ est le dual de $U$) avec $\beta(g)=(g,{}^tg^{-1})$. Par contre, le $E$ donné par la proposition \ref{rr} dépendra du morphisme $\beta$ considéré.
\end{rmq}

\begin{proof} Commençons par le lemme suivant :

\begin{lem}
Pour toute $\Z_p$-algèbre $R$, le module $I\otimes_{\Z_p}R$ s'injecte dans $R[X_{i,j}]$.
\end{lem}

\begin{proof} Pour tout $k$, notons $E_k=I\cap \Z_p[X_{i,j}]_{\leq k}$. Le module $\Z_p[X_{i,j}]_{\leq k}$ est un $\Z_p$-module libre de type fini, et $\Z_p[X_{i,j}]_{\leq k}/E_k$ est sans $p$-torsion car il s'injecte dans $\Z_p[X_{i,j}]_{1\leq i,j\leq n}/I$ qui est plat sur $\Z_p$ (par hypothèse). Donc le module $E_k$ est un $\Z_p$-module libre facteur direct dans $\Z_p[X_{i,j}]_{\leq k}$. Donc $E_k$ est un facteur direct de $\Z_p[X_{i,j}]$, par conséquent $E_k\otimes_{\Z_p}R$ s'injecte dans $R[X_{i,j}]$. Or, $I$ est la réunion des $E_k$, donc $I\otimes_{\Z_p} R$ s'injecte dans $R[X_{i,j}]$.
\end{proof}

Au passage, nous avons démontré une assertion de la proposition, à savoir que $E:=E_d$ est facteur direct.

Notons $H$ le groupe algébrique défini par $H(R)=\{g\in GL_U(R)|\eta_g(E_R)=E_R\}$. L'application naturelle $H\to GL_U$ est une immersion fermée et $H(R)=\{g\in GL_U(R)|\eta_g(E_R)\subset E_R\}$ car $E$ est facteur direct. Nous voulons montrer que $H=G$. Si $s\in GL_U(R)$ vérifie $\eta_g(E_R)=E_R$, alors pour tout $i$, $\eta(s)(f_i)\in E_R\subset I\otimes_{\Z_p}R$, donc $\eta(s)(f_i)(Id)=0$ car $Id\in G(R)$. Donc, par définition de l'action, $f_i(s^{-1})=0$ pour tout $i$, donc la famille $(f_i)$ étant une famille de générateurs de l'idéal $I$ sur $\Z_p$, nous obtenons $s^{-1}\in G(R)$, or $G$ est un groupe, donc $s\in G(R)$, ce qui montre l'inclusion $H(R)\subset G(R)$. Donc il existe un monomorphisme $H\to G$.

Montrons que c'est une immersion fermée : le morphisme $H\to GL_U$ est une immersion fermée, et $\a$ est par hypothèse une immersion fermée de $G$ dans $GL_U$. Donc, en notant $A[K]$ l'algèbre affine d'un groupe $K$, les flèches $A[GL_U]\to A[G]$ et $A[GL_U]\to A[H]$ sont surjectives, et nous avons le diagramme commutatif suivant :

\[\xymatrix{
A[G]\ar[d]&A[GL_U]\ar@{->>}[l]\ar@{->>}[dl]\\
A[H]
}\]
par conséquent, la flèche $A[G]\to A[H]$ est surjective, donc $H\to G$ est bien une immersion fermée. 

Nous allons maintenant montrer que $G$ et $H$ ont même fibre générique. Pour cela, donnons une description de $G$ semblable à celle de $H$ :

\begin{lem}
Pour toute $\Z_p$-algèbre $R$, nous avons $$G(R)=\{g\in GL_U(R)|\eta_g\big((I\otimes_{\Z_p}R)\cap R[X_{i,j}]_{\leq d}\big)=(I\otimes_{\Z_p}R)\cap R[X_{i,j}]_{\leq d}\}$$
\end{lem}

\begin{proof}
Fixons $R$, et posons $M=(I\otimes_{\Z_p}R)\cap R[X_{i,j}]_{\leq d}$. Si $s\in GL_U(R)$ vérifie $\eta_g(M)=M$, alors pour tout $i$, $\eta(s)(f_i)\in M\subset I\otimes_{\Z_p}R$, donc $\eta(s)(f_i)(Id)=0$ car $Id\in G(R)$. Donc, par définition de l'action, $f_i(s^{-1})=0$ pour tout $i$, donc la famille $(f_i)$ étant une famille de générateurs de l'idéal $I$ sur $\Z_p$, nous obtenons $s^{-1}\in G(R)$, or $G$ est un groupe, donc $s\in G(R)$, ce qui montre une inclusion.

L'inclusion réciproque sera un corollaire du lemme de Yoneda : d'abord, il suffit de montrer que $\eta_s(M)\subset M$, car $\eta$ est une action de groupe. Puis, si $s\in G(R)$ et $f\in M$, notons $P=\eta(s)(f)$. Pour toute $R$-algèbre $B$, pour tout $g\in G(B)$, nous avons $P(g)=f(s^{-1}g)=0$ car $s^{-1}g\in G(B)$ et $f\in I\otimes_{\Z_p}R$, donc $P$ est dans $I\otimes_{\Z_p}R$ (et de bon degré, donc $P\in M$) par la remarque suivante :

Soit $G=\Spec(A/I)$ un groupe algébrique au dessus d'un anneau $R$, alors si $J=\{f\in A |\text{ pour toute $R$-algebre }B,\, \forall g\in G(B), \, f(g)=0\}$, nous avons $I=J$. En effet, si $K=\Spec(A/J)$, alors $K(B)\subset G(B)$ car $I\subset J$, puis la définition de $J$ nous dit que pour tout $\phi : A/I\to B$, $J$ est inclus dans $\Ker(\phi)$, d'où se factorise en $\phi : A/J\to B$, $J$, donc $G(B)\subset K(B)$. Le lemme de Yoneda donne alors $G=K$, donc $I=J$.
\end{proof}

\begin{lem}
Soit $S/R$ une extension d'anneau, supposons que $S$ est plat sur $R$. Alors, $((I\otimes_{\Z_p}R) \cap R[X_{i,j}]_{\leq d})\otimes_R S=(I\otimes_{\Z_p}S )\cap S[X_{i,j}]_{\leq d}$.
\end{lem}

\begin{proof}
En effet, notons $M_1=I\otimes_{\Z_p}R$, $M_2= R[X_{i,j}]_{\leq d}$ et $M_3= R[X_{i,j}]$, alors nous voulons voir que $(M_1\cap M_2)\otimes_R S=(M_1\otimes_RS)\cap (M_2\otimes_RS)$. Or, nous avons la suite exacte courte de $R$-modules \[\xymatrix{
0 \ar[r] & M_1\cap M_2\ar[r] & M_3\ar[r]& M_3/M_1\oplus M_3/M_2
}\]
et nous avons supposé que $S$ est plat sur $R$, donc en tensorisant par $S$ au dessus de $R$, nous obtenons que 
\[\xymatrix{
0 \ar[r] & (M_1\cap M_2)\otimes S \ar[r]& M_3\otimes S\ar[r]& (M_3/M_1)\otimes S\oplus (M_3/M_2)\otimes S 
}\]
est une suite exacte, ce qui conclut car $(M_3/M_i)\otimes_RS=M_3\otimes_R S/M_i\otimes_RS$.
\end{proof}

Puis, $G$ et $H$ ont même fibre générique, par application directe des deux lemmes précédents. 

Il ne reste plus qu'à voir que si $H\to G$ est une immersion fermée telle que $G$ et $H$ ont même fibre générique, et $G$ plat, alors $H=G$. Nous voulons montrer que la flèche surjective $A[G]\to A[H]$ est aussi injective. $G$ étant plat, l'application naturelle $i_G : A[G] \to A[G]\otimes_{\Z_p}\Q_p$ est injective. $H$ et $G$ ayant même fibre générique, la flèche $f:A[G]\to A[H]$ donnée par l'immersion fermée induit une bijection $f_p : A[G]\otimes_{\Z_p}\Q_p \to  A[H]\otimes_{\Z_p}\Q_p$. De plus le diagramme suivant est commutatif :
\[\xymatrix{A[G]\ar[r]^-{f}\ar@{^{(}->}[d]^-{i_G}& A[H]\ar[d]^-{i_H}\\
A[G]\otimes_{\Z_p}\Q_p\ar@{^{(}->>}[r]^-{f_p}&A[H]\otimes_{\Z_p}\Q_p
}\]
donc $f$ est bien injective, donc $H=G$.
\end{proof}

En fait, lors de l'identification de $GL_U$ avec $GL_{n,\Z_p}$, nous avons identifié la représentation de $GL_U$ que sont les polyn\^omes homogènes de degré $i$ (noté $\Z_p^i[X_{i,j}]$) avec $\mathop{\rm Sym}^i(\End_U)$ qui est un sous-objet de $\End_U^{\otimes i}$, ou bien, dit autrement, $  \Z_p[X_{i,j}]_{\leq d}$ a été identifié à un sous-objet de $\oplus_{0\leq i\leq d} \End_U^{\otimes i}$. L'action de $GL_U$ sur $\Z_p^i[X_{i,j}]$ est le produit tensoriel de l'action naturelle de $GL_U$ sur $U^*$ par l'action triviale de $GL_U$ sur $U$ dans $\End_U^{\otimes i}=(U\otimes_{\Z_p}U^*)^{\otimes i}$ ; par conséquent $E$ est un sous-module de $\displaystyle\bigoplus_{0\leq i\leq d}(\underbrace{U^*\oplus\cdots \oplus U^*}_{n\text{ fois}})^{\otimes i}$ (où $n=\rg(U)$).

Appliquons ceci non pas au plongement $\a$, mais à $\a^* : G\to GL_{U^*}$, défini par $\a^*(g)={}^t\a(g)^{-1}$. Il existe alors $E^*=\displaystyle\bigoplus_{0\leq i\leq k}(\underbrace{U\oplus\cdots \oplus U}_{n\text{ fois}})^{\otimes i}\cap I^*$ un sous $\Z_p$-module (libre facteur direct) tel que $s\in \a(G)(R)\Leftrightarrow s(E^*_R)=E^*_R$ (où $I^*$ est l'idéal définissant $\a^*(G)$). 

Rassemblons tout ceci dans le lemme suivant :

\begin{prop}\label{pppxx}
Soit $G$ un groupe plat sur $\Z_p$, $U$ un $\Z_p$-module libre de rang $n$ et $\a : G\to GL_U$ une représentation qui induit une immersion fermée dans $\End_U$. Identifions $G$ avec son image. Alors, il existe un entier $k$ et un sous $\Z_p$-module $E$ (facteur direct) de $\displaystyle\bigoplus_{0\leq i\leq k}(\underbrace{U\oplus\cdots \oplus U}_{n\text{ fois}})^{\otimes i}$ laissé stable par l'action naturelle de $G(\Z_p)$ (provenant de celle de $GL_U$, notée $\eta$) tels que $G(R)=\{g\in GL(R)|\eta_g(E_R)=E_R\}$ pour toute $\Z_p$-algèbre $R$.
\end{prop}

Nous pouvons alors définir sur $M=\mathop{\bf  D_{cris} }(U)$ (ou sur $\mathop{\bf \widetilde{D}_{cris}}(U)$ suivant les cas) un groupe algébrique sur $W$, $G_M$, par : si $\eta$ est l'action naturelle de $GL_M$ sur $\displaystyle\bigoplus_{0\leq i\leq k}(\underbrace{M\oplus\cdots \oplus M}_{n\text{ fois}})^{\otimes i}$, alors $G_M(R)=\{g\in GL_M(R)|\eta_g(\overline E_R)=\overline E_R\}$ pour toute $W$-algèbre $R$. Il ne reste qu'à bien choisir $\overline E$ en liaison avec $E$, ce qui nous conduit au théorème suivant : 

\begin{theo}\label{thjj2}
Supposons $G$ lisse sur $\Z_p$. Si $\a : G\to GL_U$ induit une immersion fermée dans $\End_U$ et si la représentation de $\Gamma_{\K}$ induite sur $U$ par $\a$ (et par $\r : \Gamma_{\K}\to G(\Z_p)$) vérifie 
\begin{itemize}
\item soit elle est à poids de Hodge-Tate dans $[\![0,h]\!]$ avec $0\leq h\leq p-2$
\item soit elle est à poids de Hodge-Tate  dans $[\![-h,h]\!]$ avec $0\leq h\leq \frac{p-2}{2}$
\end{itemize}
alors, en prenant $$\overline E=\mathop{\bf  D_{cris,p} }(E\otimes_{\Z_p}\Q_p)\cap\displaystyle\bigoplus_{0\leq i\leq k}(\underbrace{\D(U)\oplus\cdots \oplus \D(U)}_{n\text{ fois}})^{\otimes i} $$ dans le premier cas, ou  $$\overline E=\mathop{\bf  D_{cris,p} }(E\otimes_{\Z_p}\Q_p)\cap\displaystyle\bigoplus_{0\leq i\leq k}(\underbrace{\mathop{\bf \widetilde{D}_{cris}}(U)\oplus\cdots \oplus {\mathop{\bf \widetilde{D}_{cris}}(U)}}_{n\text{ fois}})^{\otimes i}$$ dans le deuxième cas, il existe une bijection $\Psi\, :\, U\otimes_{\Z_p}W\to M$ qui identifie $G\times_{\Z_p}W$ et $G_M$. 
\end{theo}

\begin{rmq}
Pour qu'une application $\Psi$ bijective identifie $G\times_{\Z_p}W$ et $G_M$, il suffit de montrer que l'application naturelle induite par $\Psi$ envoie $E\otimes_{\Z_p}W$ bijectivement sur $\overline E$.
\end{rmq}

\subsection{Demonstration du théorème \ref{thjj2}}\label{dem}
Soit $h\,:\,U_R=U\otimes_{\Z_p}R\simeq M_R=M\otimes_W R$ un isomorphisme de $R$-modules, alors $h$ induit $s(h)$ de $\displaystyle\bigoplus_{0\leq i\leq k}(\underbrace{U_R\oplus\cdots \oplus U_R}_{n\text{ fois}})^{\otimes i}$ sur $\displaystyle\bigoplus_{0\leq i\leq k}(\underbrace{M_R\oplus\cdots \oplus M_R}_{n\text{ fois}})^{\otimes i}$, qui est aussi un isomorphisme.

Considérons alors $$\mathop{\bf Isom}(R)=\{h \,:\,U\otimes_{\Z_p}R\simeq M\otimes_W R| s(h)(E_R)=\overline E_R\}$$ pour $R$ une $W$-algèbre. C'est un sous-$W$-schéma de $\mathop{\bf Isom_W}(U\otimes_{\Z_p}W, M)$ (les $W$-isomorphismes de $U\otimes_{\Z_p}W$ sur $M$).

Il est non vide, car $\f_{\mathop{\bf {D}_{cris}}(U)}$ (ou $\tilde\f_{\mathop{\bf {\widetilde D}_{cris}}(U)}$ suivant les conditions sur les poids de Hodge-Tate) induit un élément de $\mathop{\bf Isom}(\O_{\widehat\E_{nr}})$. C'est une retraduction du théorème \ref{theoxxx} (ou du théorème \ref{th1})

\begin{lem}
Le schéma $\mathop{\bf Isom}\times_W \O_{\widehat\E_{nr}}$ est un torseur trivial sous $G\times_W \O_{\widehat\E_{nr}}$. 
\end{lem}

\begin{proof} $G$ agit naturellement et fidèlement à gauche sur $\mathop{\bf Isom}$ : si $f\in \mathop{\bf Isom}(R)$ et $g\in G(R)$,  \[\xymatrix{  { U \otimes_{\Z_p} R} \ar[r]^{g} & { U \otimes_{\Z_p} R} \ar[r]^{f} &  {M\otimes_W R} }\] 
est bien un isomorphisme.

Puis, l'application naturelle 
\[\xymatrix{  { \displaystyle\bigoplus_{0\leq i\leq k}(\underbrace{U\oplus\cdots \oplus U}_{n\text{ fois}})^{\otimes i} \otimes_{\Z_p} R} \ar[r]^{\eta_g} & { \displaystyle\bigoplus_{0\leq i\leq k}(\underbrace{U\oplus\cdots \oplus U}_{n\text{ fois}})^{\otimes i} \otimes_{\Z_p} R} \ar[d]^{s(f)}\\ &  {\displaystyle\bigoplus_{0\leq i\leq k}(\underbrace{M\oplus\cdots \oplus M}_{n\text{ fois}})^{\otimes i}\otimes_W R} }\] 
s'identifie naturellement à $s(f\circ g)$.

Enfin, la définition de $\mathop{\bf Isom}$, la proposition \ref{pppxx} et le fait que $s(f\circ g)=s(f)\circ \eta_g$ donnent bien $s(f\circ g)(E_R)=E_R$, donc que $f\circ g\in \mathop{\bf Isom}(R)$. Le groupe $G$ agit donc sur $\mathop{\bf Isom}$ par $(g,f)\mapsto f\circ g^{-1}$.

La fidélité provient de ce qu'un élément de $\mathop{\bf Isom}$ est un isomorphisme de modules.

Pour finir, il reste à montrer que pour $f,f'\in \mathop{\bf Isom}(R)$, il existe $g\in G(R)$ avec $f'=f\circ g^{-1}$. Autrement dit, il faut voir que $g=f'^{-1}\circ f$ est bien un élément de $G(R)$. Cela se montre de la même façon que précédement. C'est la propriété \ref{pppxx} qui est le point essentiel.\end{proof}

Nous venons donc de montrer que $\mathop{\bf Isom}\times_W \O_{\widehat\E_{nr}}$ est un $G\times_W \O_{\widehat\E_{nr}}$-espace homogène ayant un point sur $\O_{\widehat\E_{nr}}$, or $G\times_W \O_{\widehat\E_{nr}}$ est un groupe lisse, donc $\mathop{\bf Isom}\times_W \O_{\widehat\E_{nr}}$ est lisse, donc $\mathop{\bf Isom}$ aussi (car la flèche $\Spec( \O_{\widehat\E_{nr}})\to \Spec(W)$ est fidèlement plate et quasi-compact, donc c'est une application directe du corollaire 17.7.3 de EGA IV).

$\mathop{\bf Isom}$ est lisse, donc par le lemme de Hensel (cf. théorème 18.5.11 b de EGA IV), $\mathop{\bf Isom}(W)$ se surjecte (par la réduction modulo $p$) sur $\mathop{\bf Isom}(k)$. Si ce dernier est non vide, nous aurons bien montré que $\mathop{\bf Isom}(W)$ est non vide, ce qui prouvera le théorème.

$\mathop{\bf Isom}\times_W \O_{\widehat\E_{nr}}$ est lisse, donc, toujours par le lemme de Hensel, $\mathop{\bf Isom}(\O_{\widehat\E_{nr}})$ se surjecte sur $\mathop{\bf Isom}(\O_{\widehat\E_{nr}}/p)$, donc le $k$-schéma $\mathop{\bf Isom}\times \Spec(k)$ est non vide (car $\O_{\widehat\E_{nr}}/p$ est une $k$-algèbre), donc par le théorème des zéros de Hilbert, $\mathop{\bf Isom}\times \Spec(k)(k)$ est non vide (car $k$ est algébriquement clos).\qed

Remarquons qu'en nous donnant un $\Z_p$-module $N\subset M$ qui engendre $M$ comme $W$-module (c'est à dire $N\otimes_{\Z_p}W=M$), et tel que le $\Z_p$-module $N'=\overline E \cap\displaystyle\bigoplus_{0\leq i\leq k}(\underbrace{N\oplus\cdots \oplus N}_{n\text{ fois}})^{\otimes i} $ engendre $\overline E$ comme $W$-module (par exemple, avec les notations du paragraphe \ref{modu}, $N=M^{f_M}$ (et alors $N'=\overline E^{f_{\overline E}}$) ou $N=M_{\Z_p}^{u_N^{-1}f_M}$ (et alors $N'=\overline E_{\Z_p}^{u_{\overline E}^{-1}f_{\overline E}}$)), nous pouvons définir $\mathop{\bf Isom}$ sur $\Z_p$ par 
$$\mathop{\bf Isom}(R)=\{h \,:\,U\otimes_{\Z_p}R\simeq N\otimes_{\Z_p} R| s(h)(E_R)=N'_R\}$$ pour $R$ une $\Z_p$-algèbre. Alors, par un théorème de Lang (tout $H$-torseur défini sur $\F_p$ est trivial si $H$ est un groupe algébrique sur $\F_p$ connexe) , en supposant que $G$ est à fibre spéciale connexe, nous avons $\mathop{\bf Isom}(\F_p)$ non vide, donc par lissité, $\mathop{\bf Isom}(\Z_p)$ est non vide, et donc sous cette hypothèse, nous pouvons supposer que $\Psi(U)=N$. De plus, $\Psi$ identifie $G$ à une forme sur $\Z_p$ de $G_M$ (celle qui est définie à l'aide de $N'$).

\begin{cor}
Sous les hypothèses et notations précédentes, si $U'$ est un réseau de $U\otimes_{\Z_p}\Q_p$ laissé stable par l'action de $G$, alors $\Psi[\frac{1}{p}]=\Psi\otimes_{W}\K$ envoie $U'\otimes_{\Z_p}W$ sur $\D(U')$ si les poids de Hodge-Tate sont positifs, ou sur $\mathop{\bf \widetilde D_{cris}}(U')$ sinon.  
\end{cor}

\begin{proof}
Notons $M=\D(U)$ (ou $M=\mathop{\bf \widetilde D_{cris}}(U)$ si les poids de Hodge-Tate ne sont pas tous positifs). Tout d'abord, quitte à multiplier par une certaine puissance de $p$, nous pouvons supposer $U'\subset U$. Puis, $\Psi\in \mathop{\bf Isom}(W)\subset \mathop{\bf Isom}( \O_{\widehat\E_{nr}})$ et $\f_M\in \mathop{\bf Isom}( \O_{\widehat\E_{nr}})$ (respectivement, $\tilde\f_M\in \mathop{\bf Isom}( \O_{\widehat\E_{nr}})$), donc, comme $\mathop{\bf Isom}$ est un $G\times_{\Z_p}W$-espace homogène, il existe $g\in G( \O_{\widehat\E_{nr}})$ tel que $\Psi=\f_M\circ g$ (respectivement $\Psi=\tilde\f_M\circ g$). Il suffit donc (puisque $U'$ est stable par $G$) de vérifier la propriété pour $\f_M$ (ou $\tilde\f_M$), or ceci provient juste de la fonctorialité de $\f_M$ (et de $\tilde\f_M$) pour les sous-objets.
\end{proof}

\subsection{Exemples}
Remarquons que $GL_U$ se plonge naturellement par une immersion fermée $\beta$ dans $GL_{U\oplus U^*}$ où l'action sur $U$ est l'action naturelle, et l'action sur le dual $U^*$ est donnée par la transposée de l'inverse. De plus, $\beta$ provient d'une immersion fermée de $GL_U$ dans $\End_{U\oplus U^*}$, car l'image de $\beta$ est un sous-groupe fermé de $SL_{U\oplus U^*}$. Donc, si la représentation galoisienne $U$ est à poids de Hodge-Tate dans $[\![0,h]\!]$ avec $h\leq \frac{p-2}{2}$, ou dans $[\![-h,h]\!]$ avec $h\leq \frac{p-2}{4}$, alors l'immersion $\a'=\beta\circ \a$ vérifie en partie les hypothèses du théorème \ref{thjj2}. Soit alors $E$ le sous-module définissant $G$ dans $U\oplus U^*$ (de la manière décrite dans le paragraphe \ref{para}). Nous définissons de même sur $\mathop{\bf \widetilde{D}_{cris}}(U)\oplus \mathop{\bf \widetilde{D}_{cris}}(U)^*$ un groupe $G_M$ à l'aide de $$\overline{E}=\mathop{\bf  D_{cris,p} }(E\otimes_{\Z_p}\Q_p)\cap\displaystyle\bigoplus_{0\leq i\leq k}(\underbrace{\mathop{\bf \widetilde{D}_{cris}}(U)\oplus\cdots \oplus {\mathop{\bf \widetilde{D}_{cris}}(U)}}_{n\text{ fois}})^{\otimes i}$$
 
\begin{theo}\label{thjj}
Si la représentation de $\Gamma_{\K}$ induite sur $U$ par $\a$ (et par $\r$) vérifie 
\begin{itemize}
\item soit elle est à poids de Hodge-Tate dans $[\![0,h]\!]$ avec $0\leq h\leq \frac{p-2}{2}$
\item soit elle est à poids de Hodge-Tate dans $[\![-h,h]\!]$ avec $0\leq h\leq \frac{p-2}{4}$

\end{itemize}
alors, avec les notations précédentes, il existe un isomorphisme de $W$-module $\Psi\, :\, U\otimes_{\Z_p}W\to M$ qui induit une bijection $\Psi\oplus {}^t\Psi^{-1}\, :\, (U\oplus U^*)\otimes_{\Z_p}W\to M\oplus M^*$ identifiant $G\times_{\Z_p}W$ et $G_M$. En particulier, le groupe $G_M$ (plongé dans $GL_{M\oplus M^*}$) laisse stable $M$ et $M^*$.
\end{theo}

\begin{proof}
Considérons le schéma $\mathop{\bf Isom}(R)=\{h \,:\,U\otimes_{\Z_p}R\simeq M\otimes_W R| s(h)(E_R)=\overline E_R\}$ (ce n'est à priori pas le même que celui considéré dans la démonstration du théorème \ref{thjj2}, qui considère des morphismes définis sur $U\otimes U^*$). C'est un $G\times_{\Z_p}W$ espace homogène, qui a un point sur $W$ (cela se montre de la même façon que lors de la démonstration du théorème \ref{thjj2}).\end{proof}

\subsection{Données initiales pour un $\Phi$-module filtré}
Formulons ici comment les idées introduites précédemment se traduisent dans le formalisme introduit par Rapoport et Zink (cf. \cite{RZ}). Soit $G$ un groupe algébrique lisse sur $\Z_p$ et $\r$ un morphisme de groupe de $\Gamma_{\K}$ dans $G(\Z_p)$, $\mu : \G_m\to G$ un cocaractère défini sur $W$, et $b\in G(W)$. Alors, à toute représentation $\beta : G\to GL_{U}$ où $U$ est un $\Z_p$-module libre de rang fini, nous pouvons associer un objet $\I(U)$ de $\mathop{\bf MF_W'}$, défini par :
\begin{itemize}
\item le $W$-module sous-jacent est $U\otimes_{\Z_p}W=M$ ;
\item $\Fil^i(M)=\displaystyle\bigoplus_{j\geq i} M_j$ où $M_j$ est l'espace propre de poids $j$ correspondant à $\mu$ ;
\item $\p^i=p^{-i}b\circ(Id\otimes \s)$ sur $\Fil^i(M)$, c'est-à-dire si $v=\sum\limits_k v_k\otimes x_k\in \Fil^i(M)$ avec $v_k\in U$ et $x_k\in W$, $\p^i(v)=\sum\limits_k \s(x_k)\beta(b)(v_i)$.
\end{itemize}

\begin{defi}
Le triplet $(\mu,b,\beta)$ est dit admissible si $\I(U)$ est un objet de $\mathop{\bf MF_{W,tf}}$.
\end{defi}

Si $G$ est supposé lisse, si $\beta$ induit une immersion fermée de $G$ dans $\End_U$, et $\I(U)$ est un objet de $\MF$ avec $0\leq h\leq p-2$, nous pouvons faire de même qu'au paragraphe \ref{para} : $G$ est défini par $E$ un $\Z_p$-module bien choisi de $\displaystyle\bigoplus_{0\leq i\leq k}(\underbrace{U\oplus\cdots \oplus U}_{n\text{ fois}})^{\otimes i}$ (alors $E\otimes_{\Z_p}W$ sera un objet de $\mathop{\bf MF_{W,tf}}$), et sur $\V(\I(U))$ nous construisons $G_{\V(\I(U))}$ sur $\Z_p$ par son foncteur des points (si $R$ est une $\Z_p$-algèbre, $ G_{\V(\I(U))}(R)$ est le sous-groupe de $GL_{\V(\I(U))}(R)$ formée des éléments qui laissent stable $$\big(\mathop{\bf V_{cris,p}}(E\otimes_{\Z_p}\K)\cap \displaystyle\bigoplus_{0\leq i\leq k}(\underbrace{\V(M)\oplus\cdots \oplus \V(M)}_{n\text{ fois}})^{\otimes i}\big)\otimes_{\Z_p}R$$ 

\begin{theo}
Sous les conditions précédentes, avec
\begin{itemize}
\item soit $M$ est un objet de $\MF$ avec $0\leq h\leq p-2$ et $\beta$ induit une immersion fermée de $G$ dans $\End_U$,
\item soit $M$ est un objet de $\MF$ avec $0\leq h\leq \frac{p-2}{2}$,
\item soit $M$ est un objet de $\Mf$ avec $0\leq h\leq \frac{p-2}{2}$ et $\beta$ induit une immersion fermée de $G$ dans $\End_U$,
\item soit $M$ est un objet de $\Mf$ avec $0\leq h\leq \frac{p-2}{4}$,
\end{itemize}
alors il existe une bijection $\Psi\, :\, M\to \V(\I(U))\otimes_{\Z_p}W$ qui identifie $G\times_{\Z_p}W$ et $G_{\V(\I(U))}\times_{\Z_p}W$. De plus, la représentation galoisienne associée à $\V(\I(U))$ est à valeurs dans $G_{\V(\I(U))}(\Z_p)$.
\end{theo}

\begin{proof} L'existence de la bijection se montre de la m\^eme façon que pour le théorème \ref{thjj2} : nous introduisons un $G$-espace homogène $\mathop{\bf Isom}$ défini sur $\Z_p$, nous montrons qu'il est lisse sur $\Z_p$, $\mathop{\bf Isom}(k)$ est non vide par le théorème des zéros de Hilbert, et le lemme de Hensel conclut. La définition même de $G_{\V(\I(U))}$ implique que la représentation galoisienne associée à $\V(\I(U))$ est à valeurs dans $G_{\V(\I(U))}(\Z_p)$.\end{proof}

\begin{rmq}
Le théorème de Lang à propos des torseurs définis sur les corps finis implique que si la fibre spéciale de $G$ est connexe, alors le torseur $\mathop{\bf Isom}$ est trivial sur $\F_p$. Autrement dit, il a un point sur $\F_p$, que nous pouvons relever à $\Z_p$ par le lemme de Hensel. Par conséquent, si la fibre spéciale de $G$ est connexe, l'isomorphisme $\Psi$ est défini sur $\Z_p$, c'est à dire qu'il induit une bijection $\Psi : U\to \V(\I(U))$.
\end{rmq}

\bibliographystyle{alpha}
\bibliography{biblio}

\Addresses

\end{document}